\documentclass[12pt]{amsart}
\usepackage{amsmath,amsfonts,amssymb,amsxtra,amscd,enumerate,amsthm}
\usepackage{color}
\usepackage{latexsym}
\usepackage{epsfig}
\usepackage{fullpage}

\newcommand\g{\gamma}

\newcommand\e{\epsilon}
\renewcommand\l{\lambda}
\newcommand\tl{\tilde \lambda}
\newcommand\les{\lesssim}
\newcommand\ges{\gtrsim}

\newcommand\tg{\tilde g}
\newcommand\tPsi{\tilde \Psi}

\newcommand\R{\mathbb{R}}
\newcommand\C{\mathbb{C}}
\newcommand\Z{\mathbb{Z}}
\newcommand\N{\mathbb{N}}
\renewcommand\S{\mathbb{S}}
\newcommand\K{\mathcal K}
\newcommand\FH{{\mathcal F}_H}
\newcommand\FtH{{\mathcal F}_{\tilde H}}
\newcommand\F{{\mathcal F}}
\newcommand{\calO}{\mathcal O}
\newcommand\W{\mathcal{W}}
\newcommand\X{\mathcal{X}}
\newcommand\Zu{Z}

\newcommand{\V}{\tilde V}
\newcommand{\tpsi}{\tilde \psi}

\newcommand\la{\langle}
\newcommand\ra{\rangle}
\newtheorem{t1}{Theorem}
\numberwithin{t1}{section} 
\newtheorem{l1}[t1]{Lemma}
\newtheorem{p1}[t1]{Proposition}
\newtheorem{c1}[t1]{Corollary}
\newtheorem{d1}[t1]{Definition}

\newtheorem{r1}{Remark}

\numberwithin{equation}{section}

\newcommand{\Ss}{S^\sharp}
\newcommand{\Ns}{N^\sharp}
\newcommand{\WSs}{W\!S^\sharp}
\newcommand{\WNs}{W\!N^\sharp}

\newcommand{\WS}{W\!S}
\newcommand{\WN}{W\!N}
\newcommand{\ldSs}{l^2 S^\sharp}
\newcommand{\ldNs}{l^2 N^\sharp}

\newcommand{\He}{H^1_e}
\newcommand{\dHe}{\dot{H}^1_e}
\newcommand{\Hde}{H^2_e}
\newcommand{\dHde}{\dot{H}^2_e}

\begin{document}
\title[]{Near soliton evolution for equivariant Schr\"odinger Maps in two spatial dimensions}

\author{Ioan Bejenaru}
\address{ Department of Mathematics \\
  University of Chicago}

\author{ Daniel Tataru}
\address {Department of Mathematics \\
  University of California, Berkeley}

\thanks{ The first author was
  partially supported by NSF grant DMS1001676. The second author was
  partially supported by NSF grant DMS0354539}

\begin{abstract}  
  We consider the Schr\"odinger Map equation in $2+1$ dimensions, with
  values into $\S^2$. This admits a lowest energy steady state $Q$,
  namely the stereographic projection, which extends to a two
  dimensional family of steady states by scaling and rotation. We
  prove that $Q$ is unstable in the energy space $\dot H^1$. However,
  in the process of proving this we also show that within the
  equivariant class $Q$ is stable in a stronger topology $X \subset
  \dot H^1$.
\end{abstract}

\maketitle

\section{Introduction}

In this article we consider the Schr\"odinger map equation in $\R^{2+1}$
with values into $\S^2$,
\begin{equation}
u_t = u \times \Delta u, \qquad u(0) = u_0
\label{SM}\end{equation}
This equation admits a conserved energy,
\[
E(u) = \frac12 \int_{\R^2} |\nabla u|^2 dx
\]
and is invariant with respect to the dimensionless scaling
\[
 u(t,x) \to u(\lambda^2 t, \lambda x).
\]
The energy is invariant with respect to the above scaling,
therefore the Schr\"odinger map equation in $\R^{2+1}$ is 
{\em energy critical}.  

Local solutions for regular large initial data have been constructed
in \cite{SuSuBa} and \cite{Ga}.  Low regularity small data
Schr\"odinger maps were studied in several works, see \cite{Be},
\cite{Be2}, \cite{bik}, \cite{IoKe2}, \cite{IoKe3}, \cite{Ka},
\cite{KaKo}, \cite{KeNa}, \cite{NaStUh}, \cite{NaStUh2},
\cite{NaShVeZe}.  The definitive result for the small data problem was
obtained by the authors and collaborators in \cite{BIKT}. There
global well-posedness and scattering are proved for initial data which
is small in the energy space $\dot H^1$.

However, such a result cannot hold for large data. In particular 
there exists a collection of families $\mathcal Q^m$  of 
finite energy stationary solutions, indexed by integers $m \geq 1$.
To describe these families we begin with the maps $Q^m$ 
defined in polar coordinates by
\[
Q^m(r,\theta) = e^{m \theta R} \bar Q^m(r), \qquad
\bar Q^m(r)=\left( \begin{array}{ccc} h_1^m(r) \\ 0 \\ h_3^m(r) \end{array} \right),
\qquad m \in \Z \setminus\{0\}
\]
with 
\[
h_1^m(r)=\frac{2r^m}{r^{2m}+1}, \qquad h_3^m(r)=\frac{r^{2m}-1}{r^{2m} + 1}. 
\]
Here $R$ is the generator of horizontal rotations, which 
can be interpreted as a matrix or, equivalently, as the operator below
\[
 R = \left( \begin{array}{ccc} 0 & -1 & 0 \\ 1 & 0 & 0 \\ 0 & 0 & 0 \end{array}
\right)  , \qquad R u  = \overrightarrow{k} \times u 
\]
The families $\mathcal Q^m$ are constructed from $Q^{ m}$
via the symmetries of the problem, namely scaling and isometries
of the base space $\R^2$ and of the target space $\S^2$. $Q^{-m}$ generates 
the same family $\mathcal Q^m$. 
The elements of $\mathcal Q^m$ are harmonic maps from $\R^2$ into
$\S^2$, and admit a variational characterization as the unique energy
minimizers, up to symmetries, among all maps $u: \R^2 \to \S^2$
within their homotopy class.  \medskip

In the above context, a natural question is to study Schr\"odinger
maps for which the initial data is close in $\dot H^1$ to one of the
$\mathcal Q^m$ families.  One may try to think of this as a small data
problem, but in some aspects it turns out to be closer to a large data
problem.  Studying this in full generality is very difficult.  In this
article we confine ourselves to a class of maps which have some extra
symmetry properties, namely the {\em equivariant} Schr\"odinger maps.
These are indexed by an integer $m$ called the equivariance class, and
consist of maps of the form
\begin{equation} \label{equiv}
u(r,\theta) = e^{m \theta R} \bar{u}(r)
\end{equation} 
In particular the maps $Q^m$ above are $m$-equivariant.
The case $m=0$ would correspond to spherical symmetry.
Restricted to equivariant functions the energy has the form
\begin{equation} \label{energy}
 E(u) = \pi \int_{0}^\infty \left( |\partial_r \bar{u}(r)|^2 +
\frac{m^2}{r^2}(\bar{u}_1^2(r)+\bar{u}_2^2(r)) \right) r dr
\end{equation}

Intersecting the full set $\mathcal Q^{m}$ with the $m$-equivariant
class and with the homotopy class of $Q^m$ we obtain the two parameter family $\mathcal Q^m_e$ generated from $Q^m$ by rotations and scaling,
\[
 \mathcal Q^m_e = \{ Q^m_{\alpha,\lambda}; \alpha \in \R/ 2\pi\Z, \lambda \in \R^+\}, 
\qquad  Q^m_{\alpha,\lambda}(r,\theta)= e^{\alpha R}  Q^m (\lambda r,\theta)
\]
Here $Q^m_{0,1} = Q^m$. Their energy depends on $m$ as follows:
\[
E(Q^m_{\alpha,\l}) = 4\pi m:= E(\mathcal Q^m)
\]

 \medskip

 The study of equivariant Schr\"odinger maps for $m$-equivariant
 initial data close to $\mathcal Q^m_e$ was initiated by Gustafson, Kang, Tsai
 in \cite{gkt1}, \cite{gkt}, and continued by Gustafson, Nakanishi,
 Tsai in \cite{gnt}. The energy conservation suffices to confine
 solutions to a small neighborhood of $\mathcal Q^m_e$ due to the
 inequality (see \cite{gkt1})
\begin{equation}
\text{dist}_{\dot H^1} (u,\mathcal Q^m_e)^2= \inf_{\alpha,\lambda}
\| Q^m_{\alpha,\lambda} - u\|_{\dot H^1}^2 \lesssim E(u) - E(\mathcal Q^m),
\label{stab}\end{equation}
which holds for all $m$-equivariant maps $u:\R^2 \to \S^2$ in the homotopy class of
$\mathcal Q^m_e$ with $0 \leq E(u) - E(\mathcal Q^m) \ll 1$.  One can
interpret this as an orbital stability result for $\mathcal Q^m_e$.
However, this does not say much about the global behavior of solutions
since these soliton families are noncompact; thus one might have even
finite time blow-up while staying close to a soliton family.

To track the evolution of an $m$-equivariant Schr\"odinger 
map $u(t)$ along $\mathcal Q^m_e$ we use  
functions $(\alpha(t),\lambda(t))$ describing trajectories 
in $\mathcal Q^m_e$. One may be tempted to try to 
choose them as minimizers for the infimum in \eqref{stab},
but this choice turns out not to be particularly helpful.
Instead, we will allow ourselves more freedom, and be content
with any choice $(\alpha(t),\lambda(t))$ satisfying
\begin{equation}
\| u - Q_{\alpha(t),\lambda(t)}^m\|_{\dot H^1}^2 \lesssim
E(u) - E(\mathcal Q^m) 
\label{goodal}\end{equation}
\medskip

An important preliminary step in this analysis is the next result
concerning both the local wellposedness in $\dot H^1$ and the
persistence of higher regularity:
\begin{t1}
  The equation \eqref{SM} is locally well-posed in $\dot H^1$ 
  for $m$-equivari\-ant initial data $u_0$ in the homotopy class of
  $\mathcal Q^m_e$ satisfying
\[
E(u) - E(\mathcal Q^m) \ll 1
\]
If, in addition, $u_0 \in \dot{H}^2$ then $u \in L^\infty_t \dot{H}^2$. 
Furthermore, the $\dot H^1 \cap \dot H^2$ regularity persists for as long as 
the function $\lambda(t)$ in \eqref{goodal} remains in a compact set.
\label{th2}\end{t1}

This follows from Theorem 1.1 in \cite{gkt1} and Theorem $1.4$ in \cite{gkt}.  Given the above
result, the main problem remains to understand whether the steady
states $Q_{\alpha,\lambda}^m$ are stable or not; in the latter case,
one would like to understand the dynamics of the motion of the
solutions move the soliton family.  The case of large $m$ was
considered in prior work:

\begin{t1}[\cite{gkt} for $m \geq 4$, \cite{gnt} for $m=3$]
The solitons $Q_{\alpha,\lambda}^m$ are stable in the $\dot H^1$
topology within the $m$-equivariant class.
\end{t1}

In this article we begin the study of the more difficult case $m=1$,
and establish a very different type of behavior.  The soliton
$Q^1$ plays a central role in our analysis, which is why we
introduce the notation $Q:=Q^1$. Since equivariant functions are
easily reduced to their one-dimensional companion via \eqref{equiv},
we introduce the one dimensional equivariant version of $\dot H^1$,
\begin{equation}\label{defhe}
\| f \|_{\dHe}^2 = \| \partial_r f\|_{L^2(rdr)}^2 + \| r^{-1} f \|_{L^2(rdr)}^2, \quad 
\| f \|_{\He}^2=\| f \|_{\dHe}^2 + \| f \|_{L^2(rdr)}^2
\end{equation}
This is natural since for functions $u: \R^2 \to \R^2$ with $u(r,\theta) = e^{\theta R}
\bar u(r)$ we have
\[
\| u \|_{\dot H^1}= (2\pi)^\frac12 \| \bar u \|_{\dHe},\qquad 
\| u \|_{H^1}= (2\pi)^\frac12 \| \bar u \|_{\He}.
\]
For our main result we introduce a slightly stronger 
topology $X$ with the property that 
\begin{equation}
\He \subset X \subset \dHe.
\label{Xemb}\end{equation}
This is defined in Section~\ref{spectral} in terms of the 
spectral resolution of the linearized evolution around the 
soliton.  In a nutshell, the $X$ norm penalizes the behavior 
near frequency zero. Our first result below asserts that the soliton
$Q$ is stable in the $X$ topology (which applies to $\bar Q$). 
\begin{t1}
 Let $m = 1$ and $\gamma \ll 1$. Then for each 
$1$-equivariant initial data $u_0$ satisfying 
\begin{equation}
 \| \bar{u}_0 - \bar Q\|_X \leq \gamma 
\label{tdata}\end{equation}
there exists a unique global solution $u$ so that
$\bar{u} - \bar Q \in C(\R;X)$ and 
\begin{equation}
 \| \bar u - \bar Q\|_{C(\R;X)} \lesssim \gamma
\label{tsolution}\end{equation}
Furthermore, this solution has a Lipschitz dependence
on the initial data in $X$, uniformly on compact time intervals.
\label{tmain-X}\end{t1}
We also refer the reader to Theorem~\ref{tmain-XS} for a more complete 
form of this theorem.
The above result holds true if $\bar Q$ is replaced by $\bar Q_{\alpha,\l}^1$,
which implies that the solitons $Q_{\alpha,\lambda}^1$ are stable in the $X$ topology.
However, our second result asserts that the solitons
$Q_{\alpha,\lambda}^1$ are unstable in the $\dot H^1$ topology:

\begin{t1}
 For each $\epsilon,\gamma  \ll 1$ and $(\alpha,\lambda)$ so that 
\begin{equation}\label{alinst}
 |\alpha|+|\lambda -1| \approx \gamma
\end{equation}
there exists a solution $u$ as in Theorem~\ref{tmain-X} with 
the additional property that 
\begin{equation}
 \|u(0) -  Q_{\alpha,\lambda}^1\|_{\dot H^1} \lesssim \epsilon \gamma
\label{tdatai}\end{equation}
while  (recall that $Q = Q^{1}_{0,1}$)
\begin{equation}
\lim\sup_{t \to \pm \infty}  \|u - Q\|_{ \dot H^1} \lesssim |\log \epsilon|^{-1} \gamma
\label{tsolutioni}\end{equation}
\label{tmain-H}\end{t1}

We remark that, in view of \eqref{Xemb}, the solutions $u$ 
in   Theorem~\ref{tmain-X} must satisfy 
\[
 E(u) - E(\mathcal Q^1) \lesssim \gamma^2
\]
and they can move at most $O(\gamma)$ along the soliton family in the
sense of \eqref{goodal}.  For the result in  Theorem~\ref{tmain-H} we consider a
more restrictive class of solutions, for which
\[
 \| \bar u - \bar Q\|_{X} \approx \gamma
\]
while staying closer to the soliton family, (see \eqref{stab}),
\begin{equation}
 E(u) - E(\mathcal Q^1) \approx  \epsilon^2 \gamma^2
\label{tdataia}\end{equation}
Thus by \eqref{goodal} the parameters $(\alpha(t),\lambda(t))$ are
restricted to an $O( \epsilon \gamma)$ range for each $t$.  On the
other hand, \eqref{tdatai} and \eqref{tsolutioni} show that for the
solution in Theorem~\ref{tmain-H} the parameters $(\alpha(t),\lambda(t))$ 
vary by about $O(\gamma)$ along the flow.

We also remark that if in addition the initial data $u_0$ is in $ \dot
H^2$ then by Theorem~\ref{th2} the solution $u$ remains in this space at all
times. While we do not prove a uniform in time $\dot H^2$ bound, such
an estimate seems nevertheless likely to hold for solutions as in
Theorem~\ref{tmain-X}.

To better frame the context of this paper, one should compare the
above results with results for the corresponding problem for the
Wave-Maps equation in $\R^{2+1}$ with values into $\S^2$. The
equivariant families of steady states $\mathcal Q^m_e$ are the same
there, and they are also orbitally stable. However, in the case of
Wave Maps all these steady states are unstable, and blow-up may occur
in finite time for all $m$. We refer the reader to the results in
\cite{KST}, \cite{RS}, and \cite{RR}. Of special relevance to the present paper
are some of the spectral techniques developed  in \cite{KST}; we 
further develop that circle of ideas in the present paper.

{\bf Acknowledgments:}
The authors are grateful to Alexandru Ionescu,  Carlos Kenig and Wilhelm Schlag
for many useful conversations concerning the Schr\"odinger maps dynamics.

\subsection{Definitions and notations.} 
\label{defnot}
We conclude this section with few definitions and notations. However,
the reader should be aware that many objects are defined as the paper
progresses; see Section \ref{seccoulomb} for all gauge elements and
their equations, Section \ref{spectral} for the Fourier analysis and
related objects/spaces and Sections \ref{linear}-\ref{lineart} for the
functions spaces used in the analysis of the nonlinear problem.

While at fixed time our maps into the sphere are functions defined on
$\R^2$, the equivariance condition allows us to reduce our analysis to
functions of a single variable $|x|=r \in [0,\infty)$.  One such
instance is exhibited in \eqref{equiv} where to each equivariant map
$u$ we naturally associate its radial component $\bar u$.  Some other
functions will turn out to be radial by definition, see, for instance,
all the gauge elements in Section \ref{seccoulomb}.  We 
agree to identify such radial functions with the corresponding
one dimensional functions of $r$.  Some of these functions 
are complex valued, and this convention allows us to use the bar 
notation with the standard meaning, i.e. the  complex conjugate.

Even though we work mainly with functions of a single spatial variable
$r$, they originate in two dimensions. Therefore, it is natural to
make the convention that for the one dimensional functions all the
Lebesgue integral and spaces are with respect to the $rdr$ measure,
unless otherwise specified. 

For the Sobolev spaces we have introduced $\dHe$ and $\He$ in
\eqref{defhe} as the natural substitute for $\dot{H}^1$ and $H^1$. In
a similar fashion we define $\dHde$ and $\Hde$ by the norms
\[
\| f \|^2_{\dHde}= \| \partial_r^2 f \|^2_{L^2} +\| r^{-1}\partial_r f \|_{L^2}^2+ \| r^{-2}f \|_{L^2}^2,
\qquad \| f \|^2_{\dHde}=\| f \|^2_{\dHde} + \| f \|_{L^2}^2
\] 
as the as the natural substitute for $\dot{H}^2$ and $H^2$.

For a real number $a$ we define $a^+=\max\{ 0, a \}$ and $a^-=\min\{0,a\}$.

We will use a dyadic partition of $\R^2$ (or $[0,\infty)$ after the
dimensional reduction) into sets  $\{A_m\}_{m \in \Z}$ given by
\[
A_m = \{  2^{m-1} < r < 2^{m+1}\}.
\]
We will also use the notation $A_{< k}=\cup_{m < k} A_m$ as well as
$A_{> k}, A_{\geq k}$ which are similarly defined.

Two operators which are often used on radial functions are $[\partial_r]^{-1}$
and $[r \partial_r]^{-1}$ defined as
\[
[\partial_r]^{-1} f(r) = - \int_{r}^\infty f(s) ds, \qquad 
[r \partial_r]^{-1}f(r) = - \int_{r}^\infty \frac{1}s f(s) ds
\]
A direct argument shows that 
\begin{equation} \label{rdrm}
\| [r\partial_r]^{-1}f \|_{L^p} \lesssim_p \| f \|_{L^p}, \qquad 1 \leq p < \infty
\end{equation}
We also have a weighted version 
\begin{equation} \label{rdrmw}
\| w [r\partial_r]^{-1}f \|_{L^p} \lesssim_p \| w f \|_{L^p}, \qquad 1 \leq p < \infty
\end{equation}
assuming that $g(r)=w(r) r^{\frac2{p}}$ is an increasing function satisfying
\[
g(r) \leq (1-\epsilon) g(2r)
\]
for some $\e > 0$. The proof is straightforward.

\section{An outline of the paper}

Due to the complexity of the paper, an overview of the ideas 
and the organization of the paper is necessary before an in-depth reading.

\subsection{ The frame method and the Coulomb gauge}
At first sight the Schr\"odinger Map equation has little to do with
the Schr\"odinger equation. A good way to bring in the Schr\"odinger
structure is by using the frame method. Precisely, at each point
$(x,t)\in \R^{2+1}$ one introduces an orthonormal frame $(v,w)$ in
$T_{u(x,t)} \S^2$. This frame is used to measure the derivatives of
$u$, and reexpress them as the complex valued radial {\em
  differentiated fields}
\[
\psi_1 =  \partial_r u \cdot v + i \partial_r u \cdot w, 
\qquad \psi_2 = \partial_\theta u \cdot v + i \partial_\theta u \cdot w.
\]
Here the use of polar coordinates is motivated by the equivariance
condition.  Thus instead of working with the equation for $u$, one
writes the evolution equations for the differentiated fields.
The frame $(v,w)$ does not appear directly there, but only via the
real valued radial connection coefficients
\[
A_1 = \partial_r v \cdot w, \qquad A_2 =  \partial_\theta v \cdot w,
\qquad A_0 =  \partial_t  v \cdot w.
\]

A-priori the frame is not uniquely determined. To fix it one first
asks that the frame be equivariant, and then that it satisfies an
appropriate condition.  Here it is convenient to use the {\em
  Coulomb gauge}; due to the equivariance this takes a
very simple form, $A_1 =0$.  The construction of the Coulomb gauge is
the first goal in the next section.  In Proposition~\ref{p:gauge} we
prove that for $\dot H^1$ equivariant maps into $\S^2$ close to
$Q$ there exists an unique Coulomb frame $(v,w)$ which
satisfies appropriate boundary conditions at infinity, see
\eqref{bcvw}. In addition, this frame has a $C^1$ dependence on the
map $u$. 

In the Coulomb gauge the other spatial connection
coefficient $A_2$, while nonzero, has a very simple form $A_2 = u_3$.
We will also compute $A_0$ in terms of $\psi_1$, $\psi_2$ and 
$A_2$,
\begin{equation}
 A_0 = - \frac12 \left( |\psi_1|^2 - \frac{1}{r^2}|\psi_2|^2\right)
+  [r \partial_r^{-1}] \left( |\psi_1|^2 - \frac{1}{r^2}|\psi_2|^2\right)
\label{intro-a0bis}\end{equation}

\subsection{The reduced field $\psi$}
Due to the equivariance the two fields $\psi_1$ and $\psi_2$ are
not independent. Hence it is convenient to work with a single field
\[
\psi = \psi_1 - i r^{-1}
\psi_2
\]
which we will call {\em the reduced field}. The relevance of the
variable $\psi$ comes from the following reinterpretation.  If $\W$ is
defined as the vector
\[
\W = \partial_r u - \frac{m}{r} u \times R u \in T_u(\S^2)
\] 
then $\psi$ is the representation of $W$ with respect to the frame $(v,w)$. 
On the other hand, a direct computation, see for instance \cite{gkt}, leads to
\[
E(u) = \pi \int_0^\infty \left( |\partial_r \bar{u}|^2 + \frac{m^2}{r^2} |\bar{u} \times
R \bar{u}|^2 \right) rdr = \pi \| \bar \W \|_{L^2(rdr)}^2 + 4\pi m
\]
where we recall that $u(r,\theta)= e^{m \theta R} \bar{u}(r)$.
Therefore $\psi=0$ is a complete characterization of $u$ being a harmonic map.
Moreover the mass of $\psi$ is directly related to the energy of $u$ via
\begin{equation} \label{basicpsi}
\|\psi\|_{L^2}^2 = \| \bar \W \|_{L^2}^2 = \frac{E(u)-4\pi m}{\pi}.
\end{equation}

A second goal of the next section is to derive an equation for
the time evolution of  $\psi$. This is governed by  a  cubic NLS type
equation, 
\begin{equation}
(i \partial_t + \Delta -\frac{2}{r^2} ) \psi =
\left( A_0 - 2 \frac{A_2}{r^2}-\frac{1}r \Im({\psi}_2 \bar{\psi})\right)\psi.
\label{intro-psieq}\end{equation}
In addition, we show that $\psi$ is connected back to
$(\psi_2,A_2)$  via the ODE system
\begin{equation}
\partial_r A_2= \Im{(\psi \bar{\psi}_2)}+\frac{1}r |\psi_2|^2, \qquad \partial_r
\psi_2 = i A_2 \psi - \frac{1}r A_2 \psi_2
\label{intro-comp}\end{equation}
with the conservation law $A_2^2 + |\psi_2|^2 = 1$. However, this 
does not uniquely determine $(\psi_2,A_2)$ and, by extension, the
Schr\"odinger map $u$ as we are missing a suitable boundary condition. 

\subsection{ Linearizations and the operators $H$, $\tilde H$}
This is the point in our work where we specialize in the case $m=1$
and, for convenience, drop the upper-script $m$ from all elements involved,
i.e. use $h_1,h_3$ instead of $h_1^1,h_3^1$, etc. 

A key role in our analysis is played by the linearization of the
Schr\"odinger Map equation around the soliton $Q$. A solution to the
linearized flow is a function
\[
u_{lin} : \R^{2+1} \to T_Q S^2.
\]
The Coulomb frame associated to $Q$ has the form 
\[
 v_Q(\theta,r) = e^{\theta R} \bar v_Q(r), \qquad w_Q(\theta,r) = e^{\theta R}
\bar w_Q(r)
\]
with 
\[
\bar v_Q(r) = 
\left( \begin{array}{ccc} h_3(r) \\ 0 \\ - h_1(r) \end{array} \right),
\qquad 
\bar w_Q(r) = 
\left( \begin{array}{ccc} 0 \\ 1 \\ 0 \end{array} \right).
\qquad 
\]
Expressing $u_{lin}$ in this frame,
\[
 \phi_{lin} = \la u_{lin},v_Q\ra + i \la u_{lin}, w_Q\ra
\]
one obtains the Schr\"odinger type equation
\begin{equation}
(i \partial_t - H) \phi_{lin} = 0 
\label{philin}\end{equation}
where the operator $H$ acting on radial functions has the form
\[
H = -\Delta + V, \qquad V(r) = \frac{1}{r^2} - \frac{8}{(1+r^2)^2}.  
\]

On the other hand linearizing the equation \eqref{intro-psieq} around the
soliton $Q$, we obtain a linear Schr\"odinger equation of the form
\begin{equation}
(i \partial_t - \tilde H) \psi_{lin} = 0 
\label{psilin}\end{equation}
where the operator $\tilde H$ acting on radial functions has the
form
\[
\tilde H = -\Delta + \tilde V, \qquad \tilde V(r) = \frac{2}{r^2}(1-h_3( r)) = 
\frac{4}{r^2(r^2+1)}. 
\]
The operators $H$ and $\tilde H$ are conjugate operators
and admit the factorizations
\[
H =  L^* L, \qquad \tilde  H = L  L^* 
\]
where
\[
L = h_1 \partial_r h_1^{-1} = \partial_r + \frac{h_3}{r}, \qquad
L^{*}=- h_1^{-1} \partial_r h_1 -\frac1{r}= -\partial_r +
\frac{h_3-1}{r}.
\]
The linearized variables 
$\phi_{lin}$ and $\psi_{lin}$ are also conjugated variables, 
\begin{equation} \label{linL}
 \psi_{lin} = L\phi_{lin}.
\end{equation}

The operator $H$ is nonnegative and bounded from $\dot H^1$ to $\dot H^{-1}$,
but it is not positive definite; it has a zero resonance
$\phi_0$, solving $L \phi_0=0$, namely 
\[
 \phi_0(r) = \frac{2r}{1+r^2}.
\]
This corresponds to the solution $\phi_{lin}$ for \eqref{philin}
obtained by differentiating
the soliton family with respect to either parameter. A consequence 
of this is that the linear Schr\"odinger evolution \eqref{philin}
does not have good dispersive properties, a fact which  is at the heart of our
instability result.

The above heuristic linearization argument works for all $m$, with the proper adjustments.
We remark that if $m \geq 2$ then the zero resonance is replaced by a zero
eigenvalue. If $m \geq 3$ then this eigenvalue belongs to $\dot H^{-1}$,
which allows for a clean splitting of the $\dot H^1$ space into 
an eigenvalue mode, which is stationary, and an orthogonal component,
which has good dispersive properties. This leads to the stability 
results in \cite{gkt}, \cite{gnt}. If $m=2$ we expect results which are closer
to the $m=1$ case; this will be considered in subsequent work.

If $Q$ is replaced by $Q_{\alpha, \l}$ then $H$ and $\tilde H$ are
replaced by their rescaled versions $H_\lambda$ and $\tilde H_\lambda$
where $V$ and $\tilde V$ are replaced by
\[
V_\lambda = \lambda^2 V(\lambda r), \qquad \tilde V_\lambda =
\lambda^2 \tilde V(\lambda r).
\]

A first goal of Section~\ref{spectral} is to describe the spectral theory
for the linear operators $H$ and $\tilde H$. The analysis in the case
of $H$ has already been done in \cite{KST}, and it is easily obtained
via the $L$ conjugation in the case of $\tilde H$. The normalized
generalized eigenfunctions for $H$ and $\tilde H$ are denoted by
$\phi_\xi$, respectively $\psi_\xi$, and satisfy
\[
 H \phi_\xi = \xi^2 \phi_\xi, \qquad  \tilde H \psi_\xi = \xi^2 \psi_\xi,
\qquad L \phi_\xi = \xi \psi_\xi.
\]
Correspondingly we have a generalized Fourier transform $\F_H$
associated to $H$ and a generalized Fourier transform $\F_{\tilde H}$
associated to $\tilde H$.

This quickly leads to generalized eigenfunctions for the rescaled
operators $H_\lambda$ and $\tilde H_\lambda$. A considerable effort is
devoted to the study of the transition from one $\tilde H_\lambda$
frame to another. This is closely related to the transference operator
introduced in \cite{KST}.

One reason we prefer to work with the $\psi$ variable is that 
the operator $\tilde H$ has a good spectral behavior at zero,
therefore we have favorable decay estimates for the corresponding 
linear Schr\"odinger evolution \eqref{psilin}.

\subsection{ The $X$ and $LX$ spaces}

As mentioned before, a stumbling block in formulating a closed
evolution equation for $\psi$ is the need for some boundary condition
in order to insure uniqueness for the system \eqref{intro-comp}.  This
leads us to introduce a stronger topology $X \subset \dHe$ for $\bar u
- \bar Q$, and therefore also for $\psi_2 - i h_1$ and $A_2-h_3$.  Then the
relation \eqref{linL} shows that studying the Schr\"odinger map
equation in the space $X$ corresponds to studying the $\psi$ equation
\eqref{intro-wnlin-eq1} in the space $LX$ obtained by applying the
operator $L$ to functions in $X$.

Roughly speaking the space $X$ is maximal with the following
properties:
\begin{itemize}
\item[(a)] We have the embedding $X \subset \dHe$.
\item[(b)] The $X$ norm of $u$ depends only on the the $L^2$ norms 
of the dyadic pieces of $\F_H u$.
\item[(c)] The operator $L$ is surjective on $X$.
\end{itemize}

Part (b) quickly implies a similar property for $LX$ relative to the
$\tilde H$ Fourier transform $\F_{\tilde H}$. It also shows that the
linear equations \eqref{philin}, respectively \eqref{psilin} are
well-posed in $X$, respectively $LX$.

One of the  goals of Section~\ref{spectral} is to define the $X$ and $LX$ 
spaces. In particular we establish the embedding \eqref{Xemb} for $X$,
as well as a two sided embedding for $LX$, namely
\begin{equation}\label{intro-l1emb}
L^1 \cap L^2 \subset LX \subset L^2.
\end{equation}
We also establish some other simple properties of these spaces. 

A key gain due to the fact that we work in the smaller space $X$ is
that we can supplement the system \eqref{intro-comp} with a boundary
condition at infinity, namely
\begin{equation} \label{intro-bciX}
 \lim_{r \to \infty} A_2 = 1, 
\qquad \psi_2 - i h_1 \in X. 
\end{equation}
This condition is preserved dynamically along the
Schr\"odinger map flow. Together with \eqref{intro-psieq}, \eqref{intro-comp}
and \eqref{intro-a0bis} it fully describes the dynamics of $\psi$.
Most of the work in this article is devoted to the study of the 
evolution of  $\psi$.

\subsection{ The elliptic transition between $u$ and its reduced field $\psi$}
Section~\ref{elliptic} is devoted to the study of the elliptic gauge 
correspondence at fixed time between the map $u$ and its 
associated reduced field $\psi$. The main result there asserts
that this map is a local $C^1$ diffeomorphism from a neighborhood
of the soliton $\bar Q$ in $X$ to a neighborhood of $0$ in $LX$.
As an intermediate step we prove that the system \eqref{intro-comp}
with the boundary condition \eqref{intro-bciX} yields a $C^1$ map 
from $\psi$ near $0$ in $LX$ to $(\psi_2,A_2)$ near $(h_1,h_3)$ in $X$.

\subsection{ The nonlinear Schr\"odinger equation for $\psi$: Take 1 [local]}
The equation \eqref{intro-psieq} can be rewritten in the form 
\begin{equation}
(i \partial_t - \tilde H) \psi = W \psi, \qquad 
W = A_0 - 2 \frac{A_2-h_3}{r^2}-\frac{1}r \Im({\psi}_2 \bar{\psi}).
\qquad 
\label{intro-wnlin-eq1}\end{equation}
Ideally, one would hope to be able to solve this equation 
in the space $LX$ by treating the right hand side perturbatively.
This is acceptable for  short time, and it provides us with a quick 
local theory.

The first step toward this goal is achieved in Section~\ref{linear} we
consider the linear Schr\"odinger evolution \eqref{psilin} and prove
Strichartz and local energy estimates.  Based on these bounds we
introduce function spaces $l^2 \Ss\subset L^\infty L^2$, respectively
$l^2 \Ns \supset L^1 L^2$ for $L^2$ solutions, respectively for the
inhomogeneous term in the $\tilde H$ Schr\"odinger equation.
Corresponding to $LX$ data we define similar weighted norms $
\WSs\subset L^\infty  LX$, respectively $  \WNs \supset L^1 LX$. 

In the beginning of Section~\ref{nonlin} we use these spaces and a
short fixed point argument to prove small data local well-posedness
for the equation \eqref{intro-wnlin-eq1} in $LX$.  Unfortunately, such
an argument no longer works globally in time; this is due to the
failure of the local decay estimates for $A_2-h_3$. While local decay
estimates are valid for $\psi$, they do not transfer to $\psi_2-ih_1$ and
$A_2-h_3$ via the ODE \eqref{intro-comp}-\eqref{intro-bciX}.

\subsection{ The functions $(\alpha(t),\lambda(t))$}
A primary goal of this article is to track the drift of Schr\"odinger
maps along the soliton family. For this we need appropriate functions
$(\alpha(t),\lambda(t))$ so that \eqref{goodal} holds.  The role of
$(\alpha(t),\lambda(t))$ is roughly to describe the low frequency
oscillations of the Schr\"odinger map $u$ along the family of rescaled
solitons.

In the case $m \geq 3$ the parameter $\lambda$ is defined dynamically
via an orthogonality condition with respect to the $H$ eigenvalue
$\phi_0$, appropriately rescaled (see \cite{gkt}). Such a strategy
cannot work for $m = 1,2$ as in this case $\phi_0 \not \in \dot
H^{-1}$. Another alternative would be to choose $\alpha$ and $\lambda$
as the minimizers in the left hand side in \eqref{stab}. However the
above minimizer plays no other  role, and in fact it turns
out that choosing it as the "closest" harmonic map to $u(t)$ may not
be the best choice for other analytical considerations, see \cite{gkt}
or \cite{gnt}.
 
In the context of this paper, it is technically convenient to make a
choice for $(\alpha(t),\lambda(t))$ which is expressed in terms of
$(\psi_2, A_2)$ instead of $u$.  Precisely, we make a dynamic
assignation of $(\alpha(t),\lambda(t))$ via the relation
\begin{equation}
A_2(1,t) = h_3(\lambda(t)), \qquad  \psi_2(1,t) = i e^{ i \alpha(t)}
h_1(\lambda(t))
\label{intro-analdef}\end{equation}
which for a soliton simply recovers the soliton parameters.
The (small) price to pay is that we need to prove that
\eqref{goodal} holds; we do this right away in the next section.
The choice of $r=1$ above is arbitrary; different choices of $r$ 
lead to closely related functions $\lambda$.

\subsection{ The nonlinear Schr\"odinger equation for $\psi$: Take 2 [global]}
With $\lambda(t)$ defined as in \eqref{intro-analdef}, the equation
\eqref{intro-psieq} can also be rewritten in the form
\begin{equation}
(i \partial_t - \tilde H_\lambda) \psi = W_\l \psi, \qquad 
W_\l = A_0 - 2 \frac{ A_2-h_3^\lambda}{r^2}-\frac{1}r \Im({\psi}_2 \bar{\psi}).
\label{intro-wnlin-eq}\end{equation}
The advantage is that, for $\lambda$ defined as in
\eqref{intro-analdef}, the ODE \eqref{intro-comp}-\eqref{intro-bciX}
allows us to transfer local energy decay estimates from $\psi$ to
$A_2-h_3^\lambda$ as the latter vanishes now at $r=1$ instead of
infinity. This is achieved in Proposition~\ref{p:leA2}.

 The price to pay is that we now need to understand the linear
evolution
\begin{equation}
(i \partial_t - \tilde H_\lambda) \psi  = f
\label{ldep}\end{equation}
in the space $LX$, with $\lambda$ depending on time. We expect $\lambda$
to stay bounded, but this is far from being enough. Instead we introduce 
a smaller space $Z_0 \subset L^\infty$ for $\lambda-1$, defined by
\[
Z_0 = \dot H^{\frac12,1} + W^{1,1} + (L^2 \cap L^\infty).
\]
Here the last component characterizes the high frequencies (which have
good averaged decay), while the first two apply primarily for the low
frequencies (and have little decay at infinity).

This is achieved in Section~\ref{lineart}, where we consider the global in
time linear Schr\"odinger evolution \eqref{ldep} under the assumption
that $\lambda -1 $ is small in the space $Z_0$.  We construct function
spaces $\WSs[\tilde \lambda] \subset L^\infty (L X)$, respectively
$\WSs[\tilde \lambda] \supset L^1 LX$, incorporating also appropriate
dispersive information, so that the following linear bound holds for
solutions to \eqref{ldep}:
\begin{equation} \label{ldep-lin}
 \| \psi \|_{\WSs[\tilde \lambda]} \lesssim \|\psi(0)\|_{L X} + \| f\|_{\WNs[\tilde \lambda]}.
\end{equation}
Here $\tilde \lambda$ is a (somewhat arbitrary) regularization of $\lambda$
which essentially contains the low frequencies of  $\lambda$.

Section~\ref{nonlin} contains our global in time analysis of the
nonlinear equation for $\psi$. Precisely, we establish a bootstrap
estimate for the $\WSs[\tilde \lambda]$ size of $\psi$. This is obtained by
combining the linear bound \eqref{ldep-lin} with an estimate for the
nonlinearity, which has the form
\[
\|W_\lambda \psi\|_{\WNs[\tilde \lambda]} \lesssim \|\psi\|_{\WSs[\tilde \lambda]}^2.
\]
We remark that while we are able to prove a bootstrap estimate for
solutions $\psi$ to \eqref{intro-wnlin-eq}, we cannot obtain estimates
for the difference of two solutions.  Hence we can no longer treat the
nonlinearity perturbatively globally in time.

In Section~\ref{seclambda} we complement the above bootstrap estimate
for $\psi$ with a bootstrap estimate for $\|\lambda -1\|_{Z_0}$. More
precisely we show that we recover the regularity of the parameter
$\lambda(t)$ from the $\WSs[\tilde \lambda]$ regularity of $\psi$.

Finally, in Section~\ref{secboot} we prove our main stability result
in the $X$ topology in Theorem~\ref{tmain-X}. This is done via a
bootstrap argument, which uses the bootstrap estimates on $\psi$ in
$\WSs[\tilde \lambda]$, respectively for $\lambda-1$ in $Z_0$, from the previous
two sections. In addition, we use the results in
Section~\ref{elliptic} for the transition back and forth between the
Schr\"odinger map $u$ and its reduced field $\psi$.

\subsection{ The instability result}
In the final section of the paper we prove the $\dot H^1$ instability
result in Theorem~\ref{tmain-H}.  For this we  introduce a
second small parameter $\epsilon$ and look at maps $u$ for which the
reduced field $\psi$ satisfies
\[
\| \psi(t)\|_{LX} \approx \gamma, \qquad \|\psi(t)\|_{L^2} \approx \epsilon \gamma
\]
The $L^2$ smallness allows for a better control of the nonlinear effects,
and we are able to show that the $\psi$ flow is almost linear,
\[
\| \psi(t) - e^{i t \tilde H} \psi(0)\|_{LX} \lesssim |\log \epsilon|^{-1} \gamma
\]
Taking this into account, for each $(\alpha,\lambda)$ as in
\eqref{alinst} our strategy is to choose an initial data $u_0$ which
coincides with $Q^{1}_{\alpha,\lambda}$ for $r < \epsilon^{-1}$ and
with $Q=Q^{1}_{0,1}$ for larger $r$, with a smooth transition
in between. Then we are able to accurately track the Fourier transform
$\F_{\tilde H} \psi$ of $\psi$ for large $t$. The decay of the map $u$
to an $O_{\dot H^1} ( |\log \epsilon|^{-1} \gamma)$ neighborhood of
$Q$ is equivalent to the decay of $(\alpha(t),\lambda(t))$ to
an an $O ( |\log \epsilon|^{-1} \gamma)$ neighborhood of $(0,1)$,
which in turn is a consequence of cancellations due to the
oscillations in frequency for $e^{i t \tilde H} \psi(0)$ as $t$ grows
large.

\section{The Coulomb gauge representation of the equation} 
\label{seccoulomb}

In this section we rewrite the Schr\"odinger map equation for
equivariant solutions in a gauge form.  This approach originates in
the work of Chang, Shatah, Uhlenbeck~\cite{csu}. However, our analysis is
closer to the one in \cite{bik}. 

\subsection{Near soliton maps}

We first investigate some simple properties of maps $u:\R^2 \to \S^2$
which are near a soliton $Q^m_{\alpha,\lambda}$ in the sense that
\begin{equation}
\|u - Q_{\alpha,\lambda}^m\|_{\dot H^1} \leq \gamma \ll 1.
\label{l1a}\end{equation}
\begin{l1}
 Let $m \geq 1$ and $u: \R^2 \to \S^2$ be an $m$-equivariant map
which satisfies \eqref{l1a}. Then
\begin{equation}
\lim_{r \to 0} u(r,\theta) = -\vec{k}, \qquad \lim_{r \to \infty} u(r,\theta) = \vec{k}
\label{l1b}\end{equation}
and 
\begin{equation}
\| r^{-1} ( u - Q_{\alpha,\lambda}^m)\|_{L^2}+  \| u - Q_{\alpha,\lambda}^m\|_{L^\infty} \lesssim \gamma.
\label{l1c}\end{equation}
\label{l:ns}\end{l1}

\begin{proof}
 After a rescaling and a rotation the problem reduces to the case
  $\alpha = 0$ and $\lambda = 1$. We rewrite the $\dot H^1$ bound for
  $u-Q^{m}$ as in \eqref{energy}:
\begin{equation}
 \| \partial_r (\bar u - \bar Q^{m})\|_{L^2} + 
\|r^{-1} (\bar u_1-h_1^m)\|_{L^2}+
\|r^{-1} \bar u_2\|_{L^2} \lesssim \gamma 
\label{bu}\end{equation}
In particular for $u_2$ we have
\begin{equation}
 \| \partial_r \bar u_2\|_{L^2} + \|r^{-1} \bar u_2\|_{L^2} \lesssim \gamma 
\label{bu2}\end{equation}
By Sobolev type embeddings
\begin{equation}\label{sobo1}
\|f\|_{L^\infty} \lesssim  \| \partial_r f \|_{L^2} + \|r^{-1} f\|_{L^2}
\end{equation}
 this implies that 
\[
 \| \bar u_2\|_{L^\infty} \lesssim \gamma
\]
Furthermore,
\[
\frac{d}{dr} |\bar u_2|^2 = 2 \bar u_2 \partial_r \bar u_2 \in L^1(dr)
\]
therefore $|\bar u_2|^2$ is continuous and has limits as $r \to 0,\infty$.
In addition, these limits must be zero in order for the second 
left hand side norm in \eqref{bu2} to be finite.
The same argument applies for $\bar u_1 - h_1^m$. Thus we have proved 
that
\begin{equation}
 \| \bar u_2\|_{L^\infty} + \| \bar u_1 -h_1^m \|_{L^\infty} \lesssim \gamma, \qquad \lim_{r \to 0,\infty} \bar u_1, \bar u_2 = 0
\label{bu23}\end{equation}
To conclude the proof of \eqref{l1b} and \eqref{l1c} it remains 
to consider the vertical component $\bar u_3(r)$. Integrating 
the bound 
\[
 \| \partial_r (\bar u_3-h_3^m)\|_{L^2} \lesssim \gamma
\]
we obtain  (as $h_3^m(1) = 0$)
\[
 |(\bar u_3(r) - h_3^m(r))- u_3(1)| \lesssim \gamma |\log r|^\frac12
\]
The first part of \eqref{bu23} shows that $|\bar u_3(1)|^2 \lesssim
\gamma$. Since $\gamma$ is small, it follows that $\bar u_3(r)$ is
negative for say $r \in [\frac14,\frac12]$.  Since $\bar u_3$ is
continuous and, by \eqref{bu23}, cannot vanish for smaller $r$, it
follows that it stays negative for all $r \in (0,\frac12]$. Thus for
$r < 1/2$ we have
\[
\bar u_3(r) = - \sqrt{1- |\bar u_1|^2 -|\bar u_2|^2} 
\]
which by  \eqref{bu23} implies that 
\[
 |\bar u_3(r) - h_3^m(r)| \lesssim \gamma, \quad r < \frac12, 
\qquad \lim_{r \to 0} \bar u_3(r)= -1.
\]
The same argument applies for $r > 2$, where $\bar u_3$ is positive.
Integrating its $r$ derivative from either side we recover the pointwise bound for $\bar u_3-h_3^m$ for $r \in
[\frac12,2]$ and obtain
\[
 \|\bar u_3 - h_3^m\|_{L^\infty} \lesssim \gamma
\]
Finally, we consider the $L^2$ bound for $\bar u_3 - h_3^m$. 
Due to the pointwise bound, it suffices to consider $r$ close 
to $0$ and to infinity. In either case we use the equation of the 
sphere to write
\[
 |\bar u_3 - h_3^m| \lesssim |\bar u_1 - h_1^m| +|\bar u_2|
\]
and conclude by \eqref{bu}.
\end{proof}

\subsection{The Coulomb gauge} \label{CG}

 We let the differentiation operators  $\partial_0,\partial_1,\partial_2$
stand for $\partial_t, \partial_r, \partial_t$ respectively. Our strategy 
will be to replace the equation for the Schr\"odinger map $u$ with 
equations for its derivatives $\partial_1 u$, $\partial_2 u$ 
expressed in an orthonormal frame $v,w \in T_u \S^2$. To fix the 
sign in the choice of $w$, we will assume that
\[
u \times v = w
\]
 Since $u$ 
is $m$-equivariant it is natural to work with $m$-equivariant frames, i.e. 
 \[
v = e^{m \theta R} \bar{v}(r), \qquad w = e^{m \theta R} \bar{w}(r).
\]
Given such a frame  we introduce
the differentiated fields $\psi_k$ and the connection coefficients
$A_k$ by
\begin{equation}
\label{connection}
\begin{split}
\psi_k = \partial_k u \cdot v + i \partial_k u \cdot w, \qquad
A_k = \partial_k v \cdot w.
\end{split}
\end{equation}
Due to the equivariance of $(u,v,w)$ it follows that both $\psi_k$ and
$A_k$ are spherically symmetric (therefore subject to the 
conventions made in Section \ref{defnot}). Conversely, given $\psi_k$ and $A_k$
we can return to the frame $(u,v,w)$ via the ODE system:
\begin{equation}
\label{return}
\left\{ \begin{array}{l}
\partial_k u = (\Re {\psi_k}) v  + (\Im{\psi_k}) w
\cr
\partial_k v = - (\Re{\psi_k}) u + A_k w
\cr
\partial_k w = - (\Im{\psi_k}) u - A_k v
\end{array} \right.
\end{equation}

If we introduce the covariant differentiation
\[
D_k = \partial_k + i A_k, \ \ k \in \{0,1,2\}
\]
it is a straightforward computation to check the compatibility conditions:
\begin{equation} \label{compat}
D_l \psi_k = D_k \psi_l, \ \ \  l,k=0,1,2.
\end{equation}
The curvature of this connection is given by 
\begin{equation}\label{curb}
  D_l D_k - D_k D_l = \partial_l A_k - \partial_k
A_l = \Im{(\psi_l \bar{\psi}_k)}, \ \ \  l,k=0,1,2.
\end{equation}
An important geometric feature is that $\psi_2, A_2$ are closely
related to the original map. Precisely, for $A_2$ we have:
\begin{equation}
A_2 = m (k \times v) \cdot w= m k \cdot (v \times w)= m k \cdot u = m u_3
\label{a2u3}\end{equation}
and, in a similar manner,
\begin{equation}
\psi_2 = m(w_3-iv_3)
\label{psi2vw3}\end{equation}
Since the $(u,v,w)$ frame is orthonormal, the following relations
also follow: 
\[
|\psi_2|^2 = m(u_1^2 + u_2^2), \qquad
|\psi_2|^2 + A_2^2 = m^2
\]

Now we turn our attention to the choice of the $(\bar v,\bar w)$ frame
at $\theta = 0$. Here we have the freedom of an arbitrary rotation
depending on $t$ and $r$. In this article we will use the Coulomb
gauge, which for general maps $u$ has the form $\text{div } A = 0$.
In polar coordinates this is written as $\partial_1 A_1 + r^{-2} 
\partial_2 A_2 = 0$. However, in the equivariant case $A_2$ is radial,
so we are left with a simpler formulation $A_1 = 0$, or equivalently
\begin{equation}
\partial_r \bar v \cdot \bar w=0
\label{coulomb}\end{equation}
which can be rearranged into a convenient ODE as follows
\begin{equation} \label{cgeq}
\partial_r  \bar v = (\bar v \cdot \bar u) \partial_r \bar u
- (\bar v \cdot \partial_r \bar u) \bar u
\end{equation}
The first term on the right vanishes and could be omitted,
but it is convenient to add it so that the above linear ODE is solved 
not only by $v$ and $w$, but also by $u$. Then we can write an equation 
for the matrix $ \calO = (\bar v, \bar w,\bar u)$:
\begin{equation} \label{cgeq-m}
\partial_r \calO = M \calO, \qquad M = \partial_r \bar u \wedge \bar u : = 
\partial_r \bar u \otimes \bar u - \bar u \otimes \partial_r \bar u
\end{equation}
with an  antisymmetric matrix $M$. 

An advantage of using the Coulomb gauge is that it makes the
derivative terms in the nonlinearity disappear. Unfortunately,
this only happens in the equivariant case, which is why in \cite{BIKT}
we had to use a different gauge, namely the caloric gauge.

The ODE \eqref{cgeq} needs to be initialized at some point.
 A change in the initialization leads to a multiplication
of all of the $\psi_k$ by a unit sized complex number. This is 
irrelevant at fixed time, but as the time varies we need to be careful 
and choose this initialization uniformly with respect to $t$, in 
order to avoid introducing a constant time dependent potential 
into the equations via $A_0$.  Since in our results we start with
data which converges asymptotically to $Q$ as $r \to \infty$,
and the solutions continue to have this property, it is natural to
fix the choice of $\bar{v}$ and $\bar{w}$ at infinity,
\begin{equation}
\label{bcvw}
 \lim_{r \to \infty} \bar v(r,t) = \vec{i} , \qquad \lim_{r \to \infty} \bar w(r,t) =
\vec{j}
\end{equation}
Before considering the general case we begin with the solitons.
The simplest case is when $u = Q^{m}$ when the triplet 
$(\bar v,\bar w,\bar u)$ is given by
\begin{equation}\label{vwq}
\left(  \bar V^{m}, \bar W^{m},\bar Q^{m}\right)
= 
  \left(\begin{array}{ccc}
   h_3^m(r) & 0 & h_1^m(r) \cr 0 & 1 & 0 \cr
  -h_1^m(r) & 0  & h_3^m(r)
\end{array}\right)
\end{equation}
If $m=1$ we drop the superscript $m$.
More generally, if $ u = Q^{m}_{\alpha,\lambda}$ then from the above,
by rescaling and rotation, we obtain the corresponding triplet $
\left(\bar V^{m}_{\alpha,\lambda}, \bar W^{m}_{\alpha,\lambda},\bar
  Q^{m}_{\alpha,\lambda},\right)$ of the form
\[
 \!\left(\!\!\! \begin{array}{ccc}
\! h_3^m(\lambda r)\cos^2 m\alpha + \sin^2 m\alpha 
 & \! (h_3^m(\lambda r)-1) \sin m\alpha \cos m \alpha
& h_1^m(\lambda r) \cos m\alpha  \cr
 \! (h_3^m(\lambda r)-1) \sin m \alpha \cos m \alpha \!
& \!  h_3^m(\lambda r) \sin^2 m\alpha + \cos^2 m\alpha 
& h_3^m(\lambda r) \sin m\alpha  \cr
\! -h_1^m(\lambda r) \cos m\alpha & \!
- h_1^m(\lambda r) \sin m \alpha & h_3^m(\lambda r) \end{array}\!\!\! \right)\!
\]
For later reference we also note the values of $\psi_1$, $\psi_2$ and $A_2$ 
in this case:
\begin{equation} \label{solref}
\begin{split}
\psi_{\alpha,\lambda,1}^m = - m r^{-1} h_1^m(\lambda r) e^{im\alpha}, &\quad \psi_{\alpha,\lambda,2}^m= i m h_1^m(\lambda r) e^{im\alpha},
\\
 A_{\alpha,\lambda,2}^m =& m h_3^m(\lambda r).  
\end{split}
\end{equation}
To measure the regularity of the frame $(v,w)$ we use the 
Sobolev type space $\dot H^1_C$ of functions $f: \R^2 \to \R^3$,
with norm
\[
\| f\|_{\dot H^1_C} = \| \partial_{r}  \bar f \|_{L^2} + 
\|\bar f\|_{L^\infty} + \|r^{-1} \bar f_3\|_{L^2}, \qquad f(r,\theta)=e^{m\theta R}\bar f(r)
\]
The next Lemma shows that the initialization \eqref{bcvw} is  well-defined for arbitrary maps $u$ close to the soliton family:

\begin{p1} \label{lc}
a) For each  $m$-equivariant map $u: \R^2 \to \S^2$ 
satisfying \eqref{l1a} there exists an unique $m$-equivariant orthonormal frame $(v,w)$ which satisfies the Coulomb gauge condition \eqref{coulomb} and the boundary condition \eqref{bcvw}. This frame satisfies the bounds
\begin{equation}
 \| v - V^m_{\alpha,\lambda}\|_{\dot H^1_C} + \| w- W^m_{\alpha,\lambda}\|_{ \dot H^1_C}\lesssim \gamma.
\label{l1e}\end{equation}
b) Furthermore, the maps $u \to v,w$ are $C^1$ from $\dot{H^1}$ 
into $ \dot H^1_C$ as well as from $L^2 \to L^2$. 
\label{p:gauge}\end{p1}

\begin{proof}
a) To construct the Coulomb gauge we use the equation
\eqref{cgeq}. The right hand side is linear in $\bar v$ and has locally
integrable coefficients, therefore by prescribing $\bar v$ at $r = 1$
we obtain a unique solution. Also, if the relations
\begin{equation}
|\bar v|^2 = 1, \qquad \bar u \cdot \bar v=0
\label{vconstr}\end{equation}
 are enforced at $r = 1$
then they are preserved along the flow. We claim that the limit
of $\bar v(r)$ as $r \to \infty$ exists. For $\bar v_2$ and $\bar v_3$
this follows from
\[
 \| \partial_r \bar v_j\|_{L^1(dr)} \lesssim \| \partial_r \bar u\|_{L^2(rdr)}
\| r^{-1} \bar u_j\|_{L^2(rdr)}, \qquad j = 1,2.
\]
On the other hand $\lim_{r \to \infty} \bar v_3 = 0$  by orthogonality
due to the relation \eqref{l1b}.

Once we have one solution $\bar v$ to \eqref{cgeq},  a second one 
is obtained by $\bar w = \bar u \times \bar v$. Since \eqref{cgeq}
is linear, it follows that all its solutions are obtained from the 
initial one by a rotation of a fixed angle in $T_{\bar u} \S^2$.
This proves the existence and uniqueness of the desired frame
$(v,w)$ which satisfies the boundary condition \eqref{bcvw}.

We next prove the pointwise part of the bound \eqref{l1e}. From \eqref{cgeq} we obtain
\[
 \partial_r \bar v_2 = -\left(\bar v \cdot \partial_r (\bar u - \bar Q^{m})\right) \bar u_2 -  (\bar v \cdot\partial_r \bar Q^{m}) \bar u_2 
\]
Hence using \eqref{l1a} we estimate $ \| \partial_r \bar v_2\|_{L^1(dr)}
\lesssim \gamma $, which after integration shows that $ | \bar v_2|
\lesssim \gamma $. Since we also have $|\bar u_2| \lesssim \gamma$, it
follows that $|\bar w_2-1| \lesssim \gamma^2$. This in turn shows that
$|\bar w_1|+ |\bar w_3| \lesssim \gamma$.  
Then the pointwise bounds for $\bar v_1$
and $\bar v_3$ are easily obtained since $v = - u \times w$.

Next we consider the $L^2$ bounds for $\partial_r \bar v$ and
$\partial_r \bar w$.  The easy case is that of $\bar v_2$ and $\bar w_2$, for
which by \eqref{cgeq} we have
\[
 \|\partial_r \bar v_2\|_{L^2(rdr)} + \|\partial_r \bar w_2\|_{L^2(rdr)}
\lesssim \|\partial_r \bar u_2\|_{L^2(rdr)} \lesssim \gamma
\]
For $\bar v_1$ we write
\[
 \partial_r (\bar v_1 - h_3^m) =  
(\bar v \cdot \partial_r (\bar Q^m-\bar u)) \bar u_1 
- (\bar v \cdot \partial_r  \bar Q^m) (\bar u_1-h_1^m) + 
((\bar V^m - \bar v) \cdot \partial_r  \bar Q^m) h_1^m 
\]
For the first term we use \eqref{cgeq} while for the remaining terms we 
use the pointwise bounds in \eqref{l1c} and \eqref{l1e} 
for $\bar u_1 - h_1^m$, respectively
$\bar V^m - \bar v$. 
The same argument applies for $\bar v_3$, $\bar w_1$ and $\bar w_3$.

Finally, we prove the $L^2$ bounds for $\bar v_3$ and $\bar w_3$. 
This is done 
in a roundabout way using the orthogonality of the 
$(\bar u,\bar v,\bar w)$ frame.
For $\bar w_3$ we have
\[
|\bar w_3| = | \bar u_1 \bar v_2 - \bar u_2 \bar v_1| 
\lesssim |\bar u_1 - h_1^m| + h_1^m |\bar v_2| + |\bar u_2|
\]
and we conclude using the $L^2$ bounds in \eqref{bu} for $\bar
u_1-h_1^m$ and $\bar u_2$ as well as the pointwise bound for $\bar v_2$
in \eqref{l1e}.  Similarly, for $\bar v_3$ we have
\[
|\bar v_3+h_1^m| = | \bar u_1 \bar w_2 - \bar u_2 \bar w_1 + h_1^m| \lesssim 
|\bar u_1 - h_1^m| + h_1^m |\bar w_2-1| + |\bar u_2|
\]
and we conclude as before.  The proof of \eqref{l1e} is complete.

b) We now prove that the map $u \to v$ is $C^1$ from $\dot H^1$ to
$\dot H^1_C$.  Given an interval $I$ we consider a one parameter
family of maps $ u: I \times \R^2 \to \S^2$ which are smooth in all
variables and agree with $Q^{m}$ for large $r$. Then by 
ODE theory applied to \eqref{cgeq} it follows that $v$ and $w$ are
smooth in all variables away from $r=0$. The main step is to establish
the uniform bounds
\begin{equation}\label{clin}
 \| \partial_t v\|_{\dot H^1_C} \lesssim \| \partial_t u\|_{\dot H^1}
\end{equation}
Having this, the transition to more general maps $u \in C^1(I;\dot
H^1)$ is done via a standard density argument, which is omitted. We
remark that the $\dot H^{1}$ convergence for $u$ and the $\dot
H^{1}_C$ convergence for $v,w$ suffice in order to insure that the ODE
\eqref{cgeq} and the boundary condition \eqref{bcvw} survive in the
limit.

To prove \eqref{clin} we differentiate \eqref{cgeq} with respect to $t$
to obtain an ODE for the covariant time derivative of $v$, 
\[
 z = \partial_t v + (v \cdot \partial_t u) u
\]
We obtain
\begin{equation} \label{cgeqz}
\partial_r   \bar z = (\bar z \cdot \bar u) \partial_r \bar u
- (\bar z \cdot \partial_r \bar u) \bar u  + \bar f, \qquad z(\infty) = 0
\end{equation}
where
\[
 f =  (\partial_r v \cdot \partial_t u)u + (v \cdot \partial_t u) \partial_r u - (v \cdot \partial_r u) \partial_t u 
\]
For $\partial_t u$ we use the following bound:
\begin{equation} \label{ptu}
\| \partial_t \bar u\|_{L^\infty} + \| r^{-1} \partial_t \bar u\|_{L^2}
\lesssim \| \partial_t u\|_{\dot H^1}
\end{equation}
For $\partial_t \bar u_1$ and $\partial_t \bar u_2$ this follows 
directly due to the form \eqref{energy} of the $\dot H^1$ norm
for equivariant functions and to \eqref{sobo1}.
On the other hand the bound for $\partial_t \bar u_3$ is obtained 
indirectly from the orthogonality relation $\partial_t u \cdot u = 0$
(see e.g. the similar argument for $u_3$ in Lemma~\ref{l:ns}).

From the $L^2$ bound in \eqref{ptu} we obtain
\[
 \| \bar f\|_{L^1(dr)} \lesssim \| \partial_t u\|_{\dot H^1}
\]
therefore integrating \eqref{cgeqz} from infinity we 
have 
\[
 \| z\|_{L^\infty} \lesssim \| \partial_t u\|_{\dot H^1}
\]
Using the $L^\infty$ bound  in \eqref{ptu} yields
\[
 \| \bar f\|_{L^2(rdr)} \lesssim \| \partial_t u\|_{\dot H^1}
\]
which directly leads to 
\[
 \|\partial_r z\|_{L^2} \lesssim \| \partial_t u\|_{\dot H^1}
\]
Further, the orthogonality relation $z \cdot u = 0$ and \eqref{l1c}
show that
\[
 \|r^{-1} \bar z_3\|_{L^2(rdr)} \lesssim \|\bar z\|_{L^\infty}
(1 + \|r^{-1} \bar u_1\|_{L^2(rdr)}+ \|r^{-1} \bar u_2\|_{L^2(rdr)})
\lesssim \| \partial_t u\|_{\dot H^1}
\]
Thus we have proved that
\[
 \| z \|_{\dot H^1_C} \lesssim \| \partial_t u\|_{\dot H^1}
\]
Now it is easy to obtain \eqref{clin}, estimating the 
difference via \eqref{ptu}.

Finally, we prove that the map $u \to v,w$ is $C^1$ from $L^2 \to L^2$.
For this we need the following counterpart of \eqref{clin}: 
\begin{equation}\label{clina}
 \| \partial_t v\|_{L^2} \lesssim \| \partial_t u\|_{L^2}
\end{equation}
Again it suffices to consider the smooth case, since the transition to more
general maps $u \in C(I;\dot H^1) \cap C^1(I;L^2)$ is done via a
standard density argument.

We begin with 
\[
 \| \bar f\|_{L^1(rdr)} \lesssim \|\partial_t u\|_{L^2}
\]
Then integrating \eqref{cgeqz} from infinity we obtain
\[
|z(r)| \leq \int_{0}^\infty 1_{[0,s]}(r) |f(s)| ds 
\]
and by Minkowski's inequality 
\[
 \| z\|_{L^2(rdr)} 
\lesssim \int_0^\infty s |f(s)| ds
\]
The transition from $z$ to $\partial_t v$ is immediate, therefore
\eqref{clina} is proved.
\end{proof}

As a direct consequence of part (a) of the above lemma
we can describe the regularity and properties of the differentiated
fields $\psi_1$, $\psi_2$ and the connection coefficient $A_2$ at 
fixed time:

\begin{c1}
  Let  $u: \R^2 \to \S^2$ be an $m$-equivariant map
as in \eqref{l1a}. Then $\psi_1$, $\psi_2$ and $A_2$ satisfy
\eqref{compat}, \eqref{curb} for $k,l=1,2$
as well as the bounds
\[
 \| \psi_1- \psi_{\alpha,\lambda,1}^m\|_{L^2} + 
\|\psi_2 - \psi_{\alpha,\lambda,1}^m\|_{\dHe} +
\|A_2 - A_{\alpha,\lambda,2}^m\|_{\dHe} \lesssim \gamma
\]
In addition, the map $u \to (\psi_1, \psi_2, A_2)$ from 
$\dot H^1$ into the above spaces is $C^1$.
\end{c1}

A second step is to consider Schr\"odinger maps with more regularity,
i.e. as in Theorem~\ref{th2}.  For such maps,
if we make the additional decay assumption that 
$u_0 - Q^{m}_{\alpha,\lambda} \in L^2$, then this is preserved
along the flow. Hence, as a consequence of part (b) of the above lemma we 
have:

\begin{c1}\label{c:coulomb}
  Let $I$ be a compact interval, and  $u: I\times \R^2 \to \S^2$ be an $m$-equivariant map satisfying \eqref{l1a} uniformly in $I$ and which 
has the additional regularity
\[
 u - Q^{m}_{\alpha,\lambda} \in  C(I;H^2),\qquad  \partial_t u \in C(I;L^2).
\]
Then $\psi_0$, $\psi_1$, $\psi_2$ and  $A_0$, $A_2$ satisfy
the relations \eqref{compat}, \eqref{curb} for $k,l=0,1,2$
and have the additional regularity
\begin{equation} \label{Cextra}
\begin{split}
&  \psi_0, A_0 \in C(I;L^2), \psi_1- \psi_{\alpha,\lambda,1}^m \in C(I;\He), \\
& \psi_2 - \psi^{m}_{\alpha,\lambda,2},
A_2 - A^{m}_{\alpha,\lambda,2} \in C(I;\Hde).
\end{split}
\end{equation}
\end{c1}

\subsection{ Schr\"odinger maps in the Coulomb gauge}

We are now prepared to write the evolution equations for the
differentiated fields $\psi_1$ and $\psi_2$ in \eqref{connection}
computed with respect to the Coulomb gauge. To justify the following
computations we assume that $u:I \times \R^2 \to \S^2$ 
is a Schr\"odinger map as in Theorem~\ref{th2}
so that in addition $u_0-Q^m \in L^2$ . Thus the hypothesis
of Corollary~\ref{c:coulomb} is verified, and we obtain 
the additional regularity \eqref{Cextra} for $\psi_0$, $\psi_1$, $\psi_2$ and  $A_0$, $A_2$. This suffices in order to justify the 
computations below.

Writing the Laplacian in polar coordinates, a direct computation
using the formulas \eqref{connection} shows that we can rewrite the
Schr\"odinger Map equation \eqref{SM} in the form
\[
\psi_0 = i \left(D_1 \psi_1 + \frac{1}{r} \psi_1 + \frac{1}{r^2} D_2 \psi_2\right) 
\]
Applying the operators $D_1$ and $D_2$ to both sides of this equation and using the relations \eqref{compat} and \eqref{curb}, we can derive the evolution equations for $\psi_m$, $ m=1,2$:
\[
\begin{split}
\partial_t \psi_1 + i A_0 \psi_1 = & \ i \Delta \psi_1  - 2  A_1 \partial_1
\psi_1  - \partial_1 A_1 \psi_1 - \frac{1}{r} A_1 \psi_1 \\
& \ - i A_1^2 \psi_1 - i \frac{1}{r^2} A_2^2 \psi_1  - i \frac1{r^2} \psi_1
+ \frac2{r^3} A_2 \psi_2 - \frac{1}{r^2} \Im{(\psi_1 \bar{\psi}_2)} \psi_2 \\
\partial_t \psi_2 + i A_0 \psi_2 = & \ i \Delta \psi_2  - 2 A_1 \partial_1
\psi_2  - \partial_1 A_1 \psi_2 - \frac{1}{r} A_1 \psi_2 \\
& \ - i A_1^2 \psi_2 - i \frac{1}{r^2} A_2^2 \psi_2 - \Im{(\psi_2 \bar{\psi}_1)}
\psi_1 
\end{split}
\]
Under the Coulomb gauge $A_1 = 0$ these  equations become
\begin{equation*} \label{smg}
\begin{split}
\partial_t \psi_1 + i A_0 \psi_1 = & i \Delta \psi_1  - i \frac{1}{r^2} A_2^2
\psi_1  -i  \frac1{r^2} \psi_1
+  \frac2{r^3} A_2 \psi_2 - \frac{1}{r^2} \Im{(\psi_1 \bar{\psi}_2)} \psi_2 \\
\partial_t \psi_2 + i A_0 \psi_2 = & i \Delta \psi_2  - i \frac{1}{r^2} A_2^2
\psi_2 - \Im{(\psi_2 \bar{\psi}_1)} \psi_1 
\end{split}
\end{equation*}
while the relations  \eqref{compat} and \eqref{curb} become
\begin{equation} \label{comp}
\partial_r A_2= \Im{(\psi_1 \bar{\psi}_2)}, \qquad \partial_r \psi_2 = i A_2
\psi_1
\end{equation}
From the compatibility relations involving $A_0$, we obtain
\begin{equation} \label{a0}
\partial_r A_0= - \frac1{2r^2} \partial_r (r^2 |\psi_1|^2 - |\psi_2|^2)
\end{equation}
from which we derive 
\begin{equation}
 A_0 = - \frac12 \left( |\psi_1|^2 - \frac{1}{r^2}|\psi_2|^2\right)
+ [r\partial_r]^{-1} \left( |\psi_1|^2 - \frac{1}{r^2}|\psi_2|^2\right) 
\label{a0bis}\end{equation}
This is where the initialization of the Coulomb gauge 
at infinity is important. That guarantees that $A_0 \in L^2$,
while $|\psi_1|^2 - r^{-2}|\psi_2|^2 \in L^2$. Thus the integrating constant 
must be zero.

There is quite a bit of redundancy in the equations for $\psi_1$ and $\psi_2$; we eliminate this by introducing a single main variable
\[
 \psi=\psi_1 -i \frac{\psi_2}{r}
\]
A direct computation yields the equation for $\psi$:
\[
i \partial_t \psi  +  \Delta \psi = A_0 \psi - 2 \frac{A_2}{r^2} \psi +
\frac{1}{r^2} \psi + \frac{A_2^2}{r^2} \psi - \frac1{r} \Im{(\psi_2
\bar{\psi}_1)} \psi
\]
By replacing $\psi_1 = \psi + i r^{-1} \psi_2$ and using $A_2^2 + |\psi_2|^2=1$,
 we obtain the key evolution equation we work with in this paper,
\begin{equation} \label{psieq}
i \partial_t \psi  +  \Delta \psi - \frac{2}{r^2} \psi = A_0 \psi - 2
\frac{A_2}{r^2} \psi   - \frac1{r}\Im{(\psi_2 \bar{\psi})} \psi
\end{equation}
Our strategy will be to use this equation in order to obtain estimates for
$\psi$. The functions $A_2$ and $\psi_2$ are defined in terms of $\psi$
via the system of ODE's
\begin{equation}
 \label{comp1}
\partial_r A_2= \Im{(\psi \bar{\psi}_2)}+\frac{1}r |\psi_2|^2, \qquad \partial_r
\psi_2 = i A_2 \psi - \frac{1}r A_2 \psi_2
\end{equation}
derived from \eqref{comp}. If $m=1$, the boundary condition for this system
will be prescribed at infinity, and it roughly says that $(A_2,\psi_2)$
are close to $(h_3,i h_1)$ as $r$ approaches $\infty$. In the regular case
 when $u - Q \in  C(I;H^2)$ this is simply the following relation:
\begin{equation}\label{bc:psi2}
A_2 - h_3 \in L^2, \qquad \psi_2 - i h_1 \in L^2
\end{equation}
We will later prove that this condition suffices in order to uniquely
determine $\psi_2$ and $A_2$ from $\psi$. This can only work in
the $1$-equivariant case; indeed, if $m \geq 2$ then nearby solitons
cannot be differentiated in this way.

Once  $(A_2,\psi_2)$ are computed, the $A_0$
connection coefficient is given by \eqref{a0bis} which becomes now
\begin{equation} \label{aoef}
 A_0(r) =  -\frac12 |\psi|^2  + \frac{1}r \Im (\psi_2 \bar \psi)
+ [r\partial_r]^{-1}( |\psi|^2 - \frac{2}{r} \Im (\psi_2 \bar \psi)).
\end{equation}
Finally, given $\psi$, $A_2$ and 
$\psi_2$, we can return to the  Schr\"odinger map $u$
via the system \eqref{return} with the boundary condition at infinity
given by \eqref{bcvw}.

\subsection{The choice of the parameters $\alpha$, $\lambda$}

At this point we already have chosen to work $m=1$ and drop 
the upper script $m$ from $h_1^m$ and $h_3^m$. This allows us to introduce
another upper script convention
\[
h_1^\l(r)=h_1(\l r),\qquad h_3^\l(r)=h_3(\l r)
\] 
which is very useful due to the key role the parameter $\l$ plays in
our analysis.

In order to understand the way a Schr\"odinger map $u$ evolves along 
the soliton family, we need to choose a pair of time dependent functions
$\alpha(t)$, $\lambda(t)$ so that \eqref{goodal} holds. Such a choice
is not unique; we will introduce here two alternatives, show that both 
are suitable and compare them.

Our main choice is analytic, and it is motivated by the equation 
\eqref{psieq}, which we want to rewrite as a linear equation with a nonlinear
{\em perturbative} term. This is not the case in \eqref{psieq},
since $A_2$ is nonzero if $\psi = 0$. Thus we want to take the bulk
part of $A_2$ and move it into the linear part of the equation.
Since $A_2$ is initialized as $h_3$ at infinity, one may try to 
take $h_3$ as the main part of $A_2$; this leads to a nonlinear
Schr\"odinger equation governed by the operator $\tilde H$, namely
\[
 (i \partial_t - \tilde H)\psi = A_0 \psi - 2 \frac{A_2-h_3}{r^2} \psi   -
\frac1{r}\Im{(\psi_2 \bar{\psi})} \psi
\]
Unfortunately, the second term on the right, though quadratic in $\psi$,
turns out to be nonperturbative on a long time scale; the difficulty is related to the lack
of time decay of $A_2-h_3$ for $r$ in a compact set. To remedy this we
instead choose $\lambda$ so that $A_2$ is close to $h_3^\lambda$ 
for $r$ in a compact set. Precisely, our full choice of parameters
is 
\begin{equation}
A_2(1,t) = h_3^{\lambda(t)}(1), \qquad  \psi_2(1,t) = i e^{ i \alpha(t)}
h_1^{\lambda(t)}(1)
%A_2(\lambda^{-1}(t),t) = 1, \qquad  \psi_2(\lambda^{-1}(t),t) = i e^{ i \alpha(t)}
\label{analdef}\end{equation}
which matches $(A_2,\psi_2)$ with $(h_3^{\lambda(t)}, i e^{ i \alpha(t)}
h_1^{\lambda(t)}(1))$ at $r  = 1$.
The matching point  is arbitrarily chosen; any other one would do.
With these parameters, the equation \eqref{psieq} takes the form
\begin{equation} \label{psieqa}
(i \partial_t - \tilde H_{\lambda(t)}) \psi  = A_0 \psi - 2
\frac{A_2-h_3^{\lambda(t)}}{r^2} \psi   - \frac1{r}\Im{(\psi_2 \bar{\psi})} \psi
\end{equation}
With this formulation we are able to track the right hand side perturbatively.
The price we pay is that the linear part now has 
a time dependent operator $H_{\lambda(t)}$, and that in addition to 
bounds for $\psi$ we also need to bootstrap the appropriate bounds 
on the parameter $\lambda$.

An alternate  choice of the parameters $\alpha$ and $\lambda$ is 
geometric:
\begin{equation}
u(1,t) = Q_{\tilde \alpha(t),\tilde \lambda(t)}(1).
\label{geomdef}\end{equation}
This choice, somewhat related to the one in \cite{gnt}, does not play any role in our analysis, 
and is given here only for comparison purposes.
As a consequence of the pointwise part of the bounds \eqref{l1c} and \eqref{l1e} we have

\begin{c1}
Assume that $\psi$ is small in $L^2$. Then both $(\alpha(t),\lambda(t))$
and $ (\tilde \alpha(t),\tilde \lambda(t))$ satisfy the condition
\eqref{goodal}. In addition, the two sets of parameters are related by
the relations
\begin{equation}
 \lambda(t) = \tilde\lambda(t), \qquad |\alpha(t)-\tilde \alpha(t)|
\lesssim \|\psi\|_{L^2}
\label{alclose}\end{equation}
\end{c1}

We remark that the first relation in \eqref{alclose} follows directly from the identity \eqref{a2u3}. The second part of \eqref{alclose}
follows from \eqref{goodal} since a direct computation 
shows that
\[
\| Q_{\alpha,\lambda} - Q_{\tilde \alpha,\tilde \lambda}\|_{\dot H^1}
\approx \frac{|\alpha - \tilde \alpha| + | \log (\lambda /\tilde \lambda)|}
{1+ |\alpha - \tilde \alpha| + | \log (\lambda /\tilde \lambda)|}
\]

\section{Spectral analysis for  the operators 
$H$, $\tilde H$; the $X,L X$ spaces}
\label{spectral}

\subsection{Spectral theory for the operator $H$}

The spectral theory for $H$ was studied in detail by
Krieger-Schlag-Tataru in \cite{KST}.  Here we simply restate the
result in \cite{KST}, in a slightly modified setup.  The modification
is threefold. Instead of working in $L^2(dr)$, we work with
$L^2(rdr)$; this is equivalent to an $r^{-\frac12}$
conjugation. Secondly, we prefer to use $\xi^2$ instead of $\xi$ as
the spectral parameter. Finally, we include the spectral measure in
the generalized eigenfunctions.

Precisely, we consider $H$ acting as an unbounded selfadjoint operator in
$L^2(rdr)$. Then $H$ is nonnegative, and its spectrum $[0,\infty)$ is  absolutely continuous. $H$ has a zero resonance, namely $\phi_0=h_1$,
\[
 H h_1 = 0.
\]
For each $\xi > 0$ one can choose a normalized generalized eigenfunction
$\phi_\xi$,
\[
 H \phi_\xi = \xi^2 \phi_\xi. 
\]
These are unique up to a $\xi$ dependent multiplicative factor,
which is chosen as described below.

To these one associates a generalized Fourier transform $\FH$ defined by
\[
 \FH{f}(\xi)=\int_0^\infty \phi_\xi(r) f(r)  rdr
\]
where the integral above is considered in the singular sense.
This is an $L^2$ isometry, and we have the inversion formula
\[
 f(r) = \int_0^\infty \phi_\xi(r) \FH{f}(\xi)  d\xi 
\]
The functions $\phi_\xi$ are smooth with respect to both $r$ and $\xi$.
To describe them one considers two distinct regions, $r \xi \lesssim 1$ 
and $r \xi \gtrsim  1$.

In the first region $r\xi \lesssim 1$ the functions $\phi_\xi$
admit a power series expansion of the form
\begin{equation} \label{repphi}
\phi_\xi (r)= q(\xi) \left( \phi_0 + \frac{1}{r} \sum_{j=1}^\infty (r\xi)^{2j}
\phi_j(r^2)\right), \qquad r\xi \lesssim 1
\end{equation}
where $\phi_0=h_1$ and the functions $\phi_j$ are analytic and satisfy 
\begin{equation} \label{derphi}
|(r \partial_r)^\alpha \phi_j| \lesssim_\alpha \frac{C^j}{(j-1)!} \log{(1+r)}
\end{equation} 
This bound is not
spelled out in \cite{KST}, but it follows directly from the integral
recurrence formula for $f_j$'s (page 578 in the paper).
The smooth positive weight $q$ satisfies 
\begin{equation}\label{qest}
 q(\xi) \approx \left\{ \begin{array}{ll}
\displaystyle \frac{1}{\xi^\frac12 |\log \xi| }  &  \xi \ll 1 \cr\cr
\xi^{\frac32}   &  \xi \gg 1
                        \end{array} \right., \qquad
|(\xi \partial_\xi)^\alpha q| \lesssim_\alpha q
\end{equation}
Defining the weight 
\begin{equation}\label{defmk1}
 m_k^1(r)=
\left\{
\begin{array}{ll}
 \min\{1, r 2^{k} \dfrac{\ln{(1+r^2)}}{\langle k \rangle}\}  & \ k < 0 \\ \\
 \min\{1, r^3 2^{3k}\},  & \ k \geq 0
\end{array}
\right.
\end{equation}
it follows that the nonresonant part of $\phi_\xi$ satisfies
\begin{equation}\label{pointphilow}
|(\xi \partial_\xi)^\alpha (r \partial_r)^\beta \left( \phi_\xi(r) - q(\xi) \phi_0(r)\right)|
\lesssim_{\alpha\beta} 2^{\frac{k}2} m_k^1(r), \qquad \xi \approx 2^k,\  r\xi \lesssim 1 
\end{equation}

In the other region $r \xi \gtrsim  1$  we begin with the functions
\begin{equation} \label{repphi+}
\phi^{+}_\xi(r)= r^{-\frac12} e^{ir\xi} \sigma(r\xi,r), 
\qquad r\xi \gtrsim  1
\end{equation}
solving 
\[
H \phi^{+}_\xi = \xi^2 \phi^+_\xi
\]
where for $\sigma$ we have the following asymptotic expansion
\[
\sigma(q,r) \approx \sum_{j=0}^\infty q^{-j} \phi^{+}_j(r), 
\qquad \phi_0^{+}=1 , \qquad \phi_1^{+}=\frac{3i}{8} + O(\frac1{1+r^2})
\]
with
\[
\sup_{r > 0} |(r\partial_r)^k \phi^{+}_j| < \infty
\]
in the following sense
\[
\sup_{r > 0} | (r \partial r)^\alpha (q \partial_q)^\beta [\sigma(q,r)-\sum_{j=0}^{j_0}
q^{-j} \phi^{+}_{j}(r) ] | \leq c_{\alpha,\beta,j_0} q^{-j_0-1}
\]
Then we have the representation
\begin{equation} \label{phipsi}
\phi_{\xi}(r)=a(\xi) \phi^{+}_\xi(r) + \overline{a(\xi) \phi^{+}_\xi(r)}
\end{equation}
where the complex valued function $a$ satisfies
\begin{equation} \label{abound}
|a(\xi)| = \sqrt{\frac2{\pi}}, \qquad | (\xi \partial_\xi)^\alpha a(\xi)| \lesssim_\alpha 1
\end{equation}

% Note: If we want a more explicit formulation of $a$ we have
% \[
% a(\xi)= \sqrt{2} a_{old}(\xi^2) \sqrt{\rho_{old}(\xi^2)}= \sqrt{\frac2{\pi}} \frac{a_{old}(\xi^2)}{|a_{old}(\xi^2)|}
% \]
% Here $old$ stands for the notation used in \cite{KST}.
% The factor of $2$ is due to the change of variables $\xi \rightarrow \xi^2$
% in defining the new Fourier transforms. 

\subsection{Spectral theory for the operator $\tilde H$} \label{spectralth}

The spectral theory for $\tilde H$ is derived from the spectral theory 
for $H$ due to the conjugate representations
\[
 H = L^* L, \qquad \tilde H = L L^*
\]
This allows us to define generalized eigenfunctions $\psi_\xi$ for
$\tilde H$ using the generalized eigenfunctions $\phi_\xi$ for
$ H$,
\[
 \psi_\xi = \xi^{-1} L \phi_\xi, \qquad L^* \psi_\xi = \xi \phi_\xi
\]
It is easy to see that $\psi_\xi$ are real, smooth, vanish at $r = 0$
and solve
\[
 \tilde H \psi_\xi = \xi^2 \psi_\xi
\]
With respect to this frame we can define the generalized Fourier transform
adapted to $\tilde H$ by 
\[
 \FtH{f}(\xi)=\int_0^\infty \psi_\xi(r) f(r)  rdr
\]
where the integral above is considered in the singular sense.
This is an $L^2$ isometry, and we have the inversion formula
\begin{equation} \label{FTL0}
 f(r) = \int_0^\infty \psi_\xi(r) \FtH{f}(\xi)  d\xi 
\end{equation}
To see this we compute, for a Schwartz function $f$:  
\[
\begin{split}
 \FtH{Lf}(\xi) & =\! \int_0^\infty  \psi_\xi(r) L f(r)  rdr
= \!\int_0^\infty L^* \psi_\xi(r)  f(r)  rdr
\\ &= \!\int_0^\infty \xi \phi_\xi(r)  f(r)  rdr = \xi \FH{f}(\xi) 
\end{split}
\]
Hence
\[
 \| \FtH{Lf}\|_{L^2}^2 = \| \xi \FH{f}(\xi)\|_{L^2}^2 = \la H f,f\ra_{L^2(rdr)}
= \|Lf\|_{L^2}^2
\]
which suffices since $Lf$ spans a dense subset of $L^2$.

The representation of $\psi_\xi$ in the two regions $r\xi \lesssim 1$
and $r\xi \gtrsim 1$ is obtained from the similar representation of $\phi_\xi$. 
In the first region $r\xi \lesssim 1$ the functions $\psi_\xi$
admit a power series expansion of the form
\[
\psi_\xi = \xi q(\xi) \left(\psi_0(r) +   \sum_{j \geq 1} (r\xi)^{2j} 
{\psi}_j(r^2)\right)
\]
where
\[
{\psi}_j(r)= ( h_3+1 +2j) \phi_{j+1}(r) + r \partial_r \phi_{j+1}(r)
\]
From \eqref{derphi}, it follows that 
\[
 |(r \partial_r)^\alpha \psi_j| \lesssim_\alpha \frac{C^j}{(j-1)!} \log{(1+r^2)}
\]
In addition, $\psi_0$ solves $L^*\psi_0 = \phi_0$  therefore
a direct computation shows that
\[
 \psi_0 = \frac1{2} \left(\frac{(1+r^2)\log(1+r^2)}{r^2}-1 \right)
\]

In particular, defining the weights 
\begin{equation}\label{defmk}
m_k(r)=
\left\{
\begin{array}{ll}
& \min\{1, \dfrac{\ln{(1+r^2)}}{\langle k \rangle}\}, \ \ \mbox{if} \ k < 0 \\ \\
& \min\{1, r^2 2^{2k}\}, \ \ \mbox{if} \ k \geq 0
\end{array}
\right.
\end{equation}
we have the pointwise bound for $\psi_\xi$ 
\begin{equation} \label{pointtp}
| (r \partial_r)^\alpha (\xi \partial_\xi)^\beta \psi_{\xi}(r) | \lesssim_{\alpha\beta}
2^{\frac{k}2} m_k(r) , \qquad \xi \approx 2^k,\  r\xi \lesssim 1 
\end{equation}

On the other hand in the regime $r \xi \gtrsim  1$ we define 
\[
 \psi^+ = \xi^{-1} L\phi^+ 
\]
and we obtain the representation
\begin{equation} \label{psirep}
\psi_{\xi}(r)=a(\xi) \psi^{+}_\xi(r) + \overline{a(\xi) \psi^{+}_\xi(r)}
\end{equation}
For $\psi^+$ we obtain the expression
\begin{equation} \label{reppsi}
\psi^{+}_\xi(r)= r^{-\frac12} e^{ir\xi} \tilde\sigma(r\xi,r), 
\qquad r\xi \gtrsim  1
\end{equation}
where $\tilde \sigma$ has the form
\[
\tilde\sigma(q,r)  = i \sigma(q,r) -\frac12 q^{-1} \sigma(q,r)
+ \frac{\partial}{\partial q} \sigma(q,r)+ \xi^{-1} L \sigma(q,r) 
\]
therefore it has exactly the same properties as $\sigma$. In particular,
for fixed $\xi$, we obtain that
\begin{equation}
\tilde{\sigma}(r\xi,r) = i -\frac78  r^{-1}\xi^{-1} + O(r^{-2})
\end{equation}

% From the properties of $\sigma$ we derive the following key properties for
%$\tilde{\psi}_\xi$:
% 
% \begin{equation} \label{tildepsi}
% |r^{\alpha} \partial_r^\alpha \tilde{\psi}_\xi(r)| + | \xi^{\alpha} \partial_\xi^\alpha
%\tilde{\psi}_\xi(r)| \lesssim_\alpha 1
% \end{equation}

\subsection{The spaces $X$ and $LX$}

So far we have measured the Schr\"odinger map $u$ in the space $\dot H^1$ 
(which correspond to $\bar u \in \dHe$),
while the differentiated field $\psi$ is in $L^2$. The operator $L$
maps $\dHe$ into $L^2$. Conversely, if for some $f \in L^2$ we solve 
\[
Lg = f
\]
then we obtain a solution $u$ which is in $\dHe$ and satisfies
\[
 \| g \|_{\dHe} \lesssim \|f\|_{L^2}
\]
However, this solution is only unique modulo a multiple of the
resonance $\phi_0$. Furthermore, in general it does not make sense
to identify $u$ by prescribing its size at infinity. The spaces
$X$ and $LX$ are in part introduced in order to remedy this 
ambiguity in the inversion of $L$.

\begin{d1} 
a) The space $X$ is defined as the completion of the subspace 
of $L^2(rdr)$ for which the following norm is finite  
\[
\| u \|_{X} = \left( \sum_{k \geq 0} 2^{2k} \| P_k^H u \|_{L^2}^2 \right)^\frac12
+ \sum_{k < 0} \frac1{|k|} \| P_k^H u \|_{L^2}
\]
where $P^H_k$ is the  Littlewood-Paley operator localizing at
frequency $\xi \approx 2^k$ in the $H$ calculus.

b) $L X$ is the space of functions of the form $f=L
u$ with $u \in X$, with norm $\| f \|_{L X}= \| u \|_{X}$.
 Expressed in the $\tilde H$ calculus, the $LX$ norm is written as
\[
\| f \|_{LX} = \left( \sum_{k \geq 0}  \| P_k^{\tilde H} f \|_{L^2}^2
\right)^\frac12 + \sum_{k < 0} \frac{2^{-k}}{|k|} \| P_k^{\tilde H} f \|_{L^2}
\]
\end{d1}

In this article we  work with equivariant Schr\"odinger maps $u$ for which
$ \| \bar u - \bar Q\|_{X} \ll 1$. This corresponds to fields $\psi$ which satisfy
$ \|\psi\|_{LX} \ll 1$.
The simplest properties of the space $X$ are summarized as follows:

\begin{p1}
 The following embeddings hold for the space $X$:   
\begin{equation} \label{Xembt}
\He \subset  X \subset \dHe
\end{equation}
In addition for $f$ in $X$ we have the following bounds:
\begin{equation} \label{pointX}
\|\langle r \rangle^\frac12 f \|_{L^\infty} \lesssim  \| f \|_{X} 
\end{equation}
\begin{equation} \label{linX}
\left\| \frac{f}{\ln (1+r)}  \right\|_{L^2} \lesssim  \| f \|_{X}
\end{equation}
\begin{equation} \label{linX4}
\left\| \langle r \rangle^\frac12{f} \right\|_{L^4} \lesssim  \| f \|_{X}
\end{equation}

\end{p1}

\begin{proof}
We first consider bounds for frequency localized functions
in the $H$ frame:

\begin{l1} \label{pk-energy}
 Assume $f_k \in L^2$ is localized at $H$-frequency $2^k$. Then
\begin{equation} \label{pointX1}
|f_k(r)| + | r \partial_r f_k(r)|\lesssim \left( \frac{2^{k^+}}{\la k^-\ra} 
 \phi_0(r) + 2^k m_k^1(r) \right) \|f_k\|_{L^2}, \qquad  r \lesssim 2^{-k}  
\end{equation}
\begin{equation} \label{pointX2}
|r^{\frac12} f_k(r)| \lesssim  2^{\frac{k}2} \|f_k\|_{L^2} , \qquad 2^{-k} \lesssim r  
\end{equation}
\begin{equation} \label{l2dX}
\| \partial_r f_k(r)\|_{L^2(A_{\geq -k})} \lesssim  2^{k} \|f_k\|_{L^2}.
\end{equation}
\end{l1}

\begin{proof}
The Fourier inversion formula gives
\[
f_k(r) = \int p_k (\xi)\FH{f}(\xi) \phi_\xi(r)  d \xi 
\]
Then the bound \eqref{pointX1} follows from 
\eqref{pointphilow} and  the Cauchy-Schwarz inequality.
Similarly \eqref{pointX2} follows from the bound 
$ |\phi_\xi|\lesssim r^{-\frac12}$ for $ r \gtrsim  2^{-k}$.
The estimate \eqref{l2dX} on the other hand follows directly from
\begin{equation}
 \|L f_k\|_{L^2} \lesssim 2^k \|f_k\|_{L^2}
\label{lfk}\end{equation}
\end{proof}

We now prove the embedding $X \subset \dHe$.
Due to the straightforward bound
$ \| Lf\|_{L^2} \lesssim \|f\|_{X}$
and the ODE estimate
\[
 \| \partial_r f\|_{L^2}+ \|r^{-1} f\|_{L^2} \lesssim  |f(1)|+\| Lf\|_{L^2},
\]
it suffices to show that $|f(1)| \lesssim \|f\|_{X}$.  But this is
obtained by direct summation from the dyadic pointwise bounds
\eqref{pointX1} and \eqref{pointX2}.

The embedding $\He\in X$ is a consequence of the bound
\begin{equation} \label{inX}
 \| f\|_{X} \lesssim \|f\|_{L^2} + \|L f\|_{L^2}
\end{equation}
The right hand side above is in effect an equivalent norm
in $\He$. To prove \eqref{inX} we use the $L^2$ norm
of $f$ for low frequencies,
\[
 \|f_{< 0}\|_{X} \lesssim \sum_{k< 0} \frac{1}{|k|} \|f_k\|_{L^2}
\lesssim \left( \sum_{k < 0}  \| f_k \|^2_{L^2(rdr)} \right)^\frac12 \left(
\sum_{k <0}  \frac1 {k^2}\right)^\frac12 \lesssim \| f_{<0}\|_{L^2(rdr)}
\]
and the $L^2$ norm of $Lf$ for high frequencies,
\[
 \| f_{\geq 0}\|_{X} \approx \| L f_{\geq 0}\|_{L^2} \lesssim \|L f\|_{L^2} 
\]
In view of the above embedding, for \eqref{pointX} it suffices to
consider $r > 1$. Then \eqref{pointX} follows by direct summation
from \eqref{pointX1} and \eqref{pointX2}.

For \eqref{linX} it also suffices
to take $r \geq 1$. The high frequencies are bounded 
directly in $L^2$, 
\[
 \|f_{\geq 0}\|_{L^2} \lesssim \|f\|_{X}
\]
so it remains to consider a single low frequency component $f_k$.
We have
\[
\left\| \frac{f_k}{\ln (1+r)} \right\|_{L^2}^2  \lesssim
\int_{0}^{2^{-k}} \frac{|f_k|^2}{|\ln (1+r)|^2} rdr + \frac{1}{k^2}
\|f_k\|_{L^2}^2 
\]
and for the first part we use \eqref{pointX1}.

Finally we prove \eqref{linX4}. For the high frequencies $f_{\geq 0}$
we interpolate between \eqref{pointX} and the $L^2$ estimate.
It remains to consider a fixed low frequency component $f_k$.
If $r \leq 2^{-k}$ then it suffices to perform a direct computation 
based on \eqref{pointX1}. If $r \geq 2^{-k}$ then we interpolate 
between \eqref{pointX2} and the trivial $L^2$ bound.

\end{proof}

Now we turn our attention to the space $LX$.
\begin{l1} If $f \in L^2$ is localized at $\tilde{H}$- frequency $2^k$ then
\begin{equation} \label{ps0}
|  f(r) | \lesssim 2^k m_k(r)(1+2^k r)^{-\frac12}  \| f \|_{L^2}
\end{equation}

\end{l1}

\begin{proof} 
This follows from  the Fourier inversion formula \eqref{FTL0}, the 
Cauchy-Schwarz inequality and \eqref{pointtp} for $r \lesssim 2^{-k}$,
respectively the bound $|\psi_\xi| \lesssim r^{-\frac12}$ for $r \gtrsim 2^{-k}$.
\end{proof}

\begin{p1}
 The following embeddings hold for $LX$:
\begin{equation}
L^1 \cap L^2 \subset LX \subset L^2
\label{LXemb}\end{equation}
\end{p1}

\begin{proof}
The second embedding is trivial. For the first one we use the $L^2$
norm for high frequencies, and it remains to use the $L^1$ norm for low
frequencies and show that
\begin{equation}\label{l1emb}
\| P_{\leq 0} f\|_{LX} \lesssim \|f\|_{L^1}
\end{equation}
It suffices to consider the case when $f$ is a Dirac mass,
$
 f = \dfrac{1}{R} \delta_{r=R}.
$
For such $f$ we bound its Fourier transform,
\[
 |\FtH f(\xi)| \lesssim \left\{ \begin{array}{ll}
\displaystyle
\frac{\xi^\frac12}{|\ln \xi|} \ln (1+R^2)    &  \xi < R^{-1} \cr\cr
R^{-\frac12}    &  \xi > R^{-1}
\end{array} \right.
\]
Thus 
\[
 \| P_{\leq 0} f\|_{LX} \lesssim 
\sum_{k < -|\log R|}    \frac{\ln (1+R^2)}{ k^2} +
\sum_{k > -|\log R|}  R^{-\frac12} \frac{2^{-\frac{k}2}}{k}  \lesssim 1
\]
and \eqref{l1emb} follows.
\end{proof}

Based on the above results we can now establish 
multiplicative properties for $X$ and $LX$:

\begin{p1} \label{Xalg}
$X$ is an algebra and the following estimates hold:
\begin{equation}\label{XLX}
\| f Lg\|_{LX} \lesssim \|f\|_{X} \|g\|_{X}
\end{equation}
\begin{equation} \label{LXstable}
 \|h_1 Lf \|_{LX} + \|h_3 Lf \|_{LX} \lesssim \|f\|_{X} 
\end{equation}
\end{p1}

\begin{proof} From \eqref{inX} is it enough to prove that if $f,g \in X$ then $f
\cdot g \in L^2$  and $L(f \cdot g) \in L^2$. From \eqref{linX4} we have
\[
 \|fg\|_{L^2} \lesssim \|f\|_{X} \|g\|_{X}
\]
therefore it remains to show that
\[
 \| L(fg)\|_{L^2} \lesssim \|f\|_{X} \|g\|_{X}
\]
We write
\[
 L(fg) =  L f \cdot g + f \cdot L g + h_3 \frac{f \cdot g}{r}
\]
For the first term (and similarly for the second one) we have
\[
 \|L f \cdot g\|_{L^2} \lesssim \| Lf\|_{L^2} \|g\|_{L^\infty}
\]
and use \eqref{pointX}. For the third term
we have
\begin{equation} \label{fginl2}
 \| r^{-1} {f \cdot g}\|_{L^2} \lesssim \|r^{-1} f\|_{L^2} \|g\|_{L^\infty}
\lesssim \|f\|_{\dHe} \|g\|_{\dHe}
\end{equation}
and the proof of the algebra property is complete.

If $f, g \in X$ then we use \eqref{linX} to estimate
\[
\| r^{-1} f \cdot g \|_{L^1} \lesssim \| f \|_X \| g\|_X 
\]
which combined with \eqref{fginl2} and  \eqref{LXemb} implies that
\begin{equation} \label{auxfgr}
\| r^{-1} f \cdot g \|_{LX} \lesssim \| r^{-1} f \cdot g \|_{L^1 \cap L^2} \lesssim \| f \|_X \| g\|_X
\end{equation}
Now we are ready to prove \eqref{XLX}.  We have
\[
\| f \cdot Lg \|_{L^2} \lesssim \| f \|_{L^\infty} \| Lg \|_{L^2} \lesssim \| f \|_X \| g \|_X
\]
and also
\[
\|P_{\geq 0} f \cdot  Lg \|_{L^1} \lesssim \| P_{\geq 0} f \|_{L^2} \| Lg \|_{L^2} \lesssim \| f \|_X \| g \|_X
\]
which places $
P_{\geq 0} f \cdot  Lg$ in $LX$ due to \eqref{LXemb}. Also
\[
P_{\leq 0} f \cdot L P_{\geq 0} g = L (P_{\leq 0} f \cdot P_{\geq 0}g) - LP_{\leq 0} f \cdot P_{\geq 0}g
+ \frac{h_3}{r} P_{\leq 0} f \cdot P_{\geq 0}g
\]
The first term belongs to $LX$ due to the algebra property of $X$, the
second term is treated as above and the third one is estimated as in
\eqref{auxfgr}.  We are then left with estimating the low frequency
contributions $ P_{\leq 0} f L P_{\leq  0}g$. Due to the $l^1$
structure of $X$ at low frequencies, this can be reduced to the case of single
frequencies, i.e. when $f$ is replaced by $f_{k}=P_{k} f$ and $g$
by $g_{j}=P_{j} g$ with $k,j < 0$.  If $k \geq j$ then
\[
\| f_{k} L g_{j} \|_{L^1} \lesssim \| f_{k} \|_{L^2} \| Lg_{j} \|_{L^2} 
\lesssim 2^{j} | k | | j | \| f_{k} \|_{X} \| g_{j} \|_X 
\lesssim  \| f_{k} \|_{X} \| g_{j} \|_X
\]
The case $k \leq j$ is similar after moving $L$ on the lower frequency
factor,
\[
f_{k} Lg_{j} = L (f_{k} \cdot g_{j}) - L f_{k} \cdot g_{j} 
+ \frac{h_3}{r} f_{k} \cdot g_{j}.
\]
We use a similar argument to prove \eqref{LXstable}.
We have
\[
h_1 Lf = L (h_1 f) - f \partial_r h_1
\]
The expression $h_1 f$  is estimated in $\He \subset X$,
while $ f \partial_r h_1$ trivially belongs to $L^1 \cap L^2$. The 
same argument applies if $h_1$ is replaced by $h_3-1$, proving
 the second estimate in\eqref{LXstable}.  
\end{proof}

\subsection{A companion space} Here we define a Sobolev type companion
$\tilde X$ for $X$ and study some simple properties for it.  This
space will be used in Section~\ref{elliptic} in order to characterize
the regularity of the Coulomb frame $(v,w)$.

We begin with the space $[\partial_r]^{-1}l^1L^2$, defined as the
completion of  of $H^1_{comp}([0,\infty))$ with respect to the
following norm
\[
\| f \|_{[\partial_r]^{-1}l^1L^2} = \| \partial_r f \|_{l^1L^2} :=\sum_{m} \| \partial_r f \|_{L^2(A_m)} 
\]
Since $\| \partial_r f \|_{L^1(dr)} \lesssim \| f \|_{[\partial_r]^{-1}l^1L^2}$,  it follows that $f$
has limits both at $0$ and $\infty$; and since it is approximated
by functions in $H^1_{comp}([0,\infty))$, it follows that $\lim_{r \rightarrow \infty} f(r)=0$.
We also have the following inequality
\begin{equation} \label{l1l2}
\| \partial_r f \|_{L^2} + \| f \|_{L^\infty} \lesssim \| f \|_{[\partial_r]^{-1}l^1L^2}
\end{equation}
Now we can define the spaces $\tilde X$ and $\partial_r \tilde{X}$,
\[
\tilde{X}= \{ f: \chi_{r \geq 1}f \in X, \chi_{r \leq 1} f \in X + [\partial_r]^{-1}l^1L^2 \}, 
\qquad
\partial_r \tilde{X} = \{ f: f=\partial_r g, g \in \tilde X \}
\]
with the induced norms. For technical purposes only we  also introduce the norm
\[
\| f \|^2_{l^2L^\infty} :=\sum_{m} \| f \|^2_{L^\infty(A_m)} 
\]
\begin{l1} \label{indX}
 The following estimates hold:
\begin{equation}
\| f\|_{\partial_r \tilde X} \lesssim \| \chi_{r \leq 1} f\|_{l^1 L^2} +  \| r \chi_{r \geq 1} f\|_{L^2} 
\label{drtXsize}\end{equation}
\begin{equation} \label{h1X}
\| h_1 f \|_{\partial_r \tilde X}  
+\| h_3 f \|_{\partial_r \tilde{X}}\lesssim \| f \|_{\partial_r \tilde X}
\end{equation}
\begin{equation} \label{XdrX}
\| f \cdot g \|_{\partial_r \tilde X} \lesssim \| f\|_{\partial_r \tilde X} \| g \|_{\tilde X}
\end{equation}
\begin{equation} \label{LXtX}
\| f \cdot g \|_{LX} \lesssim \| f \|_{LX} \| g \|_{\tilde X}
\end{equation}

\end{l1}

\begin{proof} 
{\bf Proof of \eqref{drtXsize}.}
Define $g$ by $\partial_r g = f$, $g(\infty) = 0$. We need to estimate
$g$ in $\tilde X$.  The operator $[r \partial_r]^{-1}$ is $L^2$
bounded, therefore $ \|g\|_{L^2} \lesssim \| r f\|_{L^2}$.  Hence
by \eqref{l1l2} we obtain
\[
\|  \chi_{r \leq 1} g\|_{L^\infty} + \| \partial_r  \chi_{r \leq 1} g\|_{l^1 L^2}
+ \|  \chi_{r \geq 1} g\|_{H^1} \lesssim 
\| \chi_{r \leq 1} f\|_{l^1 L^2} +  \| r \chi_{r \geq 1} f\|_{L^2} 
\]
therefore \eqref{drtXsize} follows by definition and by \eqref{Xembt}. 

{\bf Proof of \eqref{h1X}.}  For the first term we estimate
\[
\| \chi_{r \leq 1} h_1 f\|_{l^1 L^2} +  \| r  \chi_{r \geq 1} h_1 f\|_{L^2} \lesssim \|f\|_{L^2}
\lesssim \|f\|_{\partial_r \tilde X} 
\]
and conclude by \eqref{drtXsize}. A similar argument works for the second term.

{\bf Proof of \eqref{XdrX}.} We need to show that
\begin{equation} \label{XdX}
\| f \cdot \partial_r g \|_{\partial_r \tilde{X}} \lesssim \| f \|_{\tilde X} \| g \|_{\tilde X}
\end{equation}
We write $f = f_1 + f_2$, $g = g_1 + g_2$ where $f_1,g_1 \in X$
and $f_2,g_2 \in [\partial_r]^{-1} l^1L^2$ are supported in $[0,1]$.
The expression $f \partial_r g_2$ inherits the $l^1 L^2$ bound from $g_2$.
For $f_2 \partial_r g_1$ we write 
\[
f_2 \partial_r g_1 = \partial_{r} (f_2 g_1) - g_1 \partial_r f_2. 
\]
We can bound $f_2 g_1$ in $\He \subset X$ while $g_1 \partial_r f_2$
belongs to $l^1 L^2$.  For the final term we will show 
\begin{equation} \label{XddX}
\| f_1 \cdot \partial_r g_1 \|_{\partial_r \tilde {X}} \lesssim \| f_1 \|_{X} \| g_1 \|_{ X}
\end{equation}
Starting from the simpler bound
\begin{equation} \label{l2li}
\|  f_1 \|_{l^2L^\infty} \lesssim \| f_1 \|_{\dHe}
\end{equation}
we obtain
\[
\|f_1 \cdot \partial_r g_1 \|_{l^1 L^2} \lesssim \| f_1 \|_{X} \| g_1 \|_{ X}
\]
which yields an $L^\infty$ bound for $[\partial_r]^{-1} (f_1 \cdot \partial_r g_1) $ 
and suffices for $r \lesssim 1$. For larger $r$  we consider 
a dyadic decomposition as in the proof of \eqref{XLX}.  If the first factor has high
frequency then we estimate it in $L^2$ to obtain
\[
\| P_{\geq 0} f_1 \cdot \partial_r g_1 \|_{L^1} \lesssim \| f_1 \|_{X} \| g_1 \|_{ X}
\]
Combining this with the Sobolev type bound
\begin{equation}\label{sob}
\|  [\partial_r]^{-1} h\|_{L^2} \lesssim \|h\|_{L^1}
\end{equation}
we obtain 
\[
\| [\partial_r]^{-1}(P_{\geq 0} f_1 \cdot \partial_r g_1) \|_{L^2} \lesssim \| f_1 \|_{X} \| g_1 \|_{ X}
\]
therefore 
\[
\| \chi_{r \geq 1}  
[\partial_r]^{-1}(P_{\geq 0} f_1 \cdot \partial_r g_1) \|_{\He} \lesssim \| f_1 \|_{X} \| g_1 \|_{ X}
\]
which suffices by \eqref{Xembt}.

If the second factor has high frequency then we switch them
\[
f_1 \cdot \partial_r P_{\geq 0}  g_1 = \partial_r (f_1 \cdot  P_{\geq 0}  g_1)
-  \partial_r f_1 \cdot P_{\geq 0}  g_1 
\]
and use the $X$ algebra property for the first term on the  right.

It remains to consider the low frequency interactions 
$P_{k} f_1 \cdot \partial_r P_{j}  g_1$ with $k,j < 0$. Assuming $k \leq j$,
by using \eqref{pointX1} for small $r$ and \eqref{l2dX} for large $r$ we obtain
\[
\| P_{k} f_1 \cdot \partial_r P_{j}  g_1\|_{L^1} \lesssim  \| P_k f_1 \|_{X} \| P_j g_1 \|_{ X}
\]
and conclude as before. On the other hand if $k > j$ we switch the derivative 
to the first factor and use again  the $X$ algebra property.

{\bf Proof of \eqref{LXtX}.}
We split 
\[
f g = f \chi_{r \geq 1} g + f \chi_{r \leq 1} g 
\]
For the first term we use \eqref{XLX}. The second has compact support, therefore
it suffices to estimate it in $L^2$ and use \eqref{LXemb}.
\end{proof}

\subsection{Littlewood-Paley projectors 
in the $\tilde H$ frame}

% For a given smooth dyadic partition of unity 
% \[
% 1 = \sum_{k \in \Z} \chi_k(\xi)
% \]
% with $\chi_k$ supported in  $\{\xi \approx 2^k\}$, 
% we define the Littlewood-Paley projectors 
% \[
% P_k f(r) = \int \hat{f}(\xi) \chi_k(\xi) \psi_\xi(r) d\xi.
% \]
The first aim of the following proposition is to characterize the
kernels $K_k(r,s)$ of the projectors $P_k$ in the $\tilde H$ frame.
Secondly, we consider the kernels $K_k^1(r,s)$ of the operators
$L^{-1} P_k$, which can be defined as
\[
 L^{-1} P_k := L^* \tilde H^{-1} P_k.
\]
We remark that the adjoint operators are given by
\[
( L^{-1} P_k)^* = P_k (L^*)^{-1} := P_k \tilde H^{-1} L.
\]

\begin{p1} \label{p-lp}
a) The kernel $K_k(r,s)$ of $P_k$ satisfies the bounds
\begin{equation}
|K_k(r,s)| \lesssim \frac{2^{2k} m_k(r) m_k(s)}{(1+2^k|r-s|)^{N}(1+2^k(r+s))}, 
\label{ker1}\end{equation}
\begin{equation}
| \partial_r K_k(r,s)| \lesssim \frac{(2^{3k}+ 2^{2k}r^{-1}) m_k(r) m_k(s)} {(1+2^k|r-s|)^{N}(1+2^k(r+s))}.
\label{ker2}\end{equation}

b) If $k \geq 0$ then the kernel  $K^1_k(r,s)$ of $L^{-1} P_k$ 
 satisfies the bound
\begin{equation}
|K^1_k(r,s)| \lesssim \frac{2^k m_k^1(r) m_k(s)}{ (1+2^k|r-s|)^{N}(1+2^k(r+s))}
\label{ker3}\end{equation}
If $k < 0$ then $K^1_k(r,s)$ admits a decomposition
\[
 K_k^1 = K_{k,reg}^1 + K_{k,res}^1
\]
where the regular part $K_{k,reg}^1$ satisfies \eqref{ker3}
and the resonant part $K_{k,res}^1$ has the form
\begin{equation}
 K_{k,res}^1(r,s) = h_1 (r) \chi_{2^k r \leq 1} c_k (s),
\quad |c_k (s)| \lesssim |k|^{-1} m_k(s) (1+2^k s)^{-N} 
\label{ker4}\end{equation}
and $\chi_{2^k r \leq 1}$ is a smooth bump function 
supported in $\{2^k r  \lesssim 1\}$ which equals $1$ in 
$\{2^k r  \ll 1\}$.
\end{p1}
\begin{proof}
a)  We denote by $\chi_k$ the symbol of $P_k$. This is a
smooth bump supported at $\xi \approx 2^k$, which is all that we use 
in the proof.
The kernel $K_k(r,s)$ is symmetric and has the form
\[
K_k(r,s) = \int \psi_{\xi}(r) \psi_{\xi}(s) \chi_k(\xi) d\xi.
\]
If $r,s \gtrsim  2^{-k}$ then we use the representation \eqref{psirep} to obtain
\[
K_k(r,s)=\!\! \sum_{j,l=0,1} \! 
\int \! r^{-\frac12} s^{-\frac12} e^{i((-1)^{j}r + (-1)^{l} s)\xi} a_j(\xi) a_l(\xi) 
\tilde{\sigma}_j(r\xi,r) \tilde{\sigma}_l(r\xi,r) \chi_k(\xi) d\xi
\]
where $a_0=a, a_1=\overline{a}, \tilde{\sigma}_0=\tilde{\sigma},
\tilde{\sigma}_1=\overline{\tilde{\sigma}}$.  Using stationary phase
together with the bounds on $a$ and the characterization of
$\tilde{\sigma}$ gives the bound in \eqref{ker1}. We note that the
stationary phase brings decay factors of type $(1+2^k|r-s|)^{-N}$.

We now consider the case $r \gtrsim  2^{-k}, s \lesssim 2^{-k}$
(and also, by symmetry, the case $r \lesssim 2^{-k}, s \gtrsim  2^{-k}$).
Then 
\[
K_k(r,s)= \sum_{j=0,1} \int  \psi_\xi(s) r^{-\frac12} e^{i(-1)^{j}r\xi} a_j(\xi) 
\tilde{\sigma}_j(r\xi,r)  \chi_k(\xi) d\xi
\]
The first factor $\psi_\xi(s)$ is smooth in $\xi$ on the dyadic scale;
precisely, we have the pointwise bound \eqref{pointtp}.  Then we use
stationary phase, \eqref{pointtp}, the bounds on $a$ and the
characterization of $\sigma$ to claim \eqref{ker1}.

Finally, if $r,s \lesssim 2^{-k}$, the arguments for \eqref{ker1} and
\eqref{ker2} follow directly from the pointwise bounds \eqref{pointtp}
on $\psi_\xi$.

For the estimate \eqref{ker2} we write $\partial_r = - L^* + \frac{h_3-1}r$.
Then 
\[
 | \partial_r K_k(r,s)| \lesssim | L^* K_k(r,s)| + \frac{1}{r} |K_k(r,s)|
\]
and $L^* K_k(r,s)$ is of the form $L^* K_k(r,s)= 2^{2k} K^1_k(r,s)$
with $K^1_k$ as in part (b). Hence it suffices to prove part (b)
of the proposition.

b) Since $L^* \psi_\xi = \xi \phi_\xi$, the kernel $K_k^1$ is given by 
\[
K_k^{1}(r,s) = \int \xi^{-1} \phi_{\xi}(r) \psi_{\xi}(s) \chi_k(\xi) d\xi
\]
If $r \gtrsim 2^{-k}$ then $\phi_{\xi}(r)$ is similar to $\psi_\xi(r)$
and $K_k^1$ satisfies the same bounds as $K_k$ with an additional $2^{-k}$
factor.
If $r \lesssim 2^{-k}$ and $k \geq 0$ then $\phi_\xi$ is smooth 
on the dyadic scale and has size $r \xi^{\frac32}$ therefore 
we can argue again as in case (a).
Finally if $r \lesssim 2^{-k}$ and $k < 0$ then we decompose 
$\phi_\xi$ according to \eqref{repphi} into 
\[
 \phi_\xi(r) = q(\xi) \phi_0(r) + \phi_\xi^{err}(r)
\]
where $\phi_\xi^{err}$ is smooth on the dyadic scale and has size
$q(\xi) r \xi^2 \log(1+r)$. The first term yields the resonant 
component $K^1_{k,res}$ and the second term gives the regular 
component $K^1_{k,reg}$.
\end{proof}

\subsection{ Time dependent frames and
 the transference identity}

Later in the article we need to work with a time dependent parameter 
$\lambda$, and thus with a time dependent Fourier transform 
associated to the operator $\tilde H_\lambda$. By rescaling, its normalized
generalized eigenfunctions are
\[
 \psi^\lambda_\xi(r) = \lambda^{-\frac12} \psi_{\lambda \xi}(\lambda^{-1}r),
\qquad \tilde H_\lambda \psi^\lambda_\xi = \xi^2 
\psi^\lambda_\xi 
\]
We denote the associated Fourier transform by $\F_\lambda$; this is 
an $L^2$ isometry. To study its $\lambda$ dependence we use
the transference operator $\K_\lambda$, previously introduced 
and studied in \cite{KST}:
\[
 \K_\lambda = \F_\lambda \frac{d}{d\lambda}   \F_\lambda^*
\]
By scaling it suffices to analyze the operator $\K = \K_1$.

\begin{p1} [\cite{KST}]\label{p:transference}
The operator $\K$ is a skew-adjoint Hilbert transform type operator,
whose kernel $K(\xi,\eta)$ has the form
\[
K(\xi,\eta) = \text{p.v.} \frac{F(\xi,\eta)}{\xi^2 -\eta^2}, \qquad F(\xi,\eta) = \langle \frac{1}{(1+r^2)^2} \psi_\xi, \psi_\eta \rangle
\]
where the symmetric function $F$ satisfies the following bounds
\begin{align*}
|(\xi \partial_\xi)^\alpha (\eta \partial_\eta)^\beta (\xi \partial_\xi + \eta \partial_\eta)^\sigma F(\xi,\eta)| &\lesssim 
\dfrac{\xi^\frac12 \eta^\frac12}{\la \ln \xi \ra  \la \ln \eta\ra}
&  \xi,\eta \lesssim 1 \cr
|(\xi \partial_\xi)^\alpha \partial_\eta^\beta (\xi \partial_\xi + \eta \partial_\eta)^\sigma F(\xi,\eta)| & \lesssim 
\dfrac{\xi^\frac12}{\la \ln \xi \ra (1+\eta)^N}
&  \xi \lesssim 1 \lesssim \eta \cr
 |\partial_\xi^\alpha \partial_\eta^\beta (\xi \partial_\xi +\eta \partial_\eta)^\sigma F(\xi,\eta)| &\lesssim 
\dfrac{1}{(1+|\xi-\eta|)^{N}}  & 1 \lesssim \xi,\eta
\end{align*}
where $\sigma,N \in \N$ and $\alpha, \beta \in \N $ if $|\xi -\eta| \gtrsim  \max( \xi,\eta )$ and 
$\alpha + \beta \leq 2$ if $|\xi -\eta| \ll \max( \xi,\eta )$.
\end{p1}

\begin{proof}
  We merely outline the computation, as a complete proof is given
  in \cite{KST}. Furthermore, the proof is similar to the proof of the 
next proposition, which is presented in full. Formally the kernel $K$
is given by
\[
K(\xi,\eta) = \langle \psi_{ \xi}, 
\left(\frac{d}{d\lambda} \psi_{\eta}^\lambda\right)_{|\lambda = 1}\rangle
\]
We have
\[
\left(\frac{d}{d\lambda} \psi_{\eta}^\lambda\right)_{|\lambda = 1} =
\frac{d}{d\lambda} \left(\lambda^{-1/2} \psi_{\lambda
    \eta}(\lambda^{-1}r)\right)_{|\lambda = 1} = (\eta \partial_\eta -
r \partial_r -\frac12) \psi_{\eta}
\]
therefore
\[
K(\xi,\eta) = \langle \psi_{ \xi}, 
 (\eta \partial_\eta - r \partial_r -\frac12) \psi_{\eta} \rangle
\]
Hence we obtain:
\[
\begin{split}
  (\xi^2 -\eta^2) K(\xi,\eta) = & \langle \xi^2 \psi_{ \xi},
  (\eta \partial_\eta - r \partial_r -\frac12) \psi_{\eta} \rangle
  - \langle \psi_{ \xi}, \eta^2 (\eta \partial_\eta - r \partial_r -\frac12) \psi_{\eta} \rangle \\
  = & \langle \tilde{H} \psi_{ \xi}, (\eta \partial_\eta -
  r \partial_r -\frac12) \psi_{\eta} \rangle
  - \langle \psi_{ \xi}, (\eta \partial_\eta - r \partial_r -\frac12) \eta^2 \psi_{\eta} \rangle 
  +  2 \langle \psi_{ \xi}, \eta^2 \psi_{\eta} \rangle \\
  = &- \langle  \psi_{ \xi},  [\tilde{H}, r \partial_r] \psi_{\eta} \rangle + 2 \langle \psi_{ \xi}, \eta^2 \psi_{\eta} \rangle \\
  = &- \langle \psi_\xi, \frac{8}{(1+r^2)^2} \psi_\eta \rangle
\end{split}
\]
The bounds on $F$ are derived from the representations for $\psi_\xi$ 
given by \eqref{pointtp}-\eqref{reppsi}. 
\end{proof}

\noindent
Next we consider the related problem of comparing the Fourier
transforms in nearby frames.

\begin{p1} \label{p:transfa}
Suppose that $|\lambda_1-1| \ll 1$ and $|\lambda_2-1| \ll 1$.
Then the kernel of the operator $\F_{\lambda_1} \F^*_{\lambda_2}$
has the form
\[
K_{\lambda_1 \lambda_2}(\xi,\eta) = a(\lambda_1,\lambda_2,\xi) \delta_{\xi=\eta} + \text{p.v.} \frac{2
b(\lambda_1,\lambda_2,\xi,\eta) \xi^\frac12 \eta^{\frac12}}{\xi^2-\eta^2}
\]
where $a$ and $b$ are smooth functions in all variables  satisfying 

\[
a^2(\lambda_1,\lambda_2,\xi)+b^2(\lambda_1,\lambda_2,\xi,\xi)= 1
\]
and the following size and regularity:
\begin{align}
|\partial_{\lambda_{12}}^\alpha (\xi \partial_{\xi})^\beta
(\eta \partial_\eta)^\gamma (\xi \partial_\xi + \eta \partial_\eta)^\sigma b| 
& \lesssim \frac{1}{\la \log \xi \ra \la \log \eta \ra} 
&  \xi,\eta \lesssim 1    \label{be1}
\\
|\partial_{\lambda_{12}}^\alpha (\xi \partial_{\xi})^\beta
\partial_\eta^\gamma (\xi \partial_\xi + \eta \partial_\eta)^\sigma b| 
& \lesssim \frac{1}{\la \log \xi \ra (1+\eta)^N} 
&  \xi \lesssim 1 \lesssim \eta \label{be2}
\\
|\partial_{\lambda_{12}}^\alpha \partial_{\xi}^\beta
\partial_\eta^\gamma (\xi \partial_\xi + \eta \partial_\eta)^\sigma 
b| 
& \lesssim \frac{1}{\xi^\frac12 \eta^\frac12 (1+|\eta-\xi|)^N} 
&  1 \lesssim \xi,\eta    \label{be3}
\end{align}
where $\alpha, \sigma,N \in \N$ and $\beta, \gamma \in \N $ if $|\xi -\eta| \gtrsim  \max( \xi,\eta )$  
and $\beta + \gamma \leq 2$ if $|\xi -\eta| \ll \max( \xi,\eta )$.
%(\lambda_1,\lambda_2,\xi,\eta)
\end{p1}

\begin{proof}
  Given a smooth radial function $f$ in $(0,\infty)$ which is compactly
  supported away from zero we have the following integral
  representation for $\F_{\lambda_1} \F^*_{\lambda_2} f$:
\begin{equation}
\begin{split}
\F_{\lambda_1} \F^*_{\lambda_2} f (\xi) & = \int_0^\infty \int_0^\infty \psi_\xi^{\l_1}(r) \psi_{\eta}^{\l_2}(r) f(\eta) d \eta rdr \\
& = \lim_{R \to \infty} 
\int_0^\infty \int_0^\infty \chi(r/R) \psi_\xi^{\l_1}(r) \psi_{\eta}^{\l_2}(r) f(\eta) d \eta rdr
\end{split} \label{kerrep}
\end{equation}
where $\chi$ is a smooth radial compactly supported bump function 
which equals $1$ in the unit ball.
Then the off-diagonal part of $K_{\l_1 \l_2}(\xi,\eta)$ is given by
\[
K_{\l_1 \l_2}(\xi,\eta) = \lim_{R \to \infty}
\int_0^\infty  \chi(r/R) \psi_\xi^{\l_1}(r) \psi_{\eta}^{\l_2}(r) rdr := 
\langle \psi_\xi^{\l_1},  \psi_{\eta}^{\l_2} \rangle
\]
This is meaningful if the above limit exists uniformly on compact sets
off the diagonal; that is always the case due to the asymptotics
for $\psi_\xi$ as $r \to \infty$ in \eqref{psirep}, \eqref{reppsi}.
Multiplying the previous relation by $(\xi^2-\eta^2)$ and integrating
by parts gives
\[
\begin{split}
(\xi^2-\eta^2) \langle \psi_\xi^{\l_1},  \psi_{\eta}^{\l_2} \rangle 
&= \langle \xi^2 \psi_\xi^{\l_1},  \psi_{\eta}^{\l_2} \rangle - \langle \psi_\xi^{\l_1},  \eta^2 \psi_{\eta}^{\l_2} \rangle 
 = \langle \tilde{H}_{\l_1} \psi_\xi^{\l_1},  \psi_{\eta}^{\l_2} \rangle - \langle \psi_\xi^{\l_1},  \tilde{H}_{\l_2} \psi_{\eta}^{\l_2} \rangle \\
& = \langle \psi_\xi^{\l_1},  (\tilde{H}_{\l_1} - \tilde{H}_{\l_2}) \psi_{\eta}^{\l_2} \rangle 
 = \langle \psi_\xi^{\l_1},  (\tilde{V}_{\l_1} - \tilde{V}_{\l_2}) \psi_{\eta}^{\l_2} \rangle \\
& =  \langle \psi_\xi^{\l_1},  \frac{4(\l_2^2-\l_1^2)}{(1+\l_1^2 r^2)(1+\l_2^2 r^2)} \psi_{\eta}^{\l_2} \rangle
\end{split}
\]
which leads to the following formula for $b$:
\[
b(\l_1,\l_2,\xi,\eta) =  \langle \psi_\xi^{\l_1},  \frac{4(\l_2^2-\l_1^2)}{(1+\l_1^2 r^2)(1+\l_2^2 r^2)} \psi_{\eta}^{\l_2} \rangle
\]
We note that the above computation should be done with the cutoff
$\chi(r/R)$ included, and then pass to the limit $R \to \infty$;
This computation is tedious but routine, so we omit it.  The
bounds \eqref{be1}-\eqref{be3} are obtained from this formula using
again the representation \eqref{pointtp}, \eqref{psirep} and
\eqref{reppsi} for the functions $\psi_\xi$.

Next we identify the behavior of the kernel $K_{\lambda_1\lambda_2}$
near the diagonal by using the representation in \eqref{psirep} and
\eqref{reppsi} for $r \xi \gtrsim  1$,
\[
\psi_\xi^\l(r) = \Re{\left( r^{-\frac12} e^{ir\xi} 
(i - \frac1{8 r \xi}) a(\l \xi)\right)}  + O(r^{-\frac{5}2})
\]
Since we have already identified the off-diagonal kernel of 
$K_{\lambda_1\lambda_2}$, for this purpose we can freely neglect 
any part of $K_{\lambda_1\lambda_2}$ which has a locally bounded kernel.
For large $r$ we have
\[
\begin{split}
\psi_\xi^{\l_1}(r) \psi_{\eta}^{\l_2}(r) 
& = \frac{1}r \Re{\left( e^{ir\xi} 
(i - \frac7{8 r \xi}) a(\l_1 \xi)\right)}
\Re{\left(  e^{ir\eta} 
(i - \frac7{8 r \eta}) a(\l_2 \eta)\right)} + O(\frac{1}{r^3})
\\
& = I_{res} - I_{nr} + II_{res} - II_{nr} + O(\frac{1}{r^3})
\end{split}
\]
 where
\[
I_{res} = \frac{1}{2r} \Re \left( a(\l_1 \xi)  \bar a(\l_2 \eta) e^{ir(\xi-\eta)}\right), \qquad 
I_{nr} = \frac{1}{2r} \Re \left( a(\l_1 \xi)  a(\l_2 \eta) e^{ir(\xi+\eta)}\right)
\]
\[
II_{res} =\frac{7}{16 r^2} \frac{\xi-\eta}{\xi \eta} \Im \left( a(\l_1 \xi)  
\bar a(\l_2 \eta) e^{ir(\xi-\eta)}\right),
\]
\[ 
II_{nr} =\frac{7}{16 r^2} \frac{\xi+\eta}{\xi \eta} \Im  \left( a(\l_1 \xi)  a(\l_2 \eta) e^{ir(\xi+\eta)}\right)
\]
Hence returning to  \eqref{kerrep}, for $\xi$ in a compact set
we have
\[
\F_{\lambda_1} \F^*_{\lambda_2} f (\xi) = 
\int_{0}^\infty\int_{0}^\infty (I_{res} - I_{nr} + II_{res} - II_{nr}) f(\eta) d\eta rdr
+ O(\|f\|_{L^1})
\]
In the nonresonant terms $I_{nr}$ and $II_{nr}$ the phase is 
uniformly oscillatory, so integration by parts in $r$ allows 
for a gain of arbitrarily many powers of $r^{-1}$.

In the second resonant term $II_{res}$ the phase may be stationary. However,
the factor of $\xi-\eta$ allows for one integration by parts 
in $r$ which gain an $r^{-1}$ factor, sufficient to insure absolute 
convergence in the integral. Thus we are left with 
\[
\begin{split}
\F_{\lambda_1} \F^*_{\lambda_2} f (\xi) & = 
\int_{0}^\infty \int_{0}^\infty I_{res} f(\eta) d\eta rdr
+ O(\|f\|_{L^1})
\\ 
& = \int_{0}^\infty \int_0^\infty  
\frac12 \Re \left( a(\l_1 \xi)  \bar a(\l_2 \eta) e^{ir(\xi-\eta)}\right)
f(\eta) d\eta dr+ O(\|f\|_{L^1})
\end{split}
\]
Using elementary properties of the Fourier transform and the notation $\bf H$
for the Heaviside function,
the last integral is expressed in the form
\[
\begin{split}
\frac14 \Re \iint_{\R^2} 
 & e^{ir(\xi-\eta)} a(\l_1 \xi) \overline{a}(\l_2 \eta)(1+ {\bf H}(r)) f(\eta) d\eta d r  
\\
 & =  \frac{\pi}2 \Re{(a(\l_1 \xi) \overline{a}(\l_2 \xi))}  
f(\xi)  +  \frac{\pi}2 
p.v. \int \Im{(a(\l_1 \xi) \overline{a}(\l_2 \eta))} \frac{1}{\xi-\eta} f(\eta) d\eta 
\end{split}
\]
Hence for $\xi,\eta$ in a compact set we obtain
\[
K_{\lambda_1\lambda_2}(\xi,\eta) = 
 \frac{\pi}2 \Re{(a(\l_1 \xi) \overline{a}(\l_2 \xi))} \delta(\xi-\eta)
+\frac{\pi}2 \Im{(a(\l_1 \xi) \overline{a}(\l_2 \xi))} p.v. \frac{1}{\xi-\eta}
+ O(1)
\]
Comparing this with the off-diagonal representation of $K_{\lambda_1\lambda_2}$,
we obtain the representation in the proposition 
with 
\[
a(\lambda_1,\lambda_2,\xi) = \frac{\pi}2 \Re{(a(\l_1 \xi) \overline{a}(\l_2 \xi))}, 
\qquad b(\lambda_1,\lambda_2,\xi,\xi) = \frac{\pi}2 \Im{(a(\l_1 \xi) \overline{a}(\l_2 \xi))}
 \]
Using \eqref{abound}, it then follows that on the diagonal we have
\[
a^2(\lambda_1,\lambda_2,\xi)+b^2(\lambda_1,\lambda_2,\xi,\xi)= |\frac{\pi}2 a(\l_1 \xi) \overline{a}(\l_2 \eta) |^2=1.
\]
\vspace{-.4in}

 \end{proof}

\subsection{ Compositions of Littlewood-Paley 
projectors}

We first consider dyadic bump functions in the Fourier space, and we estimate
their inverse Fourier transforms:

\begin{p1}\label{p:bumpft}
Let $k \in \Z$ and  $\chi_k$ be a unit size bump function supported
in the $\{\xi \approx 2^k\}$ dyadic region. Then for $|\lambda -1|\ll 1$
its inverse Fourier transform satisfies the bounds 
\begin{equation}
|\partial_\lambda^\alpha (r \partial_r)^\beta (\F^{*}_\lambda \chi_k)(r)| 
\lesssim_{\alpha,\beta} 2^{\frac{3k}2} m_k(r)  (1+2^k r)^{-N} 
\end{equation}\label{bumpft}
\end{p1}
\begin{proof} We have
\[
\mathcal{F}^{*}_\lambda \chi_k (r) = \int \psi_\xi^\lambda(r) \chi_{k}(\xi) d\xi
\]
If $r \lesssim 2^{-k}$, then using the pointwise estimate \eqref{pointtp}
gives \eqref{bumpft}.  If $r \geq 2^{-k}$, then using \eqref{psirep},
\eqref{reppsi} and stationary phase gives \eqref{bumpft}.
\end{proof}

The next step is to consider the composition of two dyadically
separated projectors associated with different frames.

\begin{p1}\label{p:prodtwo}
Let $j,k \in \Z$ with $|j-k| \gg 1$, and $\lambda$ in a compact subset of 
$(0,\infty)$. Then the kernels $K_{jk}(r,s)$ of $P_j P_k^\lambda$
satisfy the bounds
\begin{equation}\label{kerbiprod}
|K_{jk}(r,s)|  \lesssim \frac{2^{j+k - |j-k|-N(j^++k^+)}}
{\la j^-\ra \la k^- \ra}
\frac{m_j(r) m_k(s)}{(1+2^j r)^N (1+ 2^k s)^N} 
\end{equation}
\end{p1}

\begin{proof}
 The kernel $K_{jk}(r,s)$ of $
P_{j} \tilde P_k^\lambda$ is given by
\[
K_{jk}(r,s) = \int \psi_\xi(r) \chi_{j}(\xi) K_{1\l}(\xi,\eta)
\psi_{\eta}^\l(s) \chi_{k}(\eta) d \mu d\eta
\]
By Proposition \ref{p:transfa}, the symbol $ \chi_{j}(\xi)
K_{1\l}(\xi,\eta) \chi_{k}(\eta)$ is smooth on the dyadic scales and
has size
\[
| \chi_{j}(\xi) K_{1\l}(\xi,\eta) \chi_{k}(\eta)| \lesssim
\frac{1}{\la k^-\ra \la j^-\ra} 2^{-\frac{k+j}2 - |k-j|}
2^{-N(k^++j^+)}
\]
Hence \eqref{kerbiprod} follows from Proposition \ref{p:bumpft}.
\end{proof}

Given a frequency localized function in one frame, the above
proposition allows us to relocalize it in a different frame with good
pointwise error bounds. For this we consider a projector $\tilde P_k$ 
whose symbol $\tilde \chi_k$ equals $1$ within the support of $\chi_k$.
\begin{c1} \label{reloc}
Let $\psi \in L^2$ and $k \in \Z$. Then we have
\begin{equation}\label{relocest}
 P_k^\lambda \psi = \tilde P_k P_k^\lambda \psi+ \psi_{k}^{err}, \qquad | \psi_{k}^{err}| \lesssim  \frac{2^{k^-} 2^{-Nk^+} \la (\log (2+r) + k)^- \ra}{\la k^- \ra^2  (1+2^{k^-} r)^2} 
\|P_k^\lambda \psi\|_{L^2} 
\end{equation}
\end{c1}
\begin{proof}
We write
\[
  (1-\tilde P_k) P_k^\lambda = \sum_{|j-k| \gg 1} P_j \tilde P_k^\lambda  
P_k^\lambda 
\]
A direct estimate using \eqref{kerbiprod} gives
\[
|  P_j   P_k^\lambda \psi(r)|
 \lesssim \frac{2^{j - |j-k|-N(j^++k^+)}}
{\la j^-\ra \la k^- \ra}
\frac{m_j(r)}{(1+2^j r)^N} \|P_k^\lambda  \psi\|_{L^2} 
\]
It remains to evaluate the sum with respect to $j$ of the coefficients 
on the right. Because of the rapid decay for positive $j,k$, 
it suffices to assume that $k < 0$ and restrict the sum to $j < 0$.
Then we are left with the sum
\[
\frac{2^k}{|k|}  \sum_{j < 0} \log(1+r) \frac{2^{j-k - |j-k|} }
{j^2 (1+2^j r)^N} 
\]
There are two thresholds for $j$ in this sum, namely $-\log(1+r)$ and $k$.
If $2^k r > 1$ then the sum is given by the summand at $j = -\log(1+r)$.
Else, the sum is bounded by
\[
 \frac{2^k}{|k|}  \sum_{j = k}^{-\log(1+r)} \frac{\log(1+r)}{j^2} \frac1{(1+2^j r)^N}
\]
The bound \eqref{relocest} easily follows.
\end{proof}

Finally, we consider the product of three projectors:
\begin{p1}\label{p:prodthree}
Let $j,h,k \in \Z$  and $\lambda$ in a compact subset of 
$(0,\infty)$. 

a) Assume that $|j-k| \gg 1$ and $|h-k| \gg 1$. 
Then the kernels $K_{jkh}(r,s)$ of $P_j P_k^\lambda P_h$
can be represented as the sum of a rapidly
convergent series of terms $K^l_{jkh}(r,s)$ of the form
\begin{equation}\label{kertriproda}
\begin{split}
K^l_{jkh}(r,s) = &\  c_{jkh}\  l^{-N}   g^{l}(\lambda) \phi_j(r) \phi_h(s),  \\
c_{jkh} =& \  
\dfrac{ 2^{ - |j-k| - |k-h|-N(j^++k^++h^+)}}{\la j \ra \la k
      \ra^2 \la h \ra}\\
|\phi_j(r)| \leq &\  2^j (1+2^j r)^{-N}, \quad |\phi_h(s)| \leq  2^h (1+2^h s)^{-N}.
\end{split}
\end{equation}
with $g^l$ uniformly bounded in $C^N$.

b) Assume that either $|j-k| \lesssim 1$ and $|h-k| \gg 1$ or
$|j-k| \gg 1$ and $|h-k| \lesssim 1$. 
Then the kernels $K_{jkh}(r,s)$ of $P_j P_{k}^\l P_h$
can be represented as above but with 
\begin{equation}\label{kertriprodb}
c_{jkh} = \dfrac{ 2^{ - |j-h|-N(j^++h^+)}}{\la j \ra \la h \ra}
\end{equation}

c) Assume that  $|j-h| \lesssim 1$ and $|h-k| \lesssim 1$.
Then the operators $P_j P_{k}^\l P_h$ can be represented as 
 sum of a rapidly convergent series of the form
 \begin{equation}\label{kertriprodc}
P_j P_{k}^\l P_h=  P_j P_{k} P_h +  c_{jkh} \sum  l^{-N} 
  g^l(\lambda)  Q^l_{jkh}, \qquad c_{jkh}= \frac{1}{2^{k^+} \la k^-\ra^2} 
\end{equation}
with  $g^l$ uniformly bounded in $C^N$ and $ \|Q^l_{jkh}\|_{L^2 \to L^2} \leq 1$. 
\end{p1}

\begin{proof}
  a) We use Proposition \ref{p:transfa} to estimate the Fourier kernel
  $\hat K_{jkh}$ of $P_{j} \tilde{P}_{k}^\l \tilde{P}_{h}$, given by
\[
\hat K_{jkh}(\xi,\zeta) = \int \chi_{j}(\xi) K_{1\l}(\xi,\eta)
\chi_{k}(\eta) K_{\l1}(\eta,\zeta) \chi_{h}(\zeta) d\eta
\]
It follows that $\hat K_{jkh}$ is smooth in $\xi,\zeta$ on the dyadic
scale, smooth in $\lambda$ and has size
\begin{equation} \label{FB}
    |\hat K_{jkh}(\xi,\zeta)|  \lesssim 2^{-\frac{j+h}2}  c_{jkh}
\end{equation}
Separating variables, it suffices to consider  
kernels $\hat K_{jkh}$ of the form
\[
\hat K_{jkh}(\xi,\zeta) = \sum_l  2^{-\frac{j+h}2}
c_{jkh} l^{-N} g^l (\lambda) \chi^l_{j}(\xi) \chi^l_h(\eta)
\]
with $ \chi^l_{j}$, $\chi^l_h$ smooth dyadic bump functions, and 
$g$ smooth. Then the conclusion follows using 
the bounds for the inverse Fourier transforms of 
 $\chi^l_{j}$ and $ \chi^l_h$  given by Proposition~\ref{p:bumpft}.

b) The proof is similar to the one in case (a), with the only difference
that we need to consider the contribution of the diagonal term
in exactly one of the kernels $ K_{1\l}(\xi,\eta)$ and $ K_{\l1}(\eta,\zeta) $.
Here we take advantage of the factor $(\xi \partial_\xi + \eta \partial_\eta)^\sigma$ with
$\sigma \in \N$ in \eqref{be1}-\eqref{be3} in order to claim that if $s$ is smooth in $\eta$ then
\[
\int \frac{b(\l_1,\l_2,\xi,\eta)}{\xi^2-\eta^2} s(\eta) d\eta
\]
is smooth in $\xi$. 

c) In this case we need to allow near diagonal contributions 
from both kernels $ K_{1\l}(\xi,\eta)$ and $ K_{\l1}(\eta,\zeta) $.
For each of them we can use
Proposition~\ref{p:bumpft} to write
\[
 K_{1\l}(\xi,\eta)  =   \delta_{\xi=\eta} 
+\left(1-a(1,\lambda,\xi,\eta)\right) \delta_{\xi=\eta} 
+ p.v. \frac{2b(1,\lambda,\xi,\eta)\xi^\frac12 \eta^\frac12}{\xi^2-\eta^2} 
\]
In the region $\xi,\eta,\zeta \approx 2^k$ the functions $1-a$ and $b$ 
are smooth in the $\l$ variable and have size $\la k^- \ra^{-2} 2^{-k^+}$.
Hence separating the variable $\lambda$ and performing the 
remaining compositions we arrive at the desired conclusion.
  
\end{proof}

\subsection{Nonresonant quadrilinear forms}
\label{sec:g}
Here we prove bounds for quadrilinear expressions in nonresonant
situations. Precisely, we consider four dyadic frequencies 
\[
j   <  k_3 < k_2 = k_1, \qquad j < 0
\]
and corresponding frequencies $\xi \approx 2^j$, $\xi_l \approx 2^{k_l}, l \in \{1,2,3\}$
which are subject to one of the two additional conditions:

i) $k_2-k_3 \gg 1$ and $|\xi_1 - \xi_2| \ll 2^{k_3}$.

ii) $|k_3 - k_2| \lesssim 1$ and $|\xi_1^2 + \xi_2^2 - \xi_3^2| \ll  2^{2k_3}$. \\
In  both cases $\xi_1$ and $\xi_2$ may be close but $\xi_3$ is dyadically
separated from them.
To the quadruplet of generalized eigenfunctions $\psi_\xi$, $\psi_{\xi_1}$
$\psi_{\xi_2}$ and $\psi_{\xi_3}$ we associate the quadrilinear
expressions:
\[
 G_0(\xi_1,\xi_2,\xi_3,\xi) = \int \psi_\xi(r) \psi_{\xi_1}(r)
\psi_{\xi_2}(r) \psi_{\xi_3}(r) rdr
\]
\[
G_1(\xi_1,\xi_2,\xi_3,\xi) = 
\int_0^\infty \psi_\xi(r) \psi_{\xi_3}(r) \int_r^\infty \frac1s \psi_{\xi_1}(s)
\psi_{\xi_2}(s) ds rdr
\]
\[
G_2(\xi_1,\xi_2,\xi_3,\xi) = 
\int_0^\infty \psi_\xi(r) \psi_{\xi_1}(r) \int_r^\infty \frac1s \psi_{\xi_2}(s)
\psi_{\xi_3}(s) ds rdr
\]
Also we consider the truncated integrals
\[
 G_0^m(\xi_1,\xi_2,\xi_3,\xi) = \int \chi_{>m}(r) \psi_\xi(r) \psi_{\xi_1}(r)
\psi_{\xi_2}(r) \psi_{\xi_3}(r) rdr
\]
\[
G_1^m(\xi_1,\xi_2,\xi_3,\xi) = 
\int_0^\infty \psi_\xi(r) \psi_{\xi_3}(r) \int_r^\infty \chi_{>m}(s) \frac1s
\psi_{\xi_1}(s) \psi_{\xi_2}(s) ds rdr
\]
where $\chi_{> m}$ is a smooth approximation of the characteristic function of $[2^{m},\infty)$. 
We denote $D^{\alpha\beta\gamma\sigma}=(\xi_1 \partial_{\xi_1})^\alpha
 (\xi_3 \partial_{\xi_1})^\beta(\xi_3 \partial_{\xi_2})^\gamma(\xi_3 \partial_{\xi_3})^\sigma$
and
\begin{equation} \label{cfour}
g_{j k_1 k_2 k_3} = \frac{2^{\frac{j}2}}{|j|} \frac{\la k_3^-\ra 2^{-2k_3^+}}{2^{\frac{k_3}2}}
\end{equation}
 and  estimate these integrals as follows:

\begin{p1}
For $\xi$, $\xi_1$, $\xi_2$ and $\xi_3$ as above we have the bounds
\begin{equation}
|D^{\alpha\beta\gamma\sigma}
G_{0,1}(\xi_1,\xi_2,\xi_3,\xi)|
\lesssim_{\alpha\beta\gamma\sigma}  g_{j k_1 k_2 k_3}, 
\label{GG1}\end{equation}
\begin{equation}
|D^{\alpha\beta\gamma\sigma}
G_2(\xi_1,\xi_2,\xi_3,\xi)|
\lesssim _{\alpha\beta\gamma\sigma}  2^{k_3-k_1} g_{j k_1 k_2 k_3} .
\label{G2}\end{equation}
In addition, if $m+k_3 \geq 0$ then we have
\begin{equation}
|D^{\alpha\beta\gamma\sigma}
G_{0,1}^m(\xi_1,\xi_2,\xi_3,\xi)|
\lesssim _{\alpha\beta\gamma\sigma}  2^{k_3-k_1 -N(m+k_3)} g_{j k_1 k_2 k_3} 
\label{gg1R}\end{equation}

\label{p:g}\end{p1}

\begin{proof}
  Given the conditions (i),(ii) above and the asymptotic expansions
  for the functions $\psi_\xi$, it follows that the integral defining
  $G_0$ is oscillatory with frequencies 
$\pm \xi_1 \pm \xi_2 \pm \xi_3 \pm \xi$ of size $2^{k_3}$ and larger.  Hence the
contributions of regions $A_m$ decay rapidly for $m > -k_3$.
 Only the region $A_{<-k_3}$ has a nontrivial contribution, which we 
estimate directly.  If $k_3 < 0$ then we obtain
\[
|G_0| \lesssim \ 2^{\frac{k_1}2} 2^{\frac{k_2}2}
\frac{ 2^{\frac{k_3}2}}{|k_3|} \frac{2^{\frac{j}2}}{|j|}
\int_{r < 2^{-k_3}} \frac{|\log(1+ r^2)|^2}{1+2^{k_1} r } rdr
\approx    \frac{2^{\frac{j}2}}{|j|}\frac{|k_3|}{2^{\frac{k_3}2}}
\]
If  $k_3 \geq 0$ then there is some further gain, as $\psi_\xi$
no longer reaches the logarithmic part before the oscillatory regime.
 In that case we obtain
\[
|G_0| \lesssim 2^{\frac{k_1}2} 2^{\frac{k_2}2}2^{\frac{k_3}2} \frac{2^{\frac{j}2}}{|j|}
\int_{r  \leq 2^{-k_3}} \frac{(r 2^{k_3})^2  r^2}{1+2^{k_1} r} rdr \approx 
\frac{2^{\frac{j}2}}{|j|} 2^{-\frac{5k_2}2} 
\]
Adding the differentiation operator $D^{\alpha\beta\g\sigma}$ does not alter
the pointwise bounds used above. The estimate for the cut-off $G_0^m$
also follows from the above considerations.

In the case of $G_2$ the inner integral is oscillatory with frequencies of size
$2^{k_1}$. Hence we obtain
\[
\int_r^\infty \frac1s \psi_{\xi_2}(s) \psi_{\xi_3}(s) ds
\approx \frac1{2^{k_1}r} \psi_{\xi_2}(r) \psi_{\xi_3}(r), \quad r \gg 2^{-k_1}
\]
\[
\left |\int_r^\infty \frac1s \psi_{\xi_2}(s) \psi_{\xi_3}(s) ds\right|
\lesssim 2^{\frac{k_1}2} 2^{\frac{k_3}2} \frac{\la k_1^-\ra \la (k_1 - \ln(1+r))^{-} \ra} {2^{2 k_1^+}}
\frac{2^{2k_3^+}}{\la k_3^-\ra},  
\quad r \lesssim 2^{-k_1}
\]
Then we conclude as in the case of $G_0$.

Finally we consider $G_1$. Then we no longer want to estimate the inner integral.
Instead we integrate by parts,
\[
G_1= 
\int_0^\infty   \frac1r \int_0^r  \psi_\xi(s) \psi_{\xi_3}(s) sds \psi_{\xi_1}(r)
\psi_{\xi_2}(r) dr
\]
Now the inner integral is again oscillatory, and using the orthogonality
of $\psi_\xi$ and $\psi_{\xi_3}$ we can switch the inner integration
to $[r,\infty)$ and estimate
\[
\int_0^r  \psi_\xi(s) \psi_{\xi_3}(s) sds \approx  \frac{r}{2^{k_3}} \psi_\xi(r) \psi_{\xi_3}(r), 
\quad r \gg 2^{-k_3}
\]
\[
|\int_0^r  \psi_\xi(s) \psi_{\xi_3}(s) sds| \lesssim \frac{2^{\frac{j}2}}{|j|}\frac{2^{\frac{k_3}2} 2^{2k_3^+}}
{\la k_3^- \ra} 
r^2 \ln^2(1+r^2) , \quad r \lesssim 2^{-k_3}
\]
Then a similar argument to the one used for $G_0$ leads to the same bound, 
as the main contribution arising from the region $r \xi_3 \approx 1$ rests unchanged. 
A similar argument gives the estimate for $G_1^m$. 

\end{proof}

\section{The linear $\tilde H$ Schr\"odinger 
equation}
\label{linear}

Here we consider bounds and function spaces associated to the linear 
$\tilde H$  evolution
\begin{equation}
(i \partial_t - \tilde H) \psi = f, \qquad \psi(0) = \psi_0
\label{wlin-eq}\end{equation}
restricted to radial functions. We recall that the operator $\tilde H$
has the form
\[
 \tilde H = -\Delta+\V, \qquad \V = \frac{4}{r^2(1+r^2)}
\]
and, restricted to radial functions, admits the factorization
$\tilde H = L L^*$. In the first part of the section we 
introduce several relevant function spaces associated to this evolution,
and in the second part we prove that \eqref{wlin-eq} is well-posed
in these spaces.

\subsection{ Function spaces}
\subsubsection{ Globally defined spaces}
To measure solutions we will use the energy norm $L^\infty_t L^2_x$,
the Strichartz norm $L^4_{t,x}$ whose one-dimensional correspondents are
$L^\infty_t L^2_r$, respectively $L^4_{t,r}$. We also use the local energy norm
defined by 
\[
 \| \psi\|_{LE} = \| [r \log(2+r)]^{-1}\psi\|_{L^2}
\]
Combining these norms we define the space $S$ for solutions to
\eqref{wlin-eq} and the dual type space $N$ (precisely, $S=N^*$)
 for the inhomogeneous
term in \eqref{wlin-eq}.
\[
S = L^\infty_t L^2_r \cap L^4_{tr} \cap LE, \qquad
N = L^1_t L^2_r +  L^\frac43_{tr} +  LE^*
\]

\subsubsection{ Frequency localized spaces}
For many of our estimates we need to be more precise and work with a
dyadic Littlewood-Paley decomposition in the $\tilde H$-frequency, 
\[
 \psi =  P_k \psi
\]
To measure frequency $2^k$ waves we define a local energy
space $LE_k$,
\[
 \| \psi\|_{LE_k} =  2^{k} \| \psi\|_{L^2(A_{<-k})} + \sup_{m > -k}
2^{\frac{k-m}2} \|\psi\|_{L^2(A_{m})}
\]
as well as an adapted $L^4_k$ norm (allowed due to the radial symmetry):
\[ 
\| \psi\|_{L^4_k} = \sup_{m} \max\{2^{-\frac{m+k}2},2^{\frac{m+k}8}\}
\|\psi\|_{L^4(A_m)}.
\]
The dual norms are denoted by $LE_k^*$, respectively $L^\frac43_k$.
The frequency adapted versions of the $S$ and $N$
norms are 
\[
S_k = L^\infty_t L^2_r \cap L^4_k \cap LE_k, \qquad
N_k = L^1_t L^2_r +  L^\frac43_{k} +  LE_k^*, \qquad S_k = N_k^*
\]  
Square summing these norms we obtain the spaces
$l^2 S$ and $l^2 N$ with norms
\[
\| \psi\|_{l^2 S}^2 = \sum_{k\in \Z} \|P_k \psi\|_{ S_k}^2,
\qquad \|f\|_{l^2 N}^2 = \sum_{k\in \Z} \|P_k f\|_{N_k}^2,
\]
Given the nice bound \eqref{ker1} on the kernel of the projectors 
$P_k$ it is easy to see that these are dual spaces, thus justifying 
our notation.

We remark that in the frequency localized setting one has the usual 
Bernstein type estimates, with an additional improvement near $r=0$.
Precisely, from  the pointwise bounds \eqref{ker1} for the 
spectral projector kernels  we obtain
\begin{l1}
The following frequency localized pointwise bounds hold:
\begin{equation}
\|m_k^{-1} P_k \psi\|_{L^2_t L^\infty_r} \lesssim \|P_k \psi\|_{LE_k}
\label{point-le}\end{equation}
\begin{equation}
  \|(1+2^k r)^\frac12 m_k^{-1} P_k \psi\|_{L^\infty_{t,r}}
 \lesssim 2^k \|P_k \psi\|_{L^\infty L^2}
\label{point-e}\end{equation}
\end{l1}

\subsubsection{ The $Z$ spaces} \label{zdef}
The $Z$ spaces, which are used later in the paper for the parameter
$\lambda$ which tracks the evolution of the Schr\"odinger map along
the soliton family, are defined by
\[
Z = \dot W^{1,1} + \dot H^{\frac12,1} , \qquad Z_0=Z + L^2 \cap L^\infty
%, \qquad Z_1=Z \cap \dot H^1
\]
Concerning these spaces we need the following
\begin{l1}
The spaces  $Z$ and $Z_0$ are algebras.
\end{l1}
The proof is not very difficult and left to the reader.

\subsubsection{The $X^{s,b}$ spaces} For a function $\psi: \R \times
(0,\infty) \rightarrow \C$, we define its space-time Fourier transform
by $\hat{\psi}=\mathcal{F}_t \mathcal{F}_{\tilde H} \psi$, where
$\mathcal{F}_t$ is the time Fourier transform.  We define the
modulation localization operators $\{Q_j\}_{j \in \Z}$ by
$\widehat{Q_j\psi}(\tau,\xi)=\chi_{j}(|\tau + \xi^2|)
\hat{\psi}(\tau,\xi)$ where $\chi_{j}$ is a smooth characteristic
function of the set $\{ \tau \sim 2^j\}$.  We define $Q_{< j},
Q_{\leq j}, Q_{>j}, Q_{\geq j}$ in a similar way.

The $X^{s,b}$ type spaces 
$\dot X^{0,\pm \frac12,1}$  
and $\dot X^{0,\pm \frac12,\infty}$ 
associated to the $\tilde H$ flow are defined as
\[
\|\psi\|_{\dot X^{0,\pm \frac12,1}} = \sum_{j} 2^{\pm \frac{j}2}
\| Q_{j} \psi \|_{L^2}, \qquad 
\|\psi\|_{\dot X^{0,\pm \frac12,\infty}} = \sup_{j} 2^{\pm \frac{j}2}
\| Q_{j} \psi \|_{L^2}
\]  
These spaces play a less prominent role in this paper, as they are used only
at high modulations $j > 2k$ where $2^k$ is the frequency of $\psi$.
Precisely, we use them to define the dyadic space $S^r_k$
with norms
\[
\| \psi \|_{S_k^r} = \|\psi\|_{S_k} + \|Q_{>2k} \psi\|_{\dot X^{0,\frac12,1}+W^{1,1} L^2}
\]
We observe that due to the truncation to high modulations in the
second term above, we can replace the norm by an equivalent one
and write
\[
\| \psi \|_{S_k^r} = \|\psi\|_{S_k} + \|Q_{>2k} \psi\|_{Z L^2}
\]
Somewhat similarly, for the inhomogeneous term we define the space  $N^r_k$
by
\[
N^r_k = L^1 L^2 + Q_{>2k} \dot X^{0, -\frac12,1} 
\]
Summing up dyadic contributions in $l^2$ we obtain the spaces
$l^2 S^r$ and $l^2 N^r$:
\[
\|\psi \|_{l^2 S^r}^2 = \sum_{k} \|P_k \psi\|_{S_k^r}^2,
\qquad \|f\|_{l^2 N^r}=  \sum_{k} \|P_k f\|_{N_k^r}^2
\]
These norms, used only  on frequency $2^k$ functions, represent
a modest strengthening of the $S_k$ norms but only for 
high modulations. Their role in this paper is twofold.
On one hand, they represent all the information we are able to transfer
from the time dependent frame setting in the next section
back into the fixed frame setting; on the other hand, they are 
critically used in Section~\ref{seclambda} to recover the regularity
of the parameter $\lambda$ describing the evolution
of the Schr\"odinger map  along the soliton family.

\subsubsection{The $U^2$ and $V^2$ spaces}

Given a Hilbert space $\mathcal H$ (which in our case will be
$L^2(rdr)$), and $1 \leq p < \infty$, we define the spaces $U^p
\mathcal H$ and $V^p \mathcal H$ as follows:
\begin{itemize}
\item[a)] $U^p \mathcal H \subset L^\infty \mathcal H$ is an atomic
  space, where the atoms are step functions
\[
a = \sum_{k} 1_{[t_k,t_{k+1})} u_k, \qquad  \sum_k \|u_k\|_{\mathcal H}^p \leq 1
\]
with $t_k$ arbitrary finite increasing sequence in $[-\infty,\infty)$.

\item[b)] $V^p H \subset L^\infty H$ is the space of all right
  continuous $H$ valued functions for which the following norm is finite:
\[
\|u\|_{V^p H}^p = \sup_{t_k \nearrow} \sum_k \| u(t_{k+1}-u(t_k))\|_{H}^p 
\]
where the sup norm is over all increasing sequences $\{t_k\}$
as above. 
\end{itemize}
In our case we use the above definitions to construct the $U^p_{\tilde H} L^2$
and $V^p_{\tilde H} L^2$  spaces associated to the   $\tilde H$ flow
by
\begin{equation}
\| \psi\|_{U^p_{\tilde H} L^2}= \|  e^{-it \tilde H} \psi(t)\|_{U^p L^2},
\qquad \| \psi\|_{V^p_{\tilde H} L^2}= \|  e^{-it \tilde H} \psi(t)\|_{V^p L^2}
\end{equation}

Such spaces were introduced in the study of dispersive equations in
unpublished work of the second author. For more details we refer the reader
to \cite{KT}, \cite{BT-DNLS} and \cite{HHK}. In the context of Schr\"odinger maps
such spaces were also used in \cite{BIKT}.

We are primarily interested in the case $p=2$. There we have the embeddings
\begin{equation} \label{XUV}
\dot X^{0,\frac12,1} \subset U^2_{\tilde H} L^2 \subset V^2_{\tilde H} L^2
\subset  \dot X^{0,\frac12,\infty}
\end{equation}
Another favorable property of these spaces is that they are stable with 
respect to modulation truncations:
\begin{equation} \label{modtr}
Q_{< j} P_k: U^2_{\tilde H}L^2 \to U^2_{\tilde H}L^2,  \qquad 
Q_{< j} P_k: V^2_{\tilde H}L^2 \to V^2_{\tilde H}L^2
\end{equation}

For the inhomogeneous term we also define the space $DU^2_{\tilde H}
L^2$ as
\[
DU^2_{\tilde H} L^2 = \{ (i \partial_t - \tilde H)\psi; \psi \in  U^2_{\tilde H} L^2\} 
\]
with the induced norm. Here the derivatives are interpreted as
distributional derivatives. This satisfies 
\begin{equation}
\dot X^{0,-\frac12,1} \subset DU^2_{\tilde H} L^2  \subset \dot X^{0,-\frac12,\infty} 
\end{equation}

\subsubsection{ The sharp spaces}
Here we define our strongest dyadic spaces,  namely 
 $\Ss_k$ for frequency $2^k$ solutions, with norm
\[
 \| \psi \|_{\Ss_k}^2 =  \| \psi\|^2_{(S_k \cap V^2_{\tilde{H}} L^2)}
+ \| Q_{>2k} \psi\|^2_{W^{1,1} L^2 + \dot X^{0,\frac12,1}}
\]
as well as the space $\Ns_k$ for the inhomogeneous term,
with norm
\[
\| f \|_{\Ns_k}^2 =   \| f\|^2_{(N_k + DU^2_{\tilde{H}} L^2)}
+ \| Q_{>2k} f \|^2_{L^1  L^2 + \dot X^{0,-\frac12,1}}
\]
As before, we also define the full norms $l^2 \Ss$ and $l^2 \Ns$
by
\[
\| \psi\|_{l^2\Ss}^2 = \sum_k  \| P_k \psi\|_{l^2\Ss_k}^2,
\qquad 
\| \psi\|_{l^2\Ns}^2 = \sum_k  \| P_k \psi\|_{l^2\Ns_k}^2
\]
The $V^2$ and $DU^2$ spaces have been added in in order to allow for a
harmless transition between the high and low modulations, and also to
simplify some proofs in this section. Otherwise, the $DU^2$ norm above
plays no role. The $V^2$ space does play a role though, namely to
allow for better bounds when truncating in modulation.

The $\Ss$ and $\Ns$ type spaces are needed at two crucial points 
in the article. First, we use them to establish the well-posedness 
of the non-autonomous $\tilde H_\lambda$ Schr\"odinger flow in the 
next section. Secondly, we use them for the bootstrapping estimates 
in the nonlinear problem in Section~\ref{nonlin}.

\subsubsection{ Restrictions to compact intervals}
\label{extensions}
For the purpose of bootstrap arguments, many  of our estimates need to be 
proved first on compact time intervals $I= [0,T]$.  Thus we need to define 
our function spaces also on such intervals. This is done in a standard 
manner, in terms of extensions to the full real line, by setting 
\[
\| f\|_{\mathcal X(I)} = \inf  \{ \| f^{ext} \|_{\mathcal X}; \ \text{$f^{ext}$ extends $f$
from $I$ to $\R$} \}
\]
where $\mathcal X$ is any of the spaces previously introduced in this section.
In fact, it is only the use of the $\dot X^{0,\pm \frac12,1}$ structure at
high modulations which requires the use of extensions.

We say that an extension $f^{ext}$ of $f$ is suitable if $  \| f^{ext} \|_{\mathcal X} \sim
\| f\|_{\mathcal X(I)}$. Some ways of producing suitable extensions are described
next:
\begin{itemize}
\item For $f$ in $Z$ or $Z_1$ a suitable extension is given by 
\[
f^{ext}(t) = \left\{
\begin{array}{cc} f(a), & \quad  t \leq a \\
f(t), & \quad a \leq t \leq b \\
f(b), & \quad  b \leq t \\
\end{array}
\right.
\]
\item For $f$ in $Z_0$ a suitable extension is the zero extension.
\item For $f$ in $S$, $N$, $S_k$, $N_k$, $l^2 S$ and $l^2 N$
a suitable extension is the zero extension.

\item For $\psi$ in the spaces $S^r_k$, $l^2S^r$, $\Ss_k$, $l^2 \Ss$ 
a suitable extension is obtained by solving the homogeneous equation 
outside $I$,
\begin{equation} \label{extpsi}
\psi^{ext}(t) = \left\{
\begin{array}{cc}
e^{-i(t-a)\tilde H} \psi(0), & \quad t \leq 0 \\
\psi(t), & \quad a \leq t \leq b \\
e^{i(t-b)\tilde H} \psi(T), & \quad t \geq T
\end{array}
\right.
\end{equation}
Here a nonzero extension is required due to the high modulation structure
of the sharp spaces. This high modulation structure is of the form $Z L^2$,
therefore this extension can be thought of as  a direct counterpart of the $Z$ 
extension.

\item For the spaces $N^r_k$, $l^2N^r$, $\Ns_k$, $l^2 \Ns$ 
a suitable extension is the zero extension. This is less straightforward,
and it involves proving estimates of the type
\begin{equation} \label{xs-bcut}
\| 1_{[0,T]} \psi \|_{\ldNs} \lesssim \| \psi \|_{\ldNs}
\end{equation}
\end{itemize}
We outline the proof of \eqref{xs-bcut}. It suffices to consider its
dyadic counterpart. Of all components of the $\Ns_k$ norm, only the 
high modulation part is not trivially stable with respect to time truncations.
But at high modulation $\Ns_k$ has a $(\dot H^{-\frac12,1}+L^1)L^2$ structure,
 therefore  \eqref{xs-bcut} follows from an one-dimensional estimate
\[
\| \chi_{[0,T]} f \|_{\dot H^{-\frac12,1}+L^1} \lesssim \| f \|_{\dot H^{-\frac12,1}+L^1}
\]
In fact we only need to show $ \| \chi_{[0,T]} f \|_{H^{-\frac12,1}+L^1} \lesssim \| f \|_{H^{-\frac12,1}}$,
which can be further reduced to $ \| \chi_{[0,T]} f_k \|_{H^{-\frac12,1}+L^1} \lesssim 2^{-\frac{k}2} \| f_k \|_{L^2}$,
where $f_k = P_k f$ and $P_k$ are the standard Little-Paley projectors. It is obvious that
\[
\| P_{\gtrsim  k} ( \chi_{[0,T]} f_k) \|_{H^{-\frac12}} \lesssim 2^{-\frac{k}2} \| f_k \|_{L^2}
\]
so it is enough to show that
\[
\| P_{\ll k} ( \chi_{[0,T]} f_k) \|_{L^1} \lesssim 2^{-\frac{k}2} \| f_k \|_{L^2}
\]
But this follows from the straightforward estimate $\| P_k \chi_{[0,T]} \|_{L^2} \lesssim 2^{-\frac{k}2}$. 

\subsubsection{ Relations between spaces}
We summarize the relations between the spaces we have defined so far,
as well as some other simple properties for them, in the following
\begin{p1}
a) The following dyadic embeddings hold
\begin{equation}\label{dyadicemb}
 \Ss_k \subset S_k^r \subset S_k, \qquad N_k,N_k^r \subset \Ns_k 
\end{equation}
b) The following embeddings hold:
\begin{equation}\label{emb}
l^2 \Ss \subset l^2 S^r \subset l^2 S \subset S,
\qquad N \subset l^2 N \subset l^2 \Ns, \quad 
\end{equation}
c) Modulation localizations:
\begin{equation}\label{modcut}
\|Q_{ < j} P_k \psi\|_{\Ss_k} \lesssim \la (2k-j)^+ \ra \| P_k \psi\|_{\Ss_k}
\end{equation}
\end{p1}
\begin{proof}
  a) The first (sequence of) embeddings follows directly from the
  definitions.  The only nontrivial part of the second embedding is
  due to the contribution of the $L^\frac43$ component of $N$ at
  high modulation. There we write on dyadic pieces
\[
L^\frac43 \subset L^1 L^2 + L^2 L^1 \subset  L^1 L^2 + 
2^{k} L^2 
\]
 In the first step we can preserve the frequency localization
since by Proposition~\ref{p-lp} the spectral projectors are bounded in all
$L^p$ spaces. In the second step we use Bernstein's inequality,
which is valid in our setting due to the kernel bounds \eqref{ker1} for the 
spectral  projectors.

b) Given part (a), it remains to show that $l^2S \subset S$ and $N
\subset l^2 N$.  By duality it suffices to establish the first
embedding.  By the fixed time almost orthogonality of the $P_k$'s we
have
\[
 \| \psi\|_{L^\infty L^2}^2 \lesssim \sum_{k \in \Z} \|P_k \psi\|_{L^\infty L^2}^2
\]
The $L^4$ norms are similarly easy to add,
\[
 \| \psi\|_{L^4(A_m)} \lesssim \sum_{k \in \Z} 2^{-\frac{|m+k|}8} \|P_k \psi\|_{L^4_k}
\]
which leads to 
\[
\sum_{m \in \Z} \| \psi\|_{L^2(A_m)}^2 \lesssim \sum_{k \in \Z} \| P_k \psi\|_{L^4_k}^2 \lesssim \| \psi \|_{l^2 S}^2
\]
from which $\| \psi \|_{L^4} \lesssim \| \psi \|_{l^2 S}$ follows.
Finally we consider the local energy norms, for which we need to show that
\begin{equation}
 \| \psi\|_{LE}^2 \lesssim \sum_{k \in \Z} \|P_k \psi \|_{LE_k}^2
\label{lel2le}\end{equation}

By a direct summation in the regions $r > 2^{-k}$ and by 
summing the better bounds in \eqref{point-le} in the regions
$r < 2^{-k}$ one obtains 
\[
\begin{split}
  \|\frac{1}{r \ln{(2+r)}} \sum_k \psi_k \|_{L^2(A_j)} & \lesssim \sum_{k
    \gtrsim  -j} \| \frac{\psi_k}{r} \|_{L^2(A_j)}
  +  \sum_{k \lesssim -j} \frac{1}{\la k^- \ra} 
 \| \psi_k \|_{L^2_tL_r^\infty(A_j)} \\
  & \lesssim \sum_{k \gtrsim  -j} 2^{-\frac{j+k}2} \| \psi_k \|_{LE_k} +
  \sum_{k \lesssim -j} \frac{1}{\la k^- \ra} \| \psi_k \|_{LE_k}
\end{split}
\]
By using the following two estimates on sequences
\[
\| (\sum_{k \geq j} \frac{a_k}{\langle k \rangle})_{j \geq 1} \|_{l^2}
\lesssim \| (a_k)_{k \geq 1} \|_{l^2}, \qquad \| (\sum_{k \geq j}
2^{\frac{j-k}2} a_k )_{j \in \Z} \|_{l^2} \lesssim \| (a_k)_{k \in \Z}
\|_{l^2}
\]
we obtain the desired estimate \eqref{lel2le} in the region $r \gtrsim 
1$. A slight variation of the above argument gives also \eqref{lel2le}
in the region $r \lesssim 1$.

c) From \eqref{XUV}, $ \|Q_{l} P_k \psi \|_{\Ss_k} \lesssim \|Q_{l} P_k \psi \|_{X^{0,\frac12,1}} 
\lesssim \| P_k \psi \|_{V^2 \tilde H}
\lesssim \| P_k \psi \|_{\Ss_k}$ for any $j \leq l \leq 2k$. Then \eqref{modcut} follows, as 
there is nothing to prove if $j \geq 2k$.

\end{proof}

\subsection{Estimates for the linear $\tilde H$ 
Schr\"odinger flow}

Our main well-posedness result concerning the linear $\tilde H$ equation
is as follows:
\begin{p1}\label{p:mainS} 
The solution $\psi$ to \eqref{wlin-eq} satisfies the bound:
\begin{equation}
 \| \psi\|_{\ldSs} \lesssim \|\psi_0\|_{L^2} + \|f\|_{\ldNs}
\label{mainS} \end{equation}
\end{p1}

\begin{proof}
The bound \eqref{mainS} follows by dyadic summation from its
frequency localized version:
 \begin{equation}
 \| \psi\|_{\Ss_k} \lesssim \|\psi_0\|_{L^2} + \|f\|_{\Ns_k}
\label{mainSk} \end{equation}
whenever $\psi$ is localized at $\tilde H$-frequency $2^k$.
The proof of \eqref{mainSk} proceeds in several steps:

{\bf STEP 1: Frequency localized local energy decay.}
Here we consider functions $f$, $\psi_0$ which are 
localized at frequency $2^k$, and prove that the 
solution of \eqref{wlin-eq} obeys the following bound
\begin{equation}
  \| \psi\|_{LE_k} + 2^{-k} \| \partial_r \psi\|_{LE_k} \lesssim 
\| \psi_0\|_{L^2} + \|f\|_{LE_k^*}.
\label{lek}\end{equation}
Our approach is in the spirit of the one used by the second author in
\cite{T}, using the positive commutator method.

First we say that a sequence $\{\alpha_n\}_{n \in \Z}$ is slowly varying if
\[
|\ln{\alpha_j} - \ln{\alpha_{j-1}}| \leq 2^{-10}, \qquad \forall j \in \Z.
\]
Based on such a sequence we introduce the normed space $X_{k,\alpha}$ and
its dual $X_{k,\alpha}'$ as follows
\[
\begin{split}
\| u \|_{X_{k,\alpha}}^2 & = 2^{2k}\| u \|^2_{L^2(A_{< -k})} + 2^k \sum_{j
  \geq -k } \alpha_j 2^{-j} \| u \|^2_{L^2(A_j)} \\
\| u \|_{X_{k,\alpha}'}^2 & = 2^{-2k}\| u \|^2_{L^2(A_{< -k})} + 2^{-k}
\sum_{j \geq -k } \alpha_j^{-1} 2^{j} \| u \|^2_{L^2(A_j)}
\end{split}
\]

For all slowly varying sequences $\{\alpha_n\}_{n \in \Z}$ with $\sum_n
\alpha_n =1$, we claim that
\begin{equation} \label{ldp}
\| \psi \|_{X_{k,\alpha}}+ 2^{-k} \| \partial_r \psi\|_{X_{k,\alpha}}
 \lesssim \| \psi_0 \|_{L^\infty_t L^2_r} + \| f \|_{X'_{k,\alpha}} 
\end{equation}
Assuming that \eqref{ldp} is true, then we can consider another slowly
varying sequence $\{\beta_n\}_{n\in Z}$ with $\sum_n \beta_n=1$ and
apply the result in \eqref{ldp} for $\{\alpha_n+\beta_n\}_{n \in Z}$ to
obtain
\[
\| \psi \|_{X_{k,\alpha+\beta}}+ 2^{-k} \| \partial_r \psi\|_{X_{k,\alpha+\beta}}
\lesssim \| \psi_0 \|_{L^\infty_t L^2_r} + \| f \|_{X'_{k,\alpha+\beta}}
\] 
from which we derive the weaker estimate
\begin{equation} \label{ldp2}
\| \psi \|_{X_{k,\alpha}}+ 2^{-k} \| \partial_r \psi\|_{X_{k,\alpha }} \lesssim 
\| \psi_0 \|_{L^\infty_t L^2_r} + \| f \|_{X'_{k,\beta}}
\end{equation}

Since any $l^1$ sequence can be dominated by a slowly varying sequence
with a comparable $l^1$ size, we can drop the assumption in
\eqref{ldp2} that $\alpha$ and $\beta$ are slowly varying. By maximizing the
right-hand side with respect to $\alpha \in \l^1$ and by minimizing the
left-hand side with respect to $\beta \in l^1$, we obtain \eqref{lek}.

The remaining part of this step is devoted to the proof of \eqref{ldp}. 
We start by introducing the antisymmetric multiplier
\[
Q u = \chi(2^k r) r\partial_r u + r\partial_r (\chi(2^k r) u)
\]
where $\chi$ will be chosen to be a smooth function related to
the slowly varying sequence $\alpha_n$. Note that if the problem had a
scale invariance then one could rescale it to $k=1$ and discard the
factor of $2^k$ in the construction of $k$. But this is not the case
for \eqref{wlin-eq}.

Using the equation for $\psi$ we obtain
\[
\begin{split}
  \Re \int_0^T \langle Q \psi, f \rangle ds & = \Re \int_0^T \langle Q \psi, 
(i \partial_t - \tilde{H}) \psi \rangle ds \\
  & = - \Im \int_0^T  \langle Q \psi, \partial_t \psi \rangle - 
\Re \int \langle Q \psi, \tilde{H} \psi \rangle \\
  & = - \frac12 \Im \int_0^T \partial_t \langle Q \psi, \psi \rangle -
 \Re \int \langle Q \psi, \tilde{H} \psi \rangle \\
\end{split}
\]
which, by rearranging terms, becomes
\begin{equation} \label{idt}
- \Re \int_0^T \langle Q \psi, f \rangle ds  - \frac12 \Im \langle Q \psi, 
\psi \rangle |_0^T = \Re \int \langle Q \psi, \tilde{H} \psi \rangle 
\end{equation}
The right hand side can be expanded as follows
\[
\begin{split}
\Re \int \langle Q \psi, \tilde{H} \psi \rangle & = \Re \int \langle Q \psi, - \Delta \psi + \tilde{V} \psi  \rangle \\
& = \Re \int \langle Q \partial_r \psi, \partial_r \psi \rangle + \Re \int \langle [\partial_r,Q] \psi, \partial_r \psi \rangle 
+ \Re \int \langle Q \psi, \tilde{V} \psi \rangle \\
& = \Re \int \langle [\partial_r,Q] \psi, \partial_r \psi \rangle 
+ \frac12 \int \langle [\tilde{V},Q] \psi, \psi \rangle
\end{split}
\]
where we have used twice the antisymmetry of $Q$. We now compute the commutators and start with the easier one,
\[
\frac12[\tilde{V},Q] = - r \chi(2^k r) \partial_r \tilde{V} =  \chi(2^k r)\frac{4}{r^2+1}(\frac{1}{r^2}+\frac{1}{r^2+1}) > 0
\]
The other commutator is
\[
[\partial_r, Q] = 2(2^kr \chi'(2^k r)+ \chi(2^k r)) \partial_r  + (2^k \chi'(2^kr) + 2^{2k} r \chi''(2^kr)) 
\]
We now impose a first condition on the function $\chi$
\begin{equation} \label{chicond}
| (r\chi')'| \leq \delta(r \chi)'
\end{equation}
for some $0 < \delta \ll 1$. Using this and the Cauchy-Schwartz inequality we obtain
\[
\begin{split}
\Re \int \langle Q \psi, \tilde{H} \psi \rangle \gtrsim  \int a_k(r) 
(|\partial_r \psi|^2 - \delta 2^{2k}|\psi|^2) rdr dt + 
\frac12 \int [\tilde{V},Q] |\psi|^2 rdrdt
\end{split}
\]
where $a_k(r) =  \chi( 2^k r) +  2^k r \chi'(2^k r)$.
Hence, by  \eqref{idt} we have
\begin{equation} \label{idt1}
LHS(\eqref{idt}) \gtrsim  \int a_k(r) (|\partial_r \psi|^2 - \delta 2^{2k}|\psi|^2) rdr dt + \frac12 \int [\tilde{V},Q] |\psi|^2 rdrdt
\end{equation}
We claim that given a slowly varying sequence $\alpha_n$ and $\delta > 0$ 
we can find
$\chi$ satisfying \eqref{chicond}, so that
\begin{equation} \label{chicond2}
a_k(r) \gtrsim \frac{\alpha_{n+k}}{1+2^{n+k}}, \qquad r \approx 2^n
\end{equation}
and the following three fixed time bounds hold for functions 
localized at frequency $2^k$:
\begin{equation}
\| Q \psi\|_{L^2} \lesssim  \|\psi\|_{L^2}
 \label{lek1}\end{equation}
\begin{equation}
\| Q \psi\|_{X_{k,\alpha}} \lesssim  \|\psi\|_{X_{k,\alpha}}
\label{lek2}\end{equation}
\begin{equation}
\int_{\R} a_k(r) |\partial_r \psi|^2 rdr + \frac12 \int [\tilde{V},Q] |\psi|^2 rdrdt\gtrsim   2^{2k} 
\int_{\R} a_k(r) |\psi|^2 rdr
\label{lek3}\end{equation}

Using these three relations in the above integral estimate
we obtain
\[
\| \psi \|_{X_{k,\alpha}}^2 \lesssim 
\| \psi\|_{L^\infty L^2}^2 + \| \psi\|_{X_{x,\alpha}} \|f\|_{X_{k,\alpha}'}
\]
when all terms are restricted to the time interval $[0,T]$, but
with the a constant independent of $T$. This implies \eqref{ldp}. 

We now proceed with the construction of $\chi$ satisfying \eqref{chicond},
\eqref{chicond2} and \eqref{lek1}-\eqref{lek3}. We first increase $\alpha_n$
the so that it remains slowly varying and, in addition, satisfies

\begin{equation} \label{ren}
\alpha_n=1, \ \ \ \mbox{for} \  n \leq n_0-k; \qquad \sum_{n \geq n_0-k} \alpha_n \approx 1
\end{equation}

Here $n_0$ is a positive number to be chosen later. 
Based on this, we construct a slowly varying function $\alpha$ such that
\[
\alpha(s) \approx \alpha_n \ \ \mbox{if} \ \ s \approx 2^n
\]
and with symbol regularity
\[
|\partial^k \alpha(s)| \lesssim (1+s)^{-k} \alpha(s)
\]
Due to the first condition in \eqref{ren} we can take $\alpha$ such that
$\alpha(s)=1$ for $s \leq 2^{n_0-k}$. We then construct the function $\chi$
by
\[
s \chi(s) = \int_0^s \alpha(2^{-k} s) h(s) ds
\]
where $h_{n_0}$ is a smooth adapted variant of $ r^{-1}$,
namely $h_{n_0}(s)=1$ for $r \leq 2^{n_0}$ and $h(s) \approx 2^{n_0} s^{-1}$
for $s \geq 2^{n_0+1}$. One easily verifies the pointwise bounds

\begin{equation} \label{pchi}
\chi(s) \approx (1+2^{-n_0} s)^{-1},
\qquad |\chi^{(k)}(s)| \lesssim 2^{-kn_0} (1+2^{-n_0} s)^{-k-1}, \ \ k \leq 4
\end{equation}
Furthermore, we have 
\[
(s \chi(s))' = \alpha(2^k s) h_{n_0}(s) \gtrsim (1+2^{-n_0} s)^{-1.1}
\qquad 
|(s \chi'(s))'| \lesssim 2^{-n_0} (1+2^{-n_0} s)^{-2}
\]
It is a straightforward exercise to verify that $\chi$ satisfies
\eqref{chicond2}.  Furthermore, by taking $n_0$ large enough,
depending on $\delta$, we insure that $\chi$ satisfies also
the bound \eqref{chicond}.

Next we turn our attention to the estimates \eqref{lek1}-\eqref{lek3}.
For proving \eqref{lek1}-\eqref{lek2} we start by rearranging
\[
Q \psi = 2 \chi(2^k r) r\partial_r \psi + 2^k r \chi'(2^k r) \psi
\]
and using \eqref{pchi} we obtain $| 2^k r \chi'(2^k r)| \lesssim 1$
therefore
\[
\| 2^k r \chi'(2^k r) \psi \|_{L^2} \lesssim \| \psi \|_{L^2},
\qquad \| 2^k r \chi'(2^k r) \psi \|_{X_{k,\alpha}} \lesssim \| \psi \|_{X_{k,\alpha}}
\]
Using again \eqref{pchi}, we conclude the proof of \eqref{lek1}-\eqref{lek2}
by showing that
\begin{equation} \label{drpsik}
\| \chi(2^k r) r \partial_r \psi \|_{L^2} \lesssim \| \psi \|_{L^2},
\qquad
\| \chi(2^k r) r \partial_r \psi \|_{X_{k,\alpha}} \lesssim \| \psi \|_{X_{k,\alpha}}
\end{equation}
Since $\psi$ is localized at frequency $2^k$ we use an operator $P_k$
as in Proposition \ref{p-lp}, localizing at frequency $2^k$ and such
that $P_k \psi =\psi$. Then we use the characterization of $\partial_r
K_k(r,s)$ from \eqref{ker2} for the kernel $K_k$ of $P_k$; precisely,
by \eqref{ker2} and \eqref{pchi} we obtain
\[
| \chi(2^k r) r \partial_r K_k(r,s)| \lesssim \frac{2^{2k}}{(1+2^k|r-s|)^{N}(1+2^k(r+s))}
\]
Then \eqref{drpsik} follows since the above kernel has rapid
off-diagonal decay while the weights $\alpha_k$ are slowly varying.

 For \eqref{lek3} we claim the following estimate
\begin{equation} \label{pLstar}
\int a_k(r) | L^{*} \psi(r) |^2 rdr \gtrsim  2^{2k}\int a_k(r) |\psi(r)|^2 rdr
\end{equation}
Assuming \eqref{pLstar}, we can now complete the argument for
\eqref{lek3}. We write
\[
\partial_r \psi = - L^{*} \psi + \frac{h_3-1}{r} \psi = - L^{*} \psi - \frac{2\psi}{r(r^2+1)}
\]
which shows that 
\[
| L^{*} \psi|^2 \lesssim |\partial_r \psi|^2 + \frac{1}{r^2(r^2+1)^2}|\psi|^2
\]
Since $a_k(r) \lesssim (1+2^k r)^{-1} \approx \chi(2^k r)$, it follows that
\[
\int a_k(r) | L^{*} \psi(r) |^2 rdr \lesssim 
\int a_k(r)  |\partial_r \psi|^2 rdr + \frac12 \int [\tilde{V},Q] |\psi|^2
\]
Thus \eqref{pLstar} implies \eqref{lek}.

We finish this subsection with the proof of \eqref{pLstar}. 
For this we write $\psi$ in terms of $L^* \psi$ as
\[
\psi = P_k \psi = P_k \tilde H^{-1} L L^*\psi = (L^{-1} P_k)^* L^* \psi
\]
where the kernel $K^1_k(r,s)$ of $L^{-1} P_k$ was estimated in 
Proposition~\ref{p-lp}(b). We need to distinguish two cases:

i) $k \geq 0$. Then, by \eqref{ker3}, $K^1_k$ satisfies the symmetric bound
\[
|K^1_k(r,s)| \lesssim \frac{2^{k}}{(1+2^k|r-s|)^{N}(1+2^k(r+s))}
\]
and  \eqref{pLstar} directly follows.

ii) $k < 0$. The regular part $K^1_{k,reg}$ still satisfies the 
above bound, and causes no difficulties. For the resonant part $K^1_{k,res}$
we use \eqref{ker4} to derive the estimate
\[
| K^1_{k,res}(r,s)| \lesssim \frac{h_1(r)}{k(1+2^kr)^{N}} 
\frac{1}{(1+2^ks)^{N}}  
\]
which has a $2^{-k}$ bound from $L^2 \to L^2$ and decays rapidly 
above the $r, s = 2^{-k}$ threshold. Thus  \eqref{pLstar} again follows.

\vspace{.2in}

{\bf STEP 2: A dyadic  $L^4_t L^\infty_r$ bound for the homogeneous problem}

Here we establish the bound
\begin{equation} \label{l4li}
 2^{-\frac{k}2} \| \psi\|_{L^4_t L^\infty_r(A_{<k})}+ \sup_{j \geq - k}
2^{\frac{j}2} \| \psi\|_{L^4_t L^\infty_r(A_j)}
 \lesssim \|\psi_0\|_{L^2}
\end{equation}
for $\psi_0$ localized at frequency $2^k$.
The first term is easily bounded by interpolating between \eqref{point-e}
and \eqref{point-le}. Consider now $j \geq -k$, and $\chi_j$ a smooth 
bump function supported in $A_j$. The function $\psi_j = \chi_j \psi$
solves the equation
\[
 (i \partial_t - \tilde H) \psi_j = 2 \partial_r \chi_j \partial_r \psi+
\delta \chi_j \cdot \psi := f_j
\]
From the local energy decay estimate for $\psi$ we obtain the following bounds
\[
\|\psi_j(0)\|_{L^2} + 2^{\frac{k-j}2 }\| \psi_j\|_{L^2} + 2^{\frac{j-k}2 } \|f_j\|_{L^2} \lesssim \|\psi(0)\|_{L^2}
\]
We conjugate by $r^{\frac12}$ and set $v_j(t,r) = r^\frac12 \psi_j$.
A direct computation shows that $v_j$ solves a one dimensional Schr\"odinger
equation
\[
 (i \partial_t - \partial_r^2) v_j = r^\frac12 f_j + (r^{-2}+\tilde V(r)) v_j := g_j
\]
where 
\[
 \|v_j(0)\|_{L^2} + 2^{\frac{k-j}2 }\| v_j\|_{L^2} + 2^{\frac{j-k}2 } \|g_j\|_{L^2} \lesssim \|\psi(0)\|_{L^2}
\]
Applying the one dimensional $L^4 L^\infty$ Strichartz estimate
over each time interval of size $2^{j-k}$ we obtain
\[
 \| v_j\|_{L^4 L^\infty} \lesssim \|\psi(0)\|_{L^2}
\]
Returning to $\psi$ this yields
\[
 \| \psi \|_{L^4 L^\infty(A_j)} \lesssim 2^{-\frac{j}2} \|\psi(0)\|_{L^2}
\]
Hence \eqref{l4li} is proved.

{\bf STEP 3: The dyadic  $L^4$ bound for the homogeneous problem}

Here we establish the bound
\begin{equation} \label{l4hom}
  \| \psi\|_{L^4_k} \lesssim \|\psi(0)\|_{L^2}
\end{equation}
for $\psi$ localized at frequency $2^k$. The $L^4$ bound in 
$A_j$ with $j \leq -k$ follows from the first term in \eqref{l4li}
by H\"older's inequality. The $L^4$ bound in 
$A_j$ with $j > -k$ is obtained  by interpolating between the 
$L^4 L^\infty$ bound in the second term in \eqref{l4li},
the $L^2_{t,x}$ bound in $LE_k$ and the $L^\infty L^2$ 
energy estimate.

{\bf STEP 4: The role of the $U^p$ and $V^p$ spaces}
Here we show that
\begin{equation}\label{sdutosv}
\| \psi\|_{S_k \cap V^2_{\tilde H} L^2} \lesssim \|\psi(0)\|_{L^2} +
\|f\|_{N_k + DU^2_{\tilde H} L^2}
\end{equation}
By Steps 1 and 3 we know that for the homogeneous problem we
have
\begin{equation}\label{shom}
\|\psi\|_{S_k} \lesssim \|\psi_0\|_{L^2}, \qquad f = 0
\end{equation}
which implies \eqref{sdutosv} in this case.
By duality this shows that for the inhomogeneous problem we 
have
\begin{equation}\label{sdual}
\|\psi\|_{L^\infty L^2} \lesssim \|\psi_0\|_{L^2}+\|f\|_{N_k}
\end{equation}
Applying \eqref{shom} for each step of each $U^2$ atom, we further
obtain
\begin{equation}\label{semb}
\| \psi\|_{S_k} \lesssim \|\psi \|_{U^2_{\tilde H} L^2}
\end{equation}
which  suffices for $f \in  DU^2_{\tilde H} L^2$.
It remains to consider $f \in N_k$, which we further split
in two.

i) $f \in L^{\frac43}_k + L^1 L^2$.  
For any partition $\R = \cup I_l$
of the time interval into subintervals we have  
\[
\sum_l \|f \|_{ (L^{\frac43}_k + L^1 L^2)(I_l)}^\frac43 \lesssim 
 \|f \|_{ L^{\frac43}_k + L^1 L^2}^\frac43
\]
which combined to \eqref{sdual} yields
\[
\| \psi\|_{V^\frac43_{\tilde H} L^2} \lesssim \|f\|_{ L^{\frac43}_k + L^1 L^2}
\]
Since $V^{\frac43} \subset U^2$, the proof is concluded in this case.

ii) $f \in LE_k^*$. By Step 1 we have the $LE_k$ and $L^\infty L^2$
bounds for $\psi$.  On the other hand arguing as in case (i) above we
obtain
\[
\| \psi\|_{V^2_{\tilde H} L^2} \lesssim \|f\|_{ LE_k^*}
\]
This concludes the proof since $V^2 \subset U^4$,  and we have the
following variation of \eqref{semb},
\begin{equation}\label{semba}
\| \psi\|_{L^4_k} \lesssim \|\psi \|_{U^4_{\tilde H} L^2}
\end{equation}

{\bf STEP 5: The high modulation bound.}
Given \eqref{sdutosv}, to conclude the proof of the proposition it remains
to prove the high modulation bound
\[
\| Q_{>2k} \psi\|_{W^{1,1} L^2 + \dot X^{0,\frac12,1}} 
\lesssim \| Q_{>2k} f \|_{L^1 L^2 + \dot X^{0,-\frac12,1}} 
\]
This is straightforward; details are left for the reader.
\end{proof}

\section{ The time dependent linear evolution}
\label{lineart}

Here we consider the linear equation
\begin{equation}
(i \partial_t - \tilde H_\lambda) u = f, \qquad u(0) = u_0
\label{wlin-eq-l}\end{equation}
where
\[
 \tilde H = -\Delta+\V_\lambda , \qquad \V_\lambda = \frac{4}{r^2(1+\lambda^2
r^2)}
\]
In the space $L^2$ this can be viewed as a small perturbation of 
the $\lambda = 1$ problem in \eqref{wlin-eq}:

\begin{p1} \label{p:lin-l2}
Assume that $|\lambda-1|_{L^\infty} \ll 1$. Then the equation
\eqref{wlin-eq-l} is well-posed in $L^2$, and the 
following bound holds:
\begin{equation}
\| u\|_{\ldSs} \lesssim \|u(0)\|_{L^2} + \|f\|_{\ldNs}
\label{b:lin-l2}\end{equation}
\end{p1}
\begin{proof} Since $|V_\l - V_1| \lesssim |\l-1|(1+r^2)^{-2}$ it follows that
\[
\| (V_\l - V_1) h \|_{LE^*} \lesssim \| h \|_{LE}
\]
Therefore we can rewrite the \eqref{wlin-eq-l} as
\[
(i \partial_t - \tilde H) u = (V_\l - V_1)u + f, \qquad u(0) = u_0
\]
and treat the $(V_\l - V_1)u$ as a perturbation. The result follows then 
from \eqref{mainS}.

\end{proof}

Our main goal in this section is to study the equation
\eqref{wlin-eq-l} in the smaller space $LX$.  The condition $\|\lambda
-1\|_{L^\infty} \ll 1$ is no longer sufficient for the analysis in
$LX$. Instead we use the stronger topology $\Zu$ for $\lambda$ (see
Section~\ref{zdef}), and work with
\begin{equation}
\| \lambda -1\|_{\Zu} \leq \gamma \ll  1
\label{lambda-reg}\end{equation}

Our aim  will be achieved in two steps.
\begin{itemize} 

\item The spaces of type $\ldSs$, $\ldNs$ associated to the $\lambda = 1$ 
flow are not robust enough for  the variable $\lambda$ flow. Hence
we introduce some modified spaces $\WSs[\lambda]$, $\WNs[\lambda]$ 
adapted to the time dependent frame. To simplify some of the analysis,
we also provide some partial characterizations of functions in 
$\WSs[\lambda]$ and $\WNs[\lambda]$ with respect to the time independent frame
$\lambda = 1$.

\item We prove that if \eqref{lambda-reg} holds, then 
the evolution \eqref{wlin-eq-l} is globally well-posed for initial data in $LX$
and inhomogeneous term in $\WNs[\lambda]$, and the solution $\psi$ 
is uniformly bounded in $LX$, and further it belongs to $\WSs[\lambda]$.
\end{itemize}

\subsection{The spaces $\WSs[\lambda]$, $\WNs[\lambda]$}
Here we define the $[\l]$ type spaces as counterparts of the spaces
from the previous section which take into account time-variable Littlewood-Paley
projectors. We begin as usual with a dyadic decomposition, but with respect to the
$\lambda$ dependent frame,
\[
\psi = \sum_{k} P_k^\lambda \psi.
\]
For $\X \in \{S,N,\Ss,\Ns\}$ we define the space $l^2\X[\l]$ with norm
\[
\| \psi \|_{l^2 \X[\l]}= \left( \sum_{k} \| P_k^\l \psi \|_{l^2 \X}^2 \right)^\frac12 
\]
This gives rise to the spaces $l^2S[\l]$, $l^2N[\l]$, $l^2\Ss [\l]$, $l^2\Ns [\l]$
which correspond to $L^2$ initial data in \eqref{wlin-eq-l}.
For $LX$ initial data, on the other hand, we need to replace the $l^2$ 
dyadic summation with the same summation as in the $LX$ norm.
Hence we define
\[
\| \psi \|_{W\X[\lambda]} = \sum_{k < 0} \frac{1}{2^k |k|} \|
P_k^\lambda \psi\|_{l^2 \X} + \left( \sum_{k \geq 0} \| P_k^\lambda
  \psi\|_{l^2 \X}^2\right)^\frac12
\]
All of the above spaces $l^2\X[\lambda]$ and $W\X[\lambda]$ have their
finite time interval counterpart $l^2\X[\lambda](I)$ and
$W\X[\lambda](I)$, which are obtained by using $l^2 \X(I)$ instead of
$l^2\X$ in the above definitions. We especially remark that they are
{\bf not} obtained by restricting to $I$ a similar class of functions
over the entire real line; such a definition would be dependent on
specifying an extension of $\lambda$, which we wish to avoid. The
spaces $\WSs[\lambda] , \WNs[\lambda]$ play a fundamental role in our
analysis.

We remark that the functions $P_k^\lambda \psi$ are frequency localized in the time dependent frame
but not in the fixed $\lambda = 1$ frame. This will cause considerable
technical difficulties later on. Because of this, it is useful to
transfer as much information as possible back to the fixed frame.

\begin{p1}[{Characterizations of $\WSs[\lambda]$ and $\WNs[\lambda]$
    functions}] Suppose that $\lambda$ takes values in a compact
  subset of $(0,\infty)$. Then
  
  a) The following $S$ type norms are equivalent:
  \begin{equation}
    \| P_k^\lambda \psi\|_{l^2 S} \approx \|P_k^\lambda \psi\|_{S_k},
    \qquad \| P_k^\lambda g\|_{l^2 N} \approx \|P_k^\lambda g\|_{N_k}
    \label{le=lek}\end{equation}
  as well as
  \begin{equation}
    \|\psi\|_{l^2S[\lambda]} \approx \|\psi\|_{l^2 S}, \ \ \
 \|\psi\|_{l^2 N[\lambda]} \approx \|\psi\|_{l^2N}, \ \ \
 \|\psi\|_{\WS[\lambda]} \approx \|\psi\|_{\WS[1]}.
    \label{wsl=ws}\end{equation}
  Furthermore, we have the improved local energy decay
  \begin{equation}
    \left\| \frac{\psi}{\la r\ra^\frac12 \ln(1+r)} \right\|_{L^2} \lesssim
    \|\psi\|_{\WSs[\lambda]}
    \label{betterle}\end{equation}

  b) Assume in addition that $\lambda -1 \in \Zu$. Then the
  following inclusions hold:
  \begin{equation}\label{wsemb}
    \WSs[\lambda] \subset \WS^r[1]
  \end{equation}
\begin{equation}\label{wnemb}
    \WN^r[1] \subset \WNs[\lambda]
  \end{equation}
\end{p1}

We remark that all of the above bounds with the exception of \eqref{wsemb}
also hold trivially in any interval; this is because all the norms involved 
can be measured in an interval by taking the zero extension outside it.
The bound  \eqref{wsemb} also holds in any interval, but this is 
a more delicate matter which we will only be able to consider after we 
prove Proposition~\ref{p:lin-x} below.

\begin{proof} 
{\bf The estimate \eqref{le=lek}.}
This is trivial if $\lambda = 1$.  Otherwise, by definition,
\[
\| P_k^\lambda \psi\|_{l^2 S}^2 = \sum_{j \in \Z } \| P_{j}
P_k^\lambda \psi\|_{S_{j}}^2
\]
If $|j-k| \lesssim 1$ then the $S_k$ and $S_{j}$ norms are equivalent,
and the multipliers $P_{j}$ are bounded in $S_k$.  Thus we have
\begin{equation} \label{closek} \| P_{j} P_k^\lambda \psi\|_{S_{k}}
  \lesssim \| P_{j} P_k^\lambda \psi\|_{S_{j}} \lesssim \| P_k^\lambda
  \psi\|_{S_{k}}, \qquad |j-k| \lesssim 1
\end{equation}
Consider now the case $|j-k| \gg 1$. We write $ P_{j} P_k^\lambda=
P_{j} \tilde P_k^\lambda P_k^\lambda$.  For the kernel $K_{jk}(r,s)$ of $
P_{j} \tilde P_k^\lambda$ we use the estimate \eqref{kerbiprod}.
% \[
% |K_{jk}(r,s)| \lesssim \frac{1}{\la k^-\ra \la j^-\ra} \frac{ 2^{k+j-
%     |k-j|} 2^{-N(k^++j^+)}}{(1+2^{j} r)^{N} (1+ 2^k s)^{N}}
% \] 
Then a direct computation shows that
\begin{equation} \label{fark}
  \begin{split}
    \| P_{j} P_k^\lambda \psi\|_{S_{j}} \lesssim &\
    \frac{2^{(j-k)^-} 2^{-N(k^++j^+)}}{\langle k^- \rangle \langle j^- \rangle } \| P_k^\lambda \psi\|_{S_k},\\
    \| P_{j} P_k^\lambda \psi\|_{S_{k}} \lesssim &\
    \frac{2^{-\frac{|k-j|}2}2^{-N(k^++j^+)}} {\langle k^- \rangle
      \langle j^- \rangle} \| P_k^\lambda \psi\|_{S_k}
  \end{split}
\end{equation}
We use \eqref{closek} and \eqref{fark} to conclude the proof of the
first estimate in \eqref{le=lek}. On one hand we have
\[
\| P_k^\lambda \psi\|_{l^2 S}^2 \lesssim (1+ \sum_{j} \langle k^-
\rangle^{-2} \langle j^- \rangle^{-2}2^{-2N(k^++j^+)}) \| P_k^\lambda
\psi\|_{S_k}^2 \lesssim \| P_k^\lambda \psi\|_{S_k}^2
\]
Conversely, we denote the separation threshold by $k_0$ and compute
\[
\begin{split}
  \|P_k^\lambda \psi\|_{S_k} & \lesssim \sum_{|j-k| \leq k_0} \|P_{j}
  P_k^\lambda \psi\|_{S_{k}}
  + \sum_{|j-k| > k_0} \|P_{j} P_\mu^\lambda \psi\|_{S_k} \\
  & \lesssim c_{k_0} \| P_\mu^\lambda \psi\|_{l^2 S}
  + 2^{-\frac{k_0}2} \|P_k^\lambda \psi\|_{S_k} \\
\end{split}
\]
By appropriately adjusting $k_0$, the last term on the right can be
absorbed by the the first term $\|P_k^\lambda \psi\|_{S_k}$, thus
giving us the reverse inequality
\[
\|P_k^\lambda \psi\|_{S_k} \lesssim \| P_k^\lambda \psi\|_{l^2 S}
\]
This completes the proof of the first estimate in \eqref{le=lek}. The
second follows from a similar argument.

\medskip

{\bf The estimate \eqref{betterle}.} The proof of \eqref{betterle} is
almost identical to the proof of \eqref{emb}. The fact that $\lambda$
is not equal to $1$ makes no difference there.

\medskip

{\bf The estimates \eqref{wsl=ws}.} 
We only prove the third bound, which is more important in this 
article. The first two are similar but simpler.
The proofs of the two inclusions
are identical, so it suffices to show one of them, say $\WS[\lambda]
\subset WS[1]$. For fixed frequency $j$ we decompose
into a diagonal and an off-diagonal part
\begin{equation}\label{didi}
P_{j} \psi = \sum_{|j-k| \lesssim 1} P_{j} P_{k}^\l \psi +
 \sum_{|j-k| \gg 1} P_{j} P_{k}^\l \psi := (P_j \psi)_{\text{diag}} + 
(P_j \psi)_{\text{offd}}
\end{equation}
For the diagonal part it suffices to use the $S_j$ boundedness of
$P_{j}$.  For the off-diagonal part we use the first part of \eqref{fark} to
obtain
\[
\|(P_j \psi)_{\text{offd}}\|_{S_j} \lesssim \sum_{|k-j| \lesssim 1} \| P_{k}^\l
\psi\|_{S_k}+ \sum_{|k-j| \gg 1} \frac{2^{(j-k)^-}
  2^{-N(k^++j^+)}}{\langle k^- \rangle \langle j^- \rangle } \|
P_{k}^\l \psi\|_{S_k}
\]
To conclude it suffices to sum up the second term on the right with
respect to $j$ and the weights $2^{-j^-} \la j^-\ra^{-1}$. Indeed we
have
\[
\begin{split}
\sum_j \frac{2^{-j^-}}{\la j^- \ra} \|(P_j \psi)_{\text{offd}}\|_{S_j} \lesssim & \ 
\sum_j \frac{2^{-j^-}}{\la j^- \ra}\sum_k\frac{2^{(j-k)^-
  -N(k^++j^+)}}{\langle k^- \rangle \langle j^- \rangle } \|
P_{k}^\l \psi\|_{S_k}
\\
 \lesssim & \ \sum_k \frac{2^{-k^--Nk^+}}{\langle
  k^- \rangle^2} \| P_{k}^\l \psi\|_{S_k} \lesssim \| \psi\|_{\WS[\lambda]}
\end{split}
\]
where \eqref{le=lek} was used in the last step.

{\bf The estimate  \eqref{wsemb}.} Given \eqref{wsl=ws}, it
remains to bound the additional high modulation component in the
$S^r[1]$ norm.  We decompose $P_j \psi$ again as in \eqref{didi}.
The diagonal part  is estimated directly in $\Ss_j$.
The nontrivial part of the argument is to estimate the 
off-diagonal part. For these we decompose further 
\[
\begin{split}
  Q_{\gtrsim  2j} (P_{j}  f)_{\text{offd}} =  \sum_{|k-j| \gg 1} \sum_{h}
  Q_{\gtrsim  2j} P_{j} \tilde{P}_{k}^\l \tilde{P}_{h} P_{h}
  P_{k}^\l f
\end{split}
\]
The definition of $\WSs[\l]$ gives us good estimates on $P_{h}
P_{k}^\l f$ in $\Ss_{h}$, and, after applying the operator
$P_{j} \tilde{P}_{k}^\l \tilde{P}_{h}$, we need to estimate the
(high modulation) output  in $(X^{0,\frac12,1}+W^{1,1})L^2 = Z L^2$. 
For  the kernel $K_{jkh}$ of  $P_{j} \tilde{P}_{k}^\l \tilde{P}_{h}$
we use the representation in Proposition~\ref{p:prodthree} (a) or (b).
Hence  it suffices to consider 
kernels $K_{jkh}$ of the form
\[
K_{jkh}(r,s) = c_{jkh} g(\lambda) \phi_j(r) \phi_{h}(s) 
\]
where $|\phi_j(r)| \lesssim 2^{j}  (1+2^j r)^{-N}$, and similarly for $\phi_{h}$.
For such $K_{jkh}$ we  write
\[
 P_{j} \tilde{P}_{k}^\l \tilde{P}_{h} P_{h} f 
  P_{k}^\l f = c_{jkh} g(\lambda) \phi_j(r)  \la \phi_h, P_{h}
  P_{k}^\l f \ra 
\]
We first estimate the last inner product. Globally we use local energy to obtain
\[
\|  \la \phi_{h}, P_{h} P_{k}^\l f \ra\|_{L^2} \lesssim 
2^{-h} \|P_{h} P_{k}^\l f \|_{LE_{h}} 
\]
while at high modulation we have
\[
\|  Q_{>2h}\la \phi_{h}, P_{h} P_{k}^\l f \ra\|_{Z} \lesssim 
\| Q_{>2h} P_{h} P_{k}^\l f \|_{Z L^2}
\]
Combining the last two estimates 
we obtain
\[
\|  \la \phi_{h}, P_{h} P_{k}^\l f \ra\|_{Z} \lesssim 
 \|P_{h} P_{k}^\l f \|_{\Ss_{h}} 
\]
Due to the $Z$ algebra property this bound is not affected by
multiplication by $g(\lambda)$. Estimating $\phi_j$ in $L^2$ we 
obtain
\begin{equation}
\| Q_{>2j}  P_{j} \tilde{P}_{k}^\l \tilde{P}_{h} P_{h}  P_{k}^\l f \|_{Z L^2}
\lesssim c_{jkh}   \|P_{h} P_{k}^\l f \|_{\Ss_{h}} 
\end{equation}
where the modulation truncation $Q_{>2j}$ was harmlessly added 
at the end. Hence in order to estimate the high modulation component 
of $\|f\|_{WS^r[1]}$ we add up the dyadic pieces
\[
\begin{split}
\sum_{j}  \frac{2^{-j^-}}{\la j^-\ra} \| Q_{>2j}(P_j f)_{\text{offd}}\|_{ZL^2}
\lesssim & \  \sum_{j,k,h}  \frac{2^{-j^-}}{\la j^-\ra}  c_{jkh} 
 \|P_{h} P_{k}^\l f \|_{\Ss_{h}} 
\\
\lesssim & \  \sum_{k,h}  \frac{2^{-k^- - |k-h| - N(k^++h^+)}}{\la k^-\ra^2} 
\|P_{h} P_{k}^\l f \|_{\Ss_{h}} 
\\
\lesssim & \  \sum_{k}  \frac{2^{-k^ - - N k^+}}{\la k^-\ra^2} 
\|P_{k}^\l f \|_{\ldSs}  \lesssim \|f\|_{\WSs[\lambda]}
\end{split}
\]
which completes the proof.

\medskip

{\bf The estimate  \eqref{wnemb}.} 
We need to show that for $f$ at frequency $h < 0$ and modulation $\sigma > 2h$ 
we have
  \begin{equation} \label{inNl} \|  f \|_{\WNs[\l]} \lesssim \frac{2^{-h-\frac{\sigma}2}}{|h|  }\| f \|_{L^2}
  \end{equation}
We decompose $f$ as follows,
\[
f = f_0 + f_1 + f_2 + f_3, \qquad f_i = 
 \sum_{k,j \in \mathcal A_i} P_{j} P^\lambda_{k} \tilde{P}_h f  
\]
\[
\mathcal A_0 = \{ |j-k| \gg 1, \  |k-h| \gg 1\},
\qquad 
\mathcal A_1 = \{ |j-k| \lesssim 1, \  |k-h| \gg 1\}
\]
\[
\mathcal A_2 = \{ |j-k| \gg 1, \  |k-h| \lesssim 1\},
\qquad 
\mathcal A_3 = \{ |j-k| \lesssim 1, \  |k-h| \lesssim 1\}
\]
For indices in the set $\mathcal A_0$ we have the representation 
of $P_{j} P^\lambda_{k} \tilde{P}_h$  given in Proposition~\ref{p:prodthree}(a),
as a rapidly convergent series
of operators of the form
\[
T_{jkh} =  c_{jkh}    g(\lambda) \phi_{j}(r) \la \phi_h(s), \cdot \ra  
\]
with $c_{jkh}$ as in \eqref{kertriproda}. For the inner product we have
\[
\|\la \phi_h(s), f \ra  \|_{L^2} \lesssim  \|f\|_{L^2} 
\]
which immediately leads to
\[
\|P_{j} P^\lambda_{k} \tilde{P}_h f\|_{L^2}
+ 2^{j} \|P_{j} P^\lambda_{k} \tilde{P}_h f\|_{LE_{j}^*}
 \lesssim   
  \dfrac{ 2^{- |j-k| - |k-h|-N(j^++k^++h^+)}}{\la j^- \ra \la k^-
      \ra^2 \la h^- \ra} \|f\|_{L^2} 
\]
We will use the $LE^*_{j}$ bound at high frequencies, $2j > \sigma - 8$.
For smaller frequencies we use the $L^2$ bound at high modulation 
($\geq \sigma - 4$). For low modulations we instead obtain an $L^1 L^2$ 
bound. Here the idea is that $f$ is localized at high modulation, and the only way
to arrive to low modulations is to use high modulations of $g(\lambda)$.
Precisely we have
\[
Q_{<\sigma-4} P_{j} P^\lambda_{k} \tilde{P}_h  f = 
\sum c_{jkh}    [Q_{> \sigma -2} g(\lambda) ]\phi_{j}(r) \la \phi_h(s), f \ra  
\]
Since $g(\lambda) \in Z$, we have $\|Q_{> \sigma -2} g(\lambda)\|_{L^2}
\lesssim 2^{-\sigma/2}$. Hence we arrive at
\[
\|Q_{<\sigma-4} P_{j} P^\lambda_{k} \tilde{P}_h f\|_{L^1 L^2}
 \lesssim   2^{-\frac{\sigma}2}
  \dfrac{ 2^{- |j-k| - |k-h|-N(j^++k^++h^+)}}{\la j^- \ra \la k^-
      \ra^2 \la h^- \ra} \|f\|_{L^2} 
\]
Thus in all cases it follows that
\[
\| P_{j} P^\lambda_{k} \tilde{P}_h f\|_{\Ns_{j}}
 \lesssim   \min\{ 2^{-j},2^{-\frac{\sigma}2}\}
  \dfrac{ 2^{- |j-k| - |k-h|-N(j^++k^++h^+)}}{\la j^- \ra \la k^-
      \ra^2 \la h^- \ra} \|f\|_{L^2} 
\]
Summing up, we obtain for $f_0$
\[
\begin{split}
\| f_0\|_{\WNs[\lambda]} \lesssim & \ \sum_{j,k} \frac{2^{-j}}{\la j^-\ra}
  \min\{ 2^{-j},2^{-\frac{\sigma}2}\}
  \dfrac{ 2^{- |j-k| - |k-h|-N(j^++k^++h^+)}}{\la j^- \ra \la k^-
      \ra^2 \la h^- \ra} \|f\|_{L^2} \\
\lesssim & \ \sum_{j} \frac{2^{-j}}{\la j^-\ra}
  \min\{ 2^{-j},2^{-\frac{\sigma}2}\}
  \dfrac{ \la h-j\ra 2^{- |j-h| -N(j^+ +h^+)}}{\la j^- \ra^2  \la h^- \ra^2} \|f\|_{L^2} 
\\
\lesssim & \ 
  2^{-\frac{\sigma}2}
  \dfrac{  2^{- h -N h^+}}{ \la h^- \ra^3} \|f\|_{L^2} 
\end{split} 
\]
which is slightly stronger than needed.

For the terms in $f_1$ and $f_2$ the computation is almost
identical. Using \eqref{kertriprodb} instead of
\eqref{kertriproda}  we obtain
\[
\| P_{j} P^\lambda_{k} \tilde{P}_h f\|_{\Ns_{j}}
 \lesssim   \min\{ 2^{-j},2^{-\frac{\sigma}2}\}
  \dfrac{ 2^{- |j-h|-N(j^++h^+)}}{\la j^- \ra  \la h^- \ra} \|f\|_{L^2} 
\]
The summation with respect to $k$ is trivial in this case. 
The $j$ summation is as above, and we obtain
\[
\| f_0\|_{\WNs[\lambda]} \lesssim 
 2^{-\frac{\sigma}2}
  \dfrac{  2^{- h -N h^+}}{ \la h^- \ra^3} \|f\|_{L^2} 
\]

Finally, we consider the last component $f_4$ of $f$.
Owing to the different form of \eqref{kertriprodc}, in this case
we do not have the option of using any local energy 
decay estimate. The first part of $P_{j} P^\lambda_{k} \tilde{P}_h$
is essentially the identity, but does not depend on $\lambda$
so the high modulations are preserved,
\[
\| P_{j} P_{k} \tilde{P}_h f\|_{\Ns_{j}}
 \lesssim   \frac{2^{-h}}{\la h^-\ra} 2^{-\frac{\sigma}{2}}  \|f\|_{L^2} 
\]
 For the second
part we have the product of a time independent  $L^2$ bounded operator 
with $g(\lambda)$, so we can apply the same $L^1L^2$ estimate
as above for low modulations. We obtain
\[
\| P_{j} (P^\lambda_{k} - P_{k})\tilde{P}_h f\|_{\Ns_{j}}
 \lesssim   \frac{2^{-h}}{\la h^-\ra} 2^{-\frac{\sigma}{2}} \frac{2^{-h^+}}{\la h^-\ra^2} \|f\|_{L^2} 
\]
In both cases the $k$ and $j$ summations are trivial, so the proof is concluded.
\end{proof}

\subsection{The time dependent linear flow}

We are now ready to  consider the well-posedness of the equation 
\eqref{wlin-eq-l} in  the smaller space $L X$.

\begin{p1} \label{p:lin-x}
Let $\lambda$ be as in \eqref{lambda-reg}. 
Then the equation \eqref{wlin-eq-l} is well-posed in $L X$,
and the following estimate holds:
\begin{equation}\label{ling}
\| \psi\|_{\WSs[\lambda]} \lesssim \| \psi(0)\|_{L^2} + \|
g\|_{\WNs[\lambda]}
\end{equation}
\end{p1}

\begin{proof} 
For $\psi$ solving \eqref{wlin-eq-l} we consider its dyadic
decomposition
\[
 \psi = \sum_{k} \psi_k,   \qquad 
\qquad \psi_k = P_k^\lambda \psi
\]
To write an equation for $\psi_k$ we use the transference operator 
 described earlier in Proposition~\ref{p:transference},
\[
K_\lambda = \lambda \F_\lambda \frac{d}{d\lambda} \F^*_\lambda 
= - \lambda \frac{d}{d\lambda} \F_\lambda   \F^*_\lambda
\]
For $ P_k^\lambda = \F^*_\lambda \chi_k \F_\lambda $  we compute
\[
 \lambda \frac{d}{d\lambda} P_k^\lambda = 
\F^*_\lambda   [\K,\chi_k] \F_\lambda \psi
\]
Hence the FT of the components $\psi_k$ solve the equations
\[
 (i \partial_t -\tilde H_\lambda) \psi_k = i \frac{\lambda'}{\lambda} \mathcal
F_\lambda^{*}  [\K,\chi_k] \F_\lambda \psi
+g_k
\]
After a further dyadic decomposition on the right,
we obtain the infinite coupled system
\begin{equation}
 (i \partial_t -\tilde H_\lambda) \psi_k = \sum_{j} K_{k j}^\lambda \psi_{j} +
g_k, \qquad 
K_{kj}^\lambda v = i \frac{\lambda'}{\lambda} \F_\lambda^{*}  [\K,\chi_k] \chi_j
\F_\lambda v
\end{equation}
For the left hand side we use Proposition~\ref{p:lin-l2} to treat each
equation in this system in $L^2$ where $\tilde H_\lambda$ can be
viewed as a small perturbation of $\tilde H$.  We claim that the first
term $f_k = \sum_{j} K_{k j}^\lambda \psi_{j}$ in the right hand side
is perturbative,
\begin{equation}\label{sumall}
 \sum_{k < 0} \frac{1}{2^k |k|} \| f_k\|_{\ldNs} 
+ \left(  \sum_{k \geq 0}\| f_k\|_{\ldNs}^2\right)^\frac12 \!\!\!
\ll 
\sum_{k < 0} \frac{1}{2^k |k|} \| \psi_{k}\|_{\ldSs} 
+ \left(  \sum_{k \geq  0} \| \psi_{k}\|_{\ldSs}^2\right)^\frac12
\end{equation}
This is a consequence of the following estimate:
\begin{equation}\label{nosum}
\| K_{kj}^\lambda \psi_j \|_{\ldNs} \lesssim \gamma a_{kj} \| \psi_j \|_{\ldSs} 
\end{equation}
where
\[
a_{kj}=\dfrac{2^{-|j-k|}}{\la j^- \ra \la k^- \ra} (1+ 2^j+2^k)^{-N} \ \mbox{if} \ |k-j| \gg 1, 
\qquad a_{kj}= \frac{1}{\la k^- \ra \la j^- \ra} \ \mbox{if} \ |k-j| \lesssim 1 
\]

It is easy to see that \eqref{nosum} implies \eqref{sumall}.
Harmlessly neglecting the case when either $j > 0$ or $k > 0$, 
where we have rapid decay or $|k-j| \approx 1$ which sums directly,
it suffices to verify that for $j \leq 0$ we have
\[
\sum_{k \leq 0} \frac{1}{2^k \la k \ra} a_{kj} \lesssim \frac{1}{2^j \la j \ra}
\]
But this is straightforward.

It remains to prove \eqref{nosum}.
The Fourier kernel of $K_{kj}^\lambda$ in the $\tilde H_\lambda$ frame 
(i.e. the kernel of $\F_\lambda K_{kj}^\lambda\F_\lambda^*$)
has the form $\lambda' C_{kj}^\lambda$
where
\[
 C_{kj}^\lambda(\xi,\eta,\lambda) = 
\frac{1}{\lambda} {F(\lambda \xi,\lambda\eta)}
\frac{\chi_k(\xi)-\chi_k(\eta)}{\xi^2-\eta^2} \chi_j(\eta)
\]
The simpler case is when $|k-j| \gg 1$; then $C_{kj}^\lambda$ is dyadically
localized in its two arguments at frequency $2^k$, respectively $2^j$, has size
(see Proposition~\ref{p:transference})
\begin{equation}
 |\partial^a_\lambda (\xi \partial_\xi)^b (\eta \partial_\eta)^c 
(\xi \partial_\xi + \eta \partial_\eta)^d C_{kj}^\lambda| \lesssim \frac{2^{-\frac{k+j}2 - |k-j|}}
{\la k^- \ra \la j^- \ra} (1+ 2^k + 2^j)^{-N}.
\label{cmunusize}\end{equation}

The other case is when $|k-j| \lesssim 1$. Then $\eta$ is still localized
at frequency $2^j$ but $\xi$ ranges over all the positive real axis.
In this case we decompose smoothly
\begin{equation}
 C_{kj}^\lambda = \sum_{l} C_{kj}^{\lambda l}, \qquad  C_{kj}^{\lambda l}
(\xi,\eta) = \chi_l(\xi)C_{kj}^\lambda(\xi,\eta) 
\label{cmmdec}\end{equation}
Then $C_{kj}^{\lambda l}$ has the same size and regularity as 
$C_{lj}^\lambda$ above in the case $|l-j| \gg 1$. 

If $|l-j| \les 1$, by Proposition~\ref{p:transference}, we have that
\begin{equation} \label{cmunusized}
\begin{split}
& |\partial^a_\lambda (\xi \partial_\xi)^b (\eta \partial_\eta)^c 
(\xi \partial_\xi + \eta \partial_\eta)^d C_{kj}^{\lambda l}| \les \frac{2^{-k}}{\la k^- \ra^2}, \qquad k \les 0 \\
& |\partial^a_\lambda \partial_\xi^b \partial_\eta^c 
(\xi \partial_\xi + \eta \partial_\eta)^d C_{kj}^{\lambda l}| \les 2^{-2k} \la \xi - \eta \ra^{-N}, \qquad k \ges 0 
\end{split}
\end{equation}
where $a,d,N \in \N$ and $b + c \leq 2$.  

Our main $\ldSs \to \ldNs$ bound is:

\begin{l1}
Let $\lambda$ be as in \eqref{lambda-reg}. Let $K_{kj}^0$  be an
operator whose Fourier kernel in the $\tilde H_\lambda$ frame
has the form $\lambda' C^0_{kj}(\lambda,\xi,\eta)$
where $C^0_{kj}$ is dyadically localized in the region
$\xi \approx 2^k$, $\eta \approx 2^j$ and satisfies the 
uniform bounds \eqref{cmunusize} if $|k-j| \gg 1$ and \eqref{cmunusized} 
if $|k-j| \les 1$, then
\begin{equation}\label{kmn}
\| K_{kj}^0\|_{\ldSs \to \ldNs} \lesssim \gamma a_{kj}.
\end{equation}
\label{lemk}\end{l1}

By \eqref{cmunusize}, this implies \eqref{nosum} directly $|k-j| \gg 1$,
and after an $l$ summation corresponding to the decomposition \eqref{cmmdec}
if $|k-j| \lesssim 1$. The proof of Proposition~\ref{p:lin-x} is concluded.
\end{proof}

\begin{proof}[Proof of Lemma~\ref{lemk}]
The first step in the proof is to reduce the problem 
to the case when the Fourier kernel of $K_{kj}^0$ in the $\tilde H_1$ 
frame is as in the lemma. To switch from the $\tilde H_\lambda$ frame to the 
$\tilde H_1$ frame we write 
\[
K_{kj}^0 = \sum_{k_1,j_1} P_{k_1} K_{kj}^0 P_{j_1}:= \sum_{k_1,j_1}
K_{kj}^{k_1j_1}
\]
where the Fourier kernels $C^{k_1j_1}_{kj}$ 
of $K_{kj}^{k_1j_1}$ are the kernels of
\[
\chi_{k_1}  \F_1 K_{kj}^0 \F_1^* \chi_{j_1}= 
\chi_{k_1} (\F_1 \F_\lambda^*) (\F_\lambda K_{kj}^0 \F_\lambda^*)
(\F_\lambda \F_1^*)  \chi_{j_1}
\]
For the operator $\F_1 \F_\lambda^*$ and its adjoint above we use
Proposition~\ref{p:transfa}. Thus the bounds for $C^{k_1j_1}_{kj}$
are the bounds for $C^0_{kj}$ corrected on the right hand side with the factor
$\alpha_{kj}^{k_1j_1} = \alpha_k^{k_1} \alpha_j^{j_1}$ where
\[
\alpha_{k}^{k_1} = 
\left\{ \begin{array}{ll}
1 & |k-k_1| \lesssim 1 \cr \cr
\dfrac{2^{\frac{k-k_1}2 - |k-k_1|-N(k^++k_1^+)}}{\la k^-\ra \la k_1^- \ra} &
 |k-k_1| \gg 1
\end{array} \right.
\]
Suppose we know that the lemma applies for the operators 
$K_{kj}^{k_1j_1}$ then it is easily seen that it also applies to $K_{kj}^0$
since 
\[
\sum_{k_1,j_1} \alpha_{kj}^{k_1j_1}
\lesssim 1
\]

It remains to prove the lemma in the simpler case where 
the Fourier kernel of $ K_{kj}^0$ is given in the $\tilde H_1$ frame.

If $|k-j| \gg 1$, by separating the variables $\lambda$, $\xi$ and $\eta$ in $C_{kj}^0$
we reduce the problem to the case when the kernel $C^0_{kj}$ has the form
\[
C^0_{kj}(\lambda,\xi,\eta) =   g(\lambda) \chi_k(\xi) \chi_j(\eta)
\]
where $\chi_k$ and $\chi_j$ represent smooth unit bumps
with dyadic localization. With these notations, the operator 
$K^0_{kj}$ takes the form
\[
K^0_{kj} u_{j} = \lambda' g(\lambda) \phi_{k}  \la
\phi_{j},u_j\ra   \qquad \phi_{k} =  \F^* \chi_{k}, \ \phi_{j} =  \F^* \chi_{j}
\]
For $\phi_k$ and $\phi_j$  we have the bound \eqref{bumpft}
which we repeat here for convenience,
\begin{equation}
|\phi_k(r)| \lesssim 2^{\frac{3k}2}
(1+ 2^k r)^{-N}
\label{chimuhat}\end{equation}
 In particular we have the bounds
\begin{equation} \label{phimuinl2}
 \|\phi_k\|_{L^2} \lesssim 2^\frac{k}2, \qquad \| \phi_j \|_{L^2} \lesssim
2^\frac{j}2
\end{equation}
as well as
\begin{equation} \label{phimuinle}
 \|f \phi_k\|_{LE_k} \lesssim 2^{-\frac{k}2}\|f\|_{L^2}, \qquad \| f \phi_j\|_{
LE_k^*} \lesssim 2^{-\frac{j}2}\|f \|_{L^2}
\end{equation}
where $f$ represents a function of time.

If $|k-j| \lesssim 1$ and $k \les 0$, then we separate the variable $\l$ and write
\[
C^0_{kj}(\lambda,\xi,\eta) =   g(\lambda) G(\xi,\eta)
\]
where $G$ is dyadically localized at $\xi \approx 2^k$ and $\eta \approx 2^j$, 
and has size conditions 
\begin{equation}\label{Gdi}
|(\xi \partial_\xi)^\alpha (\eta \partial_\eta)^\beta G(\xi,\eta)| \lesssim \frac{2^{-k}}{\la k^{-} \ra^2},
\qquad \alpha+\beta \leq 2
\end{equation}

If $|k-j| \les 1$ and $k \gg 0$, then the fast decay away from the diagonal allows
us to simplify the problem to the case
\[
C^0_{kj}(\lambda,\xi,\eta) =   g(\lambda) \sum_{} G_{n}(\xi,\eta)
\]
where the sum runs over the positive integers $\approx 2^k$ and $G_n$ satisfies
\begin{equation}\label{Gdii}
|\partial_\xi^\alpha \partial_\eta^\beta G_n(\xi,\eta)| \lesssim 2^{-2k} \la \xi - n \ra^{-N} \la \eta - n \ra^{-N} ,
\qquad \alpha+\beta \leq 2
\end{equation}

Since $\lambda$ belongs to the algebra $\Zu$, we can further simplify the 
expression  $\lambda' g(\lambda)$  occurring in $K^0_{kj}$ and simply 
replace it by $\lambda'$.  By the definition of the space $\Zu$ we
have two possibilities to consider:

{\bf Case A:} $ \lambda \in W^{1,1}$. This is the easier case. Then
$\|\lambda'\|_{L^1} \lesssim \gamma$ therefore by \eqref{phimuinl2}
we obtain
\begin{equation}
 \| K^0_{kj} u_{j}\|_{L^1 L^2} \lesssim \gamma 2^{\frac{k+j}2}
 \| u_j\|_{L^\infty L^2} 
\label{linl1}\end{equation}
which suffices for \eqref{kmn}.

{\bf Case B:} $ \lambda \in \dot H^{\frac12,1}$. In this case we split 
$\lambda$ into a low frequency part and a high frequency part,
\[
 \lambda = \lambda_{\leq m_0+4} + \lambda_{> m_0+4},
\qquad m_0 = 2 \max\{ k,j\}
\]

{\bf Case B1:} The contribution of $ \lambda_{\leq m_0+4}$.
The low frequency part of $\lambda'$ satisfies a favorable $L^2$ bound 
\[
 \| \lambda'_{\leq m_0+4} \|_{L^2} \lesssim \gamma 2^{\frac{m_0}2}
\]
Suppose $|j-k| \gg 1$. Using \eqref{phimuinle} for $\phi_k$,
\eqref{phimuinl2} for $\phi_j$ and the energy of $u_j$ we
obtain
\[
 \| K^0_{k j} u_{j}\|_{l^2 LE^*} \lesssim \gamma 2^{\frac{m_0-k+j}2}
 \|u_{j}\|_{L^\infty L^2} 
\]
which is favorable if $j < k$. In the opposite case $j > k$ we use
\eqref{phimuinle} for $\phi_j$ and \eqref{phimuinl2} for $\phi_k$
 to obtain the better dual type bound
\[
 \| K^0_{k j} u_{j}\|_{L^1 L^2} \lesssim \gamma 2^{\frac{m_0+k-j}2} \|u_{j}\|_{l^2LE} 
\]
Consider now the case $|j-k| \lesssim 1$ and $k \les 0$.  There we have
\[
K^0_{k j} u_{j} = \lambda'_{\leq m_0+4}(t)  f,\qquad 
\hat f(t,\xi) = \int G(\xi,\eta) \hat u_j(t,\eta) d \eta
\]
For $\hat f$ we can use \eqref{Gdi} to estimate up to two derivatives
at fixed time
\[
|(\xi \partial_\xi)^\alpha \hat f(t,\xi)| \lesssim \frac{2^{-\frac{k}2}}{\la k^- \ra^2} \|u_j(t)\|_{L^2}, \qquad 
\alpha \leq 1
\]
Then a variation of Proposition~\ref{p:bumpft}
shows that
\[
 |f(t,r)| \lesssim \frac{2^{k}}{\la k^{-} \ra^2} m_k(r) (1+2^k r)^{-\frac52}  \|u_j(t)\|_{L^2}
\]
which allows us to estimate
\[
\| K^0_{k j} u_{j}\|_{LE^*_k}  \lesssim  \| \lambda'_{\leq m_0+4} \|_{L^2}
\frac{2^{-k}}{\la k^{-} \ra^2} \|u_j(t)\|_{L^\infty L^2} \lesssim \gamma a_{kj} \|u_j(t)\|_{L^\infty L^2} 
\]
In the case $|j-k| \lesssim 1$ and $k \gg 0$ we have
\[
K^0_{k j} u_{j} = \lambda'_{\leq m_0+4}(t) \sum_n f_n,\qquad 
\hat f_n(t,\xi) = \int G_n(\xi,\eta) \hat u_j(t,\eta) d \eta
\]
For $\hat f_n$ we can use \eqref{Gdii} to estimate up to two derivatives
at fixed time
\[
| \partial_\xi^\alpha \hat f_n(t,\xi)| \lesssim 2^{-\frac{3k}2} \la \xi - n \ra^{-N} \|u_j(t)\|_{L^2}, \qquad 
\alpha \leq 2
\]
Then a variation of Proposition~\ref{p:bumpft}
shows that
\[
 |f_n(t,r)| \lesssim m_k(r) (1+2^k r)^{-\frac52}  \|u_j(t)\|_{L^2}
\]
which allows us to estimate
\[
\| K^0_{k j} u_{j}\|_{LE^*_k} \lesssim  \| \lambda'_{\leq m_0+4} \|_{L^2} \sum_n \| f_n \|_{LE^*_k}
 \lesssim  \gamma 2^k \sum_{n} 2^{-2k} \|u_j(t)\|_{L^\infty L^2} \lesssim \gamma \|u_j(t)\|_{L^\infty L^2} 
\]
where we have used that the range of summation has cardinal $\approx 2^k$.

{\bf Case B2:} The contribution of $ \lambda_{m}$, 
$m > m_0+4$. The idea in this case
is that a large modulation for $\lambda$  forces 
either a large modulation in the input or a large modulation in the
output. Precisely,  if $|k-j| \gg 1$ then we have 
\[
\begin{split}
\lambda'_m \phi_k \la \phi_j u_j \ra = &\
Q_{> m-4} ( \lambda'_m \phi_k \la \phi_j, u_j \ra) 
+  \lambda'_m \phi_k \la \phi_j, Q_{>m-4} u_j \ra
\\
& \ - Q_{> m-4} ( \lambda'_m \phi_k \la \phi_j,  Q_{>m-4}u_j \ra) 
\end{split}
\]
The first (as well as the last) term is at high modulation so it
suffices to bound it in $L^2$,
\[
\| Q_{> m-4} ( \lambda'_m \phi_k \la \phi_j, u_j \ra) \|_{L^2}
\lesssim 2^{\frac{k+j}2} \| \lambda'_m\|_{L^2} \|u_j\|_{L^\infty L^2}
\lesssim \gamma 2^{\frac{k+j}2} 2^{\frac{m}2} \|u_j\|_{L^\infty L^2}
\]
On the other hand in the second term the function $u_j$ is restricted to 
high modulations, where we have a good $L^2$ bound:
\[
\|  \lambda'_m \phi_k \la \phi_j, Q_{>m-4} u_j \ra\|_{L^1 L^2}
\lesssim 2^{\frac{k+j}2} \| \lambda'_m\|_{L^2} \| Q_{>m-4} u_j\|_{ L^2}
\lesssim \gamma 2^{\frac{k+j}2} \| u_j\|_{\Ss_j}
\]
A similar argument also applies for $|k-j| \lesssim 1$.
The proof of the lemma is concluded.
\end{proof}

The global in time result in the previous proposition easily implies
its compact interval counterpart:

\begin{c1} \label{c:lin-x}
Let $T > 0$ and $\lambda$ so that
\begin{equation}
\| \lambda -1\|_{\Zu([0,T])} \lesssim \gamma \ll 1  
\end{equation}
Then the solution of \eqref{wlin-eq-l} in $[0,T]$ satisfies:
\begin{equation}\label{linl}
\| \psi\|_{\WSs[\lambda][0,T]} \lesssim \| \psi(0)\|_{L^2} + \|
g\|_{\WNs[\lambda][0,T]}
\end{equation}
\end{c1}

\begin{proof}
Consider an admissible extension $\lambda^{ext}$  for $\lambda$ in $\Zu$,
so that
\[
\| \lambda^{ext} -1\|_{\Zu} \lesssim \gamma \ll 1  
\]
Consider also the zero extension $g^{ext}$ of $g$. This satisfies
\[
\|  g^{ext}\|_{\WNs[\lambda^{ext}]} \lesssim  \|g\|_{\WNs[\lambda][0,T]}
\]
Now solve \eqref{wlin-eq-l} with $\lambda^{ext}$, $g^{ext}$  instead of $\lambda,g$.
By the previous proposition this yields a global solution $\psi$ satisfying
\[
\| \psi\|_{\WSs[\lambda^{ext}]} \lesssim \| \psi(0)\|_{L^2} + \|
g\|_{\WNs[\lambda][0,T]}
\]
The conclusion follows by restricting $\psi$ to the time interval $[0,T]$.
\end{proof}

Proposition~\ref{p:lin-x} also allows us to prove the interval counterpart to
\eqref{wsemb}:

\begin{c1}
Let $I$ be an interval and  $\lambda$
with
\[
\| \lambda -1\|_{\Zu(I)} \lesssim \gamma \ll 1
\]
 Then the
  following inclusion holds:
  \begin{equation}\label{wsembI}
    \WSs[\lambda](I) \subset \WS^r[1](I)
  \end{equation}
\end{c1}

\begin{proof}
Set $I = [a,b]$.
Consider an admissible extension $\lambda^{ext}$  for $\lambda$ in $\Zu$.
Given $\psi \in \WSs[\lambda](I)$, we extend it to the real line 
as a solution to 
\[
(i \partial_t - \tilde H_{\lambda^{ext}}) \psi = 0 \qquad \text{in } \R \setminus I,
\]
matching Cauchy data at $t = a,b$. By the previous corollary we have
\[
\| \psi\|_{ \WSs[\lambda^{ext}](-\infty,a)} +
 \| \psi\|_{ \WSs[\lambda^{ext}](b,\infty)} \lesssim \| \psi(a)\|_{LX} + \|\psi(b)\|_{LX}
\lesssim \|\psi\|_{  \WSs[\lambda](I) }
\]
Due to the Cauchy data matching at $t = a,b$ this implies the global bound
\[
\| \psi\|_{ \WSs[\lambda^{ext}]} \lesssim \|\psi\|_{  \WSs[\lambda](I) }
\]
Now we apply \eqref{wsemb} to obtain 
\[
\| \psi\|_{ \WS^r[\lambda^{ext}]} \lesssim \|\psi\|_{  \WSs[\lambda](I) }
\]
The conclusion follows by restricting the LHS to $I$.
\end{proof}

\subsection{ The autonomous vs nonautonomous flow}

Here, under a suitable $L^2$ smallness condition, we show 
that the solution to the non-autonomous  homogeneous equation 
\begin{equation}
(i \partial_t - \tilde H_\lambda) \psi = 0, \qquad \psi(0) = \psi_0
\label{wlin-h-na}\end{equation}
stays close to the solution of the corresponding autonomous 
 homogeneous equation 
\begin{equation}
(i \partial_t - \tilde H) \tpsi = 0, \qquad \tpsi(0) = \psi_0
\label{wlin-h-a}\end{equation}

\begin{p1} \label{p:comp}
Let $\lambda$ be as in \eqref{lambda-reg}.  Suppose that 
\begin{equation}
\| \psi(0)\|_{LX} \leq 1, \qquad 
\| \psi(0)\|_{L^2} \leq \epsilon \ll 1.
\end{equation}
Then the solutions to \eqref{wlin-h-na} and \eqref{wlin-h-a}
stay close,
\begin{equation}
\| \psi -\tpsi\|_{L^\infty LX} \lesssim |\log \epsilon|^{-1}
\end{equation}
\end{p1}
\begin{proof}
From the $L^2$ bound for the initial data we obtain
\[
\| \psi\|_{S} \lesssim \epsilon, \qquad \| \tpsi\|_{S} \lesssim \epsilon
\]
Because of the $L^2$ bounds above we have 
\[
\| P_{\gtrsim [\log \epsilon]}^\lambda \psi\|_{L^\infty LX} \lesssim |\log \epsilon|^{-1},
\qquad 
\| P_{\gtrsim  [\log \epsilon]} \tpsi \|_{L^\infty LX} \lesssim |\log \epsilon|^{-1}
\]
so it remains to consider the low  frequencies. For $k < \log \epsilon$
we will compare 
\[
\psi_k = P_k^\lambda \psi, \qquad \tpsi_k = P_k \tpsi 
\]
With the notations from the proof of Proposition \ref{p:lin-x}, we have the
following system for $ \{\psi_k\}$ 
\[
(i \partial_t - \tilde H_\lambda) \psi_k = \sum_{j} K_{k j}^\l \psi_j := g_k
\]
By \eqref{nosum}, using the $L^2$ bound for frequencies larger than $\epsilon$ 
and the $LX$ bound for frequencies smaller than $\epsilon$ we obtain 
\begin{equation}
\sum_{k < \log \epsilon} \frac{1}{2^k  |k| } \|g_k\|_{l^2 \Ns} \lesssim 
\frac{1}{|\log \epsilon|}
\label{comp:rhs}\end{equation}

For the initial data we claim to have a similar relation
\begin{equation}
\sum_{k < \log \epsilon} \frac{1}{2^k  |k|} \| \psi_k(0) - \tpsi_k(0)\|_{L^2} 
\lesssim  \frac{1}{|\log \epsilon|}
\label{comp:data1}\end{equation}
To prove this we  write 
\[
\mathcal{F}_{\tilde{H}}  (\psi_k(0) - \tpsi_k(0))
= \sum_{j,h}  ( P_j P_k \tilde P_h  -  P_j P_k^\l \tilde P_h ) P_h \psi(0)
\]
and use Proposition~\ref{p:prodthree} to estimate each term,
\[
\begin{split}
LHS\eqref{comp:data1} \lesssim & \!\!  \sum_{k < \log \epsilon}\!\!  \frac{1}{2^k  |k| } \sum_{j,h}
 c_{jkh} \|P_h\psi(0)\|_{L^2}
\approx \!\! \sum_h \!\! \sum_{k < \log \epsilon} \frac{1}{2^k  |k| } c_{kkh}\|P_h\psi(0)\|_{L^2}
\\
\approx  & \  \sum_h  \sum_{k < \log \epsilon} \frac{1}{2^k  |k| } \frac{2^{-|k-h|-Nh^+}}{|k| \la h \ra}
\|P_h\psi(0)\|_{L^2}
\\
\approx  & \    \sum_{h < \log \epsilon} \frac{1}{2^h  |h^2| }\|P_h\psi(0)\|_{L^2}
+  \sum_{h > \log \epsilon} \frac{1}{|\log \e|} \frac{2^{-Nh^+}}{2^h  \la h \ra }\|P_h\psi(0)\|_{L^2}
\end{split}
\]
Hence \eqref{comp:data1} follows.

Finally we consider the effect of the change in the potential,
\begin{equation}
\|(V_\lambda - V) \psi_k\|_{LE^*} \lesssim \frac{1}{|k|} 
\|\psi_k\|_{LE_k} \lesssim \frac{1}{|k|} 
\|\psi_k\|_{l^2 LE}
\label{comp:data2}\end{equation}
(see \eqref{le=lek} and \eqref{point-le}).
Thus, comparing $\psi_k$ and $\tpsi_k$ along the $\tilde H$ flow we obtain
\begin{equation}
\sum_{k < \log \epsilon} \frac{1}{2^k  |k|} \| \psi_k - \tpsi_k\|_{S}
\lesssim \frac{1}{|\log \epsilon|}
\label{compar}
\end{equation}
We need to turn this into an $L^\infty LX$ bound.
We will use only the $L^\infty L^2$ part of the $S$ norm. At fixed time we write
\[
\| \psi_k - \tpsi_k\|_{LX} \lesssim  \frac{1}{2^k |k|}
\| \tilde P_k ( \psi_k -\tpsi_k)\|_{L^2} + \| (1- \tilde P_k)\psi_k\|_{LX}
\]
For the second term we use  Proposition~\ref{p:transfa}. We obtain
\[
\| \psi_k - \tpsi_k\|_{LX} \lesssim  \frac{1}{2^k |k|}
\|  \psi_k -\tpsi_k\|_{L^2} + \frac{1}{2^k |k|^2}\| \psi_k\|_{L^2}
\]
which combined with \eqref{compar} leads to 
\[
\sum_{k < \log \epsilon}  \| \psi_k - \tpsi_k\|_{LX} \lesssim 
\frac{1}{|\log \epsilon|}
\]
The proof is concluded.

\end{proof}

\section{Analysis of the gauge elements in $X,LX$}
\label{elliptic}

In Section \ref{seccoulomb} we have studied the forward transition from the 
Schr\"odinger map $u$ to its Coulomb gauge $v,w$, its coordinates $\psi_1, \psi_2,A_2$
and finally to its reduced field $\psi$ in the setup where $\|u -Q \|_{\dot H^1}
\leq \gamma \ll 1$, which corresponds to $\psi \in L^2$.  However, the reverse process 
is not uniquely determined in this context. The easiest way to see this
is that if $\psi=0$ then all we can say is that $u$ is one of the solitons $Q_{\alpha,\l}$. 
In some sense, this is the only possible ambiguity.

Here we consider again the transition between $u$ and its reduced
field $\psi$ but in the more regular setting where $\bar u- \bar Q \in X$ and
$\psi \in LX$. An advantage in doing this is that it allows us to 
impose a natural boundary condition at infinity for the system \eqref{comp1},
namely 
\begin{equation} \label{bciX}
 \lim_{r \to \infty} A_2 = 1,  \qquad \psi_2 - i h_1 \in X 
\end{equation}
We will see that on one hand this condition is dynamically preserved 
along the Schr\"odinger map flow, while, on the other hand, it allows for 
an unique identification of $u$ in terms of $\psi$. We remark that, in view 
of \eqref{Xemb}, this condition is satisfied for maps $u$ for which
$u-Q \in L^2$. The main result of this section is as follows:
\begin{t1} \label{thex} a) Let $u : \R^2 \rightarrow \S^2$ be an
  $1$-equivariant map which satisfies $\| u - Q \|_{\dot
    H^1} \ll 1$ and $\| \bar{u} - \bar Q \|_{X} \leq \g \ll 1$.
  Then the Coulomb gauge constructed in Section \ref{seccoulomb}
  satisfies the additional properties
\begin{equation} \label{frinX}
\| \bar v - \bar V \|_{ \tilde X} +  \| \bar w - \bar W \|_{ \tilde X} \leq \g
\end{equation} 
\begin{equation} \label{bfinX}
\| \bar v_3 - \bar V_3 \|_{X} +  \| \bar w_3 - \bar W_3 \|_{X} \leq \g
\end{equation} 
\begin{equation} \label{psi2a2}
\| \psi_2 - ih_1 \|_{X} +  \| A_2 - h_3 \|_{X} \leq \g
\end{equation} 
\begin{equation} \label{psiinX}
\| \psi \|_{LX} \lesssim \g
\end{equation}
Furthermore, the map from $\bar u-\bar Q \in X$ to $\psi \in LX$ is of class $C^1$.

b) Let $\psi : \R^2 \rightarrow \C$ be a function  which satisfies 
$\| \psi \|_{LX} \leq \g \ll 1$.  Then there exists an unique $1$-equivariant 
map $u:\R^2 \to \S^2$ satisfying 
\begin{equation}
\| \bar u - \bar Q \|_X \lesssim \g
\end{equation}
so that $\psi$ is the reduced field for $u$, and the map from
$\psi \in LX$ to $\bar u- \bar Q \in X$ is of class $C^1$. Furthermore, the uniqueness 
of $u$ is also valid in the class of maps $u$ with 
$\|u-Q\|_{\dot H^1} \ll 1$ which satisfy the additional qualitative condition $u-Q \in L^2$. 
\end{t1}

As a quick reminder, $\bar V, \bar W, \bar Q$ were introduced in \eqref{vwq}, 
and $\bar V_3, \bar W_3$ stand for the 
third component of the vectors $\bar V$, respectively $\bar W$.

The plan of this section is as follows. We first prepare for the proof
with an ODE result which will be applied later to the system for the
orthogonal matrix $\calO = (v,w,u)$ in both parts (a) and (b). Then we
prove part (a) in two stages, beginning with the ODE construction of
$\calO$ and continuing with the algebraic derivation of $\psi_2, A_2$ 
and $\psi$. Finally, we prove part (b) also in two stages, namely 
the recovery of $(\psi_2,A_2)$ via the ODE system \eqref{comp1}
with the boundary condition \eqref{bciX}, and then a second ODE 
construction for the matrix $\calO$.

\subsection{ An ODE result}

\begin{l1} \label{lZ}
Consider the ODE
\begin{equation} \label{mp}
\partial_r Z = N Z + F, \qquad \lim_{r \rightarrow \infty} Z(r) =  0
\end{equation}
If $N$ is small in $ \partial_r \tilde X$  then the above equation has a unique solution
$Z \in \tilde X$ satisfying
\begin{equation}
\| Z \|_{\tilde X} \lesssim  \| F \|_{\partial_r \tilde X}
\end{equation}
Furthermore, the map from $N,F \in  \partial_r \tilde X$ to $Z \in \tilde X$ is analytic.
\end{l1}

\begin{proof} The solution $Z$ is obtained via a Picard iteration in the space $\tilde X$. 
 Indeed, the results of Lemma \ref{indX} show that  
\[
\| N Z + F \|_{\partial_r X}  \lesssim 
 \| N \|_{\partial_r \tilde X} \|  Z \|_{\tilde X} + \| F \|_{\partial_r \tilde X} 
\]
and the convergence of the iterations is insured due to the smallness
of $\| N \|_{\partial_r \tilde X}$.  
\end{proof}

\subsection{The transition from $u$ to $(v,w)$} 
\label{utovw}
We use the equation \eqref{cgeq-m} for the matrix $\calO =  ( \bar v, \bar w,\bar u)$,
namely
\begin{equation} \label{sysM}
\partial_r \calO = M(\bar u) \calO, \qquad \calO(\infty) = I_3, \qquad M(\bar u) = \partial_r \bar u \wedge \bar u
\end{equation}
If $u = Q$ them $M(\bar u)$ has the form
\begin{equation}
M(\bar Q) =
\left( 
\begin{array}{lll}
0 & 0 &  -\frac{h_1}r  \\
 0 & 0 & 0 \\
\frac{h_1}r & 0 & 0
\end{array}
\right) \label{E}
\end{equation}
For the difference we claim that
\begin{equation}\label{Mdiff}
\| M(\bar u) - M(\bar Q) \|_{\partial_r \tilde X} \lesssim  
\| \bar{u}-\bar Q \|_X
\end{equation}
Indeed, we write
\[
M(\bar u) - M(\bar Q) = \partial_r ( \bar u -\bar Q) \wedge (\bar u
-\bar Q) + 2 \partial_r ( \bar u -\bar Q) \wedge \bar Q + \partial_r
\left( \bar Q \wedge (\bar u - \bar Q)\right).
\]
For the first term we use \eqref{XdrX} and  for the second we use \eqref{h1X}.
For the third term we have
\[
 \bar Q \wedge (\bar u - \bar Q) = \vec k \wedge (\bar u - \bar Q) 
+ ( \bar Q - \vec k) \wedge (\bar u - \bar Q)
\]
where the first term is trivially in $\tilde X$ while the second belongs to $\He \subset X$
due to the $r^{-1}$ decay of $\bar Q - \vec k$ at infinity. Hence \eqref{Mdiff} is proved.

Returning to \eqref{sysM}, we start with the solution $\calO_0$  for the case 
$u = Q$, which is given by 
\begin{equation}\label{Fsol}
\calO_0=
\left( 
\begin{array}{ccc}
h_3 & 0  & h_1 \\
0  & 1 & 0 \\
-h_1 & 0 & h_3
\end{array}
\right), \qquad
\calO_0^{-1}= \calO_0^t
\end{equation} 
Then we write the solution to \eqref{sysM} is of the form
\begin{equation} \label{Yexpr}
\calO(r) = \calO_0(r)(I+Y(r))
\end{equation}
where $Y$  solves  the differential equation
\begin{equation} \label{Ysysup}
\partial_r Y = N Y + G, \qquad \calO(\infty) = 0 
 \qquad N=G=\calO_0^{-1} (M(\bar u) - M(\bar Q))  \calO.
\end{equation}
The bound \eqref{Mdiff} combined with \eqref{h1X} shows that we can
apply Lemma~\ref{lZ} for $Y$.  The bound \eqref{frinX} follows after
another application of \eqref{h1X}. 

For the extra improvement in \eqref{bfinX}
we still need to estimate $\|r^{-1} \bar v_3\|_{L^2} $ and  $\|r^{-1} \bar w_3\|_{L^2} $.
Consider for instance the latter. Writing
\[
\bar w_3= \bar u_1 \bar v_2 - \bar u_2 \bar v_1= (\bar u_1 -h_1) \bar v_2 
-  \bar u_2 \bar v_1 + h_1 \bar v_2 
\]
the desired bound easily follows.

\subsection{ The transition from $(u,v,w)$ to 
$\psi_2,A_2$ and $\psi$}
By \eqref{psi2vw3}, the bound \eqref{psi2a2} for $\psi_2$ is exactly \eqref{bfinX}, while
by \eqref{a2u3}, the bound \eqref{psi2a2} for $A_2$ follows from the hypothesis. 

It remains to consider $\psi$, which is  represented as 
\[
\psi = \W \cdot v+ i \W \cdot w, \qquad \W=\partial_r \bar u - \frac{1}{r} \bar u \times R \bar u 
\]
In view of the bound \eqref{bfinX} for $(v,w)$ and of the
$LX$ multiplicative estimates   \eqref{LXstable} and \eqref{LXtX}, it suffices 
to show that  
\begin{equation}
\| \W \|_{LX} \lesssim \gamma
\end{equation}
 Since $\W$ vanishes if $u = Q$, we can write
\[
\begin{split}
\W \!= & \partial_r (\bar u - \bar Q) - \frac{1}{r} (\bar u -\bar Q) \times \! R (\bar u - \bar Q) 
-  \frac{1}{r} \bar Q \times R (\bar u - \bar Q) 
-  \frac{1}{r} (\bar u -\bar Q) \times R \bar Q
\\
= & L (\bar u - \bar Q) - \frac{1}{r} (\bar u -\bar Q) \times R (\bar u - \bar Q) 
+ \tilde \W
\end{split}
\]
The first term is in $LX$ by definition and the second belongs to the smaller
space $L^1 \cap L^2$ by  \eqref{pointX} and \eqref{linX}. It remains to consider 
the last component
\[
\tilde \W = - \frac{h_3}{r} (\bar u - \bar Q) -  \frac{1}{r} \bar Q \times R (\bar u - \bar Q) 
-  \frac{1}{r} (\bar u -\bar Q) \times R \bar Q
\]
A direct computation shows that the components of $\tilde \W$ contain the 
expressions $r^{-1} h_1 (\bar u_3 - h_3)$,  $r^{-1} h_1 (\bar u_1 - h_1)$
and  $r^{-1} h_3 (\bar u_3 - h_3)$; we will estimate all of them in $L^1 \cap L^2$.
The $L^2$ bound is obtained directly from the $\dot H^1$ norm of $u-Q$.
The $L^1$ bound for the first two expressions is a consequence of \eqref{linX}.
This also applies to the third expression but only for $r \lesssim 1$
On the other hand for $r \gg 1$  we can use the equation of the sphere 
to obtain
\[ 
|\bar u_3 - h_3| \lesssim (u_1-h_1)^2 + u_2^2 + h_1 |u_1-h_1| 
\]
at which point we can use again \eqref{linX}.

\subsection{The transition from $\psi$ to $(\psi_2,A_2)$.}

This is achieved by solving the ODE system \eqref{comp1} with the boundary condition
\eqref{bciX} at infinity.  We note that by Proposition~\ref{lc}, if $u-Q \in L^2$ then
we have $\psi_2 - i h_1 \in L^2$, which implies \eqref{bciX}.
For convenience we recall
\eqref{comp1} here:
\begin{equation}
 \label{comp1rep}
\left\{ \begin{array}{l}
\partial_r A_2= \Im{(\psi \bar{\psi}_2)}+\frac{1}r |\psi_2|^2, \cr \cr
\partial_r \psi_2 = i A_2 \psi - \frac{1}r A_2 \psi_2
\end{array} \right.
\end{equation}
We are only interested in solutions which belong to the sphere
\begin{equation}
A_2^2 + |\psi_2|^2 = 1
\label{sphere}\end{equation}
A straightforward computation shows that this sphere is invariant with
respect to the \eqref{comp1rep} flow. Thus given any Cauchy data on
this sphere at any point $r_0 \in \R^+$ and any $\psi \in L^2$, there
exists an unique global solution to this ODE. Our challenge here is to
instead prescribe the asymptotic behavior at infinity via
\eqref{bciX}. To achieve this we will take advantage of the additional
information that $\psi \in LX$.  We state our main result here
separately for later use:

\begin{p1} \label{constr}
Assume that $\psi \in LX$, small. Then the system \eqref{comp1rep}
admits a unique solution $(\psi_2,A_2)$  which satisfies \eqref{bciX}.
Furthermore, this solution satisfies the bound
\begin{equation} 
\| \psi_2 - ih_1\|_{X} + \| A_2-h_3\|_{ X} \lesssim \|\psi\|_{LX}, 
\label{compsol}\end{equation}
and it has Lipschitz dependence on $\psi$,
\begin{equation} 
\| \psi_2 -\tilde{\psi}_2\|_{X} + \| A_2-\tilde{A_2}\|_{ X} \lesssim \|\psi -
\tilde{\psi}\|_{LX} .
\label{compsol2}\end{equation}

In addition, the above solution satisfies the following $\dHe$ bounds:
\begin{equation} 
\| \psi_2 -ih_1\|_{\dHe} + \| A_2-h_3\|_{\dHe} \lesssim  \| L^{-1} \psi \|_{\dHe}.  
\label{compsolh1}\end{equation}

\end{p1}
We remark that from \eqref{compsol} and \eqref{sphere} one can get better decay
for $A_2-h_3$ both near $0$ and near infinity.

\begin{proof} We carry out this proof in several steps:

  {\bf Step 1:} Here we assume that a solution $(\psi_2,A_2)$ to
  \eqref{comp1rep} which satisfies the boundary condition \eqref{bciX}
  exists, and we study further its a-priori regularity. In what
  follows $C$ will denote a large positive constant which may vary
  from line to line.

Since $\psi_2- i h_1 \in X$ then, by \eqref{linX} and \eqref{pointX},
 we must have
\begin{equation}
 |\psi_2| \leq  C \la r\ra^{-\frac12}, \qquad \|\frac{r^\frac12}{
\log(1+r)}\psi_2\|_{L^2(dr)} < \infty
\label{psi2i}\end{equation}
and similar bounds for $A_2-h_3$. By virtue of the compatibility relation
\eqref{sphere}, we can improve the bounds for $A_2-h_3$ to
\begin{equation}
 |A_2 - h_3| \leq C \la r\ra^{-1}, 
\qquad \|\frac{r^\frac12\la r\ra^{\frac12}} {\log(1+r)}(A_2 - h_3)\|_{L^2(dr)} <
\infty
\label{A2i}\end{equation}

We rewrite the second equation in \eqref{comp1rep} as
\begin{equation}
L \psi_2 =   i  \psi + f, \qquad
f = i(A_2-1) \psi + \frac{h_3-A_2}{r} \psi_2
\label{eqi}\end{equation}
For large $r$ it suffices to consider this equation, since 
$A_2$ is uniquely determined as $A_2 = \sqrt{1-|\psi_2|^2}$
due to \eqref{sphere}.
Since $\psi \in L^2$, using also \eqref{psi2i} and \eqref{A2i}
we obtain the decay of $f$ at infinity,
\begin{equation}
\| r^{\frac12} \la r \ra f\|_{L^2(dr)} + \|\frac{r \la r\ra^\frac12}{\log(1+r)}
f\|_{L^1(dr)} < \infty
\label{fi}\end{equation}
In particular by \eqref{LXemb} it follows that $f \in LX$. As $\psi_2
- h_1 \in X$, the solution $\psi$ to \eqref{eqi} must have the form
\[
 \psi_2 =i h_1 + i g + \Psi, \qquad g = L^{-1} \psi \in X, \qquad \Psi = L^{-1} f \in X.
\]
Since $f$ has the better decay at infinity given by \eqref{fi}, we can
express $\Psi$ in the integral form
\[
\Psi(r) = - h_1 \int_{r}^\infty h_1(r_1)^{-1} f(r_1) dr_1 
\]
By \eqref{fi} this integral is absolutely convergent, and 
we have the pointwise bound
\begin{equation}
 |\Psi(r)| \leq C \frac{\log (2+r)}{\la r \ra^\frac32}
\label{psili}\end{equation}

We can recast the equation \eqref{eqi} as an equation for $\Psi$,
\begin{equation}
\Psi = N(\Psi,\psi) 
\label{Psig}\end{equation}
where the nonlinear expression on the right has the integral form
\[
 N(\Psi,\psi) = - h_1(r) \int_{r}^\infty  h_1^{-1} 
(i(A_2-1) \psi + \frac{h_3-A_2}{s} (i h_1+ig+\Psi)) ds, 
\]
and $A_2$ and $g$ are dependent variables given by
\[
 A_2 = \sqrt{1-
|ih_1+ig +\Psi|^2}, \qquad g = L^{-1} \psi
\]
We have proved that solving the system \eqref{comp1rep} with the 
boundary condition \eqref{bciX} is equivalent to solving 
the equation \eqref{Psig} for $\Psi \in X$ satisfying 
the decay condition \eqref{psili}.

{\bf Step 2:}
Here we will use the contraction principle to show that for $r$
near infinity, $r \in [R,\infty)$ 
there exists a unique solution $\Psi$ to  \eqref{Psig}
satisfying \eqref{psili}, which depends in a Lipschitz manner on 
$\psi \in LX$. For this we will prove that the nonlinearity $N$ satisfies the Lipschitz bound
\begin{equation}
\| N(\Psi,\psi) -  N(\tilde \Psi,\tilde{\psi})
\|_{L^\infty_{r^\frac54}([R,\infty))} \leq \frac{C_0}{ R^{\frac14}}
( \|\psi-\tilde{\psi} \|_{LX} +  \|\Psi-\tPsi
\|_{L^\infty_{r^\frac54}([R,\infty))} ) 
\label{lipPsig}\end{equation}
where
\[
C_0 = C_0(\| \psi \|_{LX}, \| \tilde{\psi}\|_{LX}, \|\Psi\|_{
L^\infty_{r^{\frac54}}[R,\infty)} ,
\|\tPsi\|_{ L^\infty_{r^{\frac54}}[R,\infty)})
\]
For the existence and Lipschitz dependence part we start with $\psi$
satisfying $\|\psi\|_{X} \leq 1$, choose $R$ so that 
\[
   R^{-\frac14} C_0(1,1,1,1) \leq \frac14
\]
and apply the contraction principle in the unit ball in 
$L^\infty_{r^{\frac54}}[R,\infty)$. This yields a unique solution $\Psi$
which satisfies 
\begin{equation}\label{Psiex}
\|\Psi\|_{L^\infty_{r^\frac54}([R,\infty))} \lesssim  \| \psi \|_{LX}
\end{equation}
with the Lipschitz dependence 
\begin{equation} \label{Psilip}
\|\Psi-\tilde \Psi\|_{L^\infty_{r^\frac54}([R,\infty))} \lesssim  
\| \psi - \tilde \psi \|_{LX}
\end{equation}
For the uniqueness part we use \eqref{lipPsig} but with a 
larger $\tilde R$ chosen so that 
\[
 \tilde R^{-\frac14} C_0(\| \psi \|_{LX}, \| \tilde{\psi}\|_{LX}, \|\Psi\|_{
L^\infty_{r^{\frac54}}[R,\infty)} ,
\|\tPsi\|_{ L^\infty_{r^{\frac54}}[R,\infty)}) < 1
\]
It suffices to prove uniqueness on a smaller interval $[\tilde R,\infty)$
since the equation \eqref{Psig} is equivalent to the original ODE system 
\eqref{comp1rep}, for which uniqueness holds for $r$ in a compact interval in
$(0,\infty)$.

We now continue with the proof of \eqref{lipPsig}. From the
formulas for $A_2, \tilde{A}_2$, we have 
\[
 |A_2 -\tilde A_2| \lesssim r^{-\frac12} (|g -\tg|+ |\Psi-\tPsi|) 
\]
\[
|A_2 -1| \lesssim h_1^2 + r^{-\frac12}(|g| + |\Psi|),
 \quad |A_2-h_3| \lesssim r^{-\frac12}(|g| + |\Psi|)
\]
which implies that (recall the definition of $f$ from \eqref{eqi})
\begin{equation}
\begin{split}
 |f - \tilde f|  \lesssim \ & |A_2-\tilde{A_2}| |\tilde{\psi}| +
 |A_2-1||\psi-\tilde{\psi}|  \\  & + r^{-1} |A_2-h_3| |\psi_2 -\tilde{\psi}_2|+
 r^{-1} |A_2-\tilde{A}_2||\tilde{\psi}_2| \\
  \lesssim & \ \delta f_1 + \delta f_2 + \delta f_3
\end{split} \label{deltaf}
\end{equation}
where
\[
 \delta f_1 = r^{-\frac12} |\psi| (|g -\tg|+ |\Psi-\tPsi|),
\qquad \delta f_2 = (h_1^2
+ r^{-\frac12}(|\tilde{g}| + |\tilde{\Psi}|)) |\psi - \tilde{\psi}| 
\]
\[
  \delta f_3 =   r^{-2}(|g -\tg|+|\Psi-\tPsi|)
\]
Correspondingly we derive a bound for $ N(\Psi,\psi) -  N(\tilde \Psi,\tilde{\psi})$,
\[
|N(\Psi,\psi) -  N(\tilde \Psi,\tilde{\psi})| \lesssim\delta N_1 + \delta N_2 +
 \delta N_3, \ \ \  \delta N_i(r) =  h_1(r) \int_r^\infty h_1^{-1}(s) \delta f_i(s)  ds
\]

Now we successively consider the three contributions.
For $\delta f_3$ by \eqref{pointX}  we have the pointwise bound
\[
|\delta f_3(r)| \lesssim r^{-\frac52} (\|g
-\tg\|_{X}+\|\Psi-\tPsi\|_{L^\infty_{r^\frac54}([R,\infty))}) 
\]
therefore its  contribution $\delta N_3$ satisfies:
\[
\begin{split}
\delta N_3(r) \lesssim r^{-\frac32} (\|g
-\tg\|_{X}+\|\Psi-\tPsi\|_{L^\infty_{r^\frac54}([R,\infty))}) 
\end{split}
\]

For $\delta f_2$ we use instead \eqref{linX}, 
to get an $L^1$ bound
\[
\| r^\frac12 |\log r|^{-1} \delta f_2\|_{L^1(rdr)} \lesssim  
(\|\tg\|_X+\|\Psi\|_{L^\infty_{r^\frac54}([R,\infty))}) 
\| \psi - \tilde{\psi} \|_{L^2}
\]
% Then we use the trivial estimate
% \[
% \int_r^\infty s^{-\alpha} |h(s)| ds \lesssim r^{-\alpha} \| h \|_{L^2(rdr)}, \quad \alpha > 0
% \]
which leads to
\[
\delta N_2(r) \lesssim r^{-\frac32}|\log r|^{-1} 
(\|\tg\|_X+\|\Psi\|_{L^\infty_{r^\frac54}([R,\infty))}) 
\| \psi - \tilde{\psi} \|_{L^2}
\]

Finally, for $\delta f_1$ we use again \eqref{linX} to estimate 
\[
\|r^{\frac12} |\log r|^{-1} \delta f_1\|_{L^1(rdr)} \lesssim (\|g-\tilde{g}\|_{X} + \| \Psi - \tilde{\Psi}
\|_{L^\infty_{r^\frac54}([R,\infty))}) \|\psi\|_{L^2}
\]
which yields
\[
\delta N_1(r) \lesssim r^{-\frac32} |\log r|^{-1} (\|g-\tilde{g}\|_{X} + \| \Psi - \tilde{\Psi}
\|_{L^\infty_{r^\frac54}([R,\infty))}) \|\psi\|_{L^2}
\]
The proof of \eqref{lipPsig} is concluded. 

{\bf Step 3:} Now we consider the solution $\Psi$ to \eqref{Psig}  obtained in
the previous step in the interval $[R,\infty)$, and we supplement 
the pointwise bounds \eqref{Psiex} and \eqref{Psilip} with 
$\He$ bounds
\begin{equation}\label{PsiexH1}
\|\Psi\|_{\He[R,\infty)} \lesssim  \| \psi \|_{LX}
\end{equation}
\begin{equation} \label{PsilipH1}
\|\Psi-\tilde \Psi\|_{\He[R,\infty)} \lesssim  \| \psi - \tilde \psi \|_{LX}
\end{equation}
In view of the embedding $\He \subset X$, these will be useful later
to establish the Lipschitz dependence in $X$. Returning to 
$(\psi_2,A_2)$ these bounds imply that
\begin{equation}\label{aaa1}
\|(\psi_2- g)- (\tilde \psi_2 - \tilde g)\|_{\He[R,\infty)} 
\lesssim  \| \psi - \tilde \psi \|_{LX}
\end{equation}
respectively
\begin{equation}\label{aaa2}
\| A_2- \tilde A_2 \|_{\He[R,\infty)} 
\lesssim  \| \psi - \tilde \psi \|_{LX}
\end{equation}

The $L^2$ part of
\eqref{PsiexH1} and \eqref{PsilipH1} follow trivially from
\eqref{Psiex} and \eqref{Psilip}. Consider now the estimate for the
$L^2$ norm for $\partial_r \Psi$.  Given that we already have an $L^2$
bound for $\Psi$, we can freely replace $\partial_r \Psi$ by $L\Psi = f$. It
remains to show that
\begin{equation} \label{l2f}
\| f-\tilde{f} \|_{L^2} \lesssim \|\psi-\tilde{\psi}\|_{LX}+
\|\Psi-\tilde{\Psi}\|_{L_{r^\frac54}^\infty([R,\infty))}  
\end{equation}
Consider the bound \eqref{deltaf} for $\| f-\tilde{f} \|_{L^2}$.  The
estimate for $\delta f_3$ proved in the previous step is already good
enough. We only need to revisit the bounds on $\delta f_2$ and $\delta
f_1$, which we do by using \eqref{pointX} instead of \eqref{linX}.
For $\delta f_2$ we obtain
\[
\| r \delta f_2\|_{L^2(rdr)} \lesssim  
(\|\tg\|_X+\|\Psi\|_{L^\infty_{r^\frac54}([R,\infty))}) 
\| \psi - \tilde{\psi} \|_{L^2}
\]
Similarly for $\delta
f_1$ we have 
\[
\| r \delta f_1\|_{L^2(rdr)}
\lesssim (\|g-\tilde{g}\|_{X} + \| \Psi - \tilde{\Psi}
\|_{L^\infty_{r^\frac54}([R,\infty))}) \|\psi\|_{L^2}
\]
Both bounds are much stronger than we need.

{\bf Step 4:} Here we prove the large $r$ part of \eqref{compsolh1},
which with respect to $\Psi$ takes the form
\begin{equation} \label{claimN1}
\| N(\Psi,\psi) \|_{\dHe([R,\infty))} \lesssim \| g \|_{\dHe}
\end{equation}
This  implies 
that the solution $\psi_2$ constructed above in $[R,\infty)$ satisfies
\begin{equation} \label{hdot11}
\| \psi_2 -i h_1 \|_{\dHe([R,\infty))}+
\| A_2 -h_3 \|_{\dHe([R,\infty))}   \lesssim \| g \|_{\dHe}
\end{equation}
In proving \eqref{claimN1} we can assume that the following bounds are 
valid:
\[
\| g\|_{X} + \|\Psi\|_{L^\infty_{r^\frac54}([R,\infty))} \lesssim 1
\]

To establish \eqref{claimN1}, we use the following pointwise estimate
on $f$: 
\[
\begin{split}
|f| & \lesssim |A_2-1||\psi| + \frac{|h_3-A_2|}{r}||\psi_2| \\
& \lesssim (h_1^2+r^{-\frac12}(|g|+|\Psi|)) |\psi| + r^{-\frac32}(|g|+|\Psi|)
 (h_1 + |g| + |\Psi|)
\end{split}
\]
Then for $\Psi=N(\Psi, \psi)= L^{-1} f$, using\eqref{linX}, we obtain
\[
\begin{split}
\| \frac{\Psi}{r} \|_{L^2([R,\infty))}  \lesssim &  \int_R^{\infty} s |f(s)| ds \\
\lesssim & \int_R^\infty (s^{-2}+s^{-\frac12}(|g(s)|+|\Psi(s)|)) |\psi(s)| sds \\
& + \int_R^\infty s^{-\frac32}(|g(s)|+|\Psi(s)|) (h_1 + |g(s)| + |\Psi(s)|)  sds \\
 \lesssim & \  R^{-\frac12 + \epsilon} (\| \psi \|_{L^2} + \| \frac{g}{r} \|_{L^2} + \| \frac{\Psi}{r} \|_{L^2([R,\infty))})
\end{split}
\]

By taking $R$ large and using $\| \psi \|_{L^2} \lesssim \| g \|_{\dHe}$  we obtain
\[
\| \frac{\Psi}{r} \|_{L^2([R,\infty))} \lesssim \| g \|_{\dHe}
\]
From this estimate and the above pointwise bound for $f$, 
it also follows that
\[
\| f \|_{L^2([R,\infty))} \lesssim \| \psi \|_{L^2} + \| \frac{g}{r} \|_{L^2}+\| \frac{\Psi}{r} \|_{L^2([R,\infty))} \lesssim \| g \|_{\dHe}
\]
Finally from the last two estimates we obtain
\[
\| \partial_r \Psi \|_{L^2([R,\infty))} \lesssim \| g \|_{\dot{H}^1} + \| f \|_{L^2} \lesssim \| g \|_{\dHe}
\]
which concludes the proof of \eqref{claimN1} and of the full characterization
of $\Psi$ on $[R,\infty)$.

{\bf Step 5: The bounds for $\psi_2, A_2$ on $I = [R^{-1},R]$.}
On the interval $[R^{-1},R]$ we can no longer 
use only the second equation in \eqref{comp1rep}. However, in 
this interval there is no singularity so a standard ODE analysis
allows us to extend the solution with Lipschitz pointwise bounds.
Precisely, a straightforward application of Gronwall's inequality
shows that as long as $\|\psi\|_{L^2}, \|\tilde \psi \|_{L^2} \ll 1$
we have
\[
\begin{split}
\| (\psi_2,A_2) - (\tilde \psi_2,\tilde A_2)\|_{L^\infty(I)}
 \lesssim  
\| \psi-\tilde \psi\|_{L^2(I)}
  + | (\psi_2,A_2)(R) - (\tilde \psi_2,\tilde A_2)(R)| 
\end{split}
\]
Reusing this  in \eqref{comp1rep} we can also estimate the  
$r$ derivatives,
\[
\| \partial_r (\psi_2,A_2) -\partial_r (\tilde \psi_2,\tilde A_2)\|_{L^2(I)}
 \lesssim  
\| \psi-\tilde \psi\|_{L^2(I)}
  + | (\psi_2,A_2)(R) - (\tilde \psi_2,\tilde A_2)(R)| 
\]
Estimating the second term on the right by \eqref{aaa1} and  \eqref{aaa2} we obtain
\begin{equation} \label{aaa3}
\| (\psi_2,A_2) - (\tilde \psi_2,\tilde A_2)\|_{\He(I)}
 \lesssim  
\| \psi-\tilde \psi\|_{LX}
\end{equation}
If $\tilde \psi = 0$ then using \eqref{hdot11} instead we get
\begin{equation} \label{aaa3bis}
\| (\psi_2,A_2) - (ih_1,h_3)\|_{\He(I)}
 \lesssim  
\| L^{-1} \psi \|_{\dHe}
\end{equation}

{ \bf Step 6 : The bounds for $(\Psi_2,A_2)$ on $(0, R^{-1}]$.}  On
the interval $(0,R^{-1}]$, $A_2$ is expected to be negative, so we can
use again only the second equation in \eqref{comp1rep} with $A_2 = -
\sqrt{1-|\psi_2|^2}$.  We repeat the fixed point argument as we did on
$[R,\infty)$. We rewrite the second equation in \eqref{comp1rep} as
\begin{equation}
L \psi_2 =  - i  \psi + f, \qquad
f = i(A_2+1) \psi + \frac{h_3-A_2}{r} \psi_2
\label{eqi2}\end{equation}
We introduce $\psi_2 - i h_1+ i g = \Psi$ and rewrite the problem as 
\begin{equation} \label{prob}
\Psi(r) =  \frac{h_1(r)}{h_1(R^{-1})} \Psi(R^{-1}) + N(\Psi,\psi)
\end{equation}
where the nonlinearity $N$ is defined as
\[
 N(\Psi,\psi) = - h_1(r) \int_{r}^{R^{-1}}  h_1(s)^{-1} f(s) ds, \qquad f = 
i(A_2+1) \psi + \frac{h_3-A_2}{s} (ih_1- ig+\Psi)
\]
with $A_2$ and $g$ as dependent variables,
\[
 A_2 = -\sqrt{1-|h_1-ig +\Psi|^2}, \qquad g = L^{-1} \psi.
\]
The value for $\Psi(R^{-1}) = (\psi_2 - ih_1 - i g)(R^{-1})$ is
collected from  Step 5 and satisfies
\begin{equation}
|\Psi(R^{-1})-\tilde\Psi(R^{-1})|  \lesssim \| \psi -\tilde \psi\|_{LX} \ll 1,
\qquad  \Psi(R^{-1}) \lesssim \|g\|_{\dHe} \ll 1.
\label{Psidata}\end{equation}
For this new nonlinearity $N$ we claim the following bound:
\begin{equation}
\| N(\Psi,\psi)-N(\tPsi,\tilde{\psi}) \|_{\dHe((0,R])} \leq \frac12 ( \|
\Psi -\tPsi\|_{\dHe(0,R]} +
\| g - \tilde{g} \|_{\dHe(0,R]})
\label{lipPsig2}\end{equation}
under the assumption that 
\begin{equation} \label{ass2}
 \| \psi \|_{LX}, \| \tilde{\psi} \|_{LX}, \|\Psi\|_{\dHe(0,R]} ,
\|\tPsi\|_{\dHe(0,R]}, R^{-1} \ll 1
\end{equation}

As in the large $r$ case, the Lipschitz bound \eqref{lipPsig2} and
\eqref{Psidata} allows us use the contraction principle to obtain a
solution $\Psi \in \dHe(0,R^{-1}]$ to \eqref{prob} satisfying
\[
\| \Psi \|_{\dHe((0,R])} \lesssim \| L^{-1} \psi \|_{\dHe}, \quad \| \Psi
-\tilde{\Psi} \|_{\dHe((0,R])} 
\lesssim \| \psi-\tilde{\psi} \|_{LX}
\]
Returning to $\psi_2$ and $A_2$ this gives
\begin{equation} \label{aaa4}
\| \psi_2 - i h_1 \|_{\dHe((0,R])} \lesssim \| L^{-1} \psi \|_{\dot{H}^1}, \quad \| \psi_2
-\tilde{\psi}_2 \|_{\dHe((0,R])} 
\lesssim \| \psi-\tilde{\psi} \|_{LX}
\end{equation}
\begin{equation} \label{aaa5}
\| A_2 - h_3 \|_{\dHe((0,R])} \lesssim \| L^{-1} \psi \|_{\dHe}, 
\quad \| A_2 -\tilde{A}_2 \|_{\dHe((0,R])} 
\lesssim \| \psi-\tilde{\psi} \|_{LX}
\end{equation}

It remains to establish \eqref{lipPsig2}. We start with the inequalities
\[
|A_2 - \tilde{A}_2| \lesssim |g-\tg| + |\Psi - \tPsi|, \quad
|A_2+1| + |A_2-h_3| \lesssim r + |g| + |\Psi|
\]
which are derived from \eqref{ass2} and the formulas from $A_2, \tilde{A}_2$.
From these estimates we derive a pointwise bound for $f$,
\[
\begin{split}
|f-\tilde{f}| & \lesssim (|g-\tg| + |\Psi - \tPsi|)|\psi| + (r + |\tg| +
|\tilde{\Psi}|)|\psi-\tilde{\psi}| \\
& + (r +|g|+|\Psi| + |\tg| + |\tilde{\Psi}|)\frac{|g-\tg| + |\Psi
- \tPsi|}{r}
\end{split}
\]
This directly leads to the $L^2$ bound 
\[
\| f - \tilde f\|_{L^2} \lesssim  
C ( \|g-\tg \|_{\dHe((0,R])} +  \|\Psi-\tPsi \|_{\dHe((0,R])} )
\]
\[
C = R^{-1} + \|g\|_{\dHe} + \| \Psi \|_{\dHe((0,R])}+ \|\tilde{g}\|_{\dHe} + \| \tilde{\Psi} \|_{\dHe((0,R])}
\]
For small $r$ we have $h_1(r) \sim r$ therefore
\[
|N(\Psi,\psi) - N(\tilde \Psi,\tilde \psi)| \lesssim r [r \partial_r]^{-1}|f -\tilde f|
\]
Hence combining the above $L^2$ bound for $f -\tilde f$ with the Hardy
estimate \eqref{rdrm} and with \eqref{ass2} we obtain
\[
\begin{split}
 \| r^{-1}(N(\Psi,g)-N(\tilde \Psi,\tilde g)) \|_{L^2((0,R])} 
 \lesssim C ( \|g-\tg \|_{\dHe((0,R])} +  \|\Psi-\tPsi \|_{\dHe((0,R])} )
\end{split}
\]
Finally, using
\[
\partial_r N(\Psi,g)= -\frac{h_3}{r} N(\Psi,g) + f
\]
we also bound $\partial_r N(\Psi,g)$ in $L^2$, completing the proof of 
 \eqref{lipPsig2}.

{\bf Step 7: Conclusion} 

In the end, based on  \eqref{aaa1}, \eqref{aaa2}, \eqref{aaa3}, \eqref{aaa4},
\eqref{aaa5} and \eqref{Xembt},
we upgrade the solution constructed above to $\psi_2 - i h_1, A_2 - h_3 \in X$
with the bounds \eqref{compsol}-\eqref{compsol2}. In addition \eqref{compsolh1} follows
from \eqref{claimN1}, \eqref{aaa3bis} and \eqref{aaa4}.

\end{proof}

\subsection{The transition from $\psi,\psi_2$ 
and $A_2$ to $(u,v,w)$}
 
To achieve this we use the system  \eqref{return}, 
which we recast in a matrix form as an equation 
for $\calO=(v,w,u)$ as follows
\begin{equation} \label{sysM2}
\partial_r \calO = \calO R(\psi) 
\end{equation}
with 
\[
R=\left( 
\begin{array}{lll}
0 & 0 &   \Re \psi_1 \\
0 & 0 &  \Im \psi_1 \\
- \Re \psi_1 &  - \Im \psi_1 &  0
\end{array}
\right)
\]
If $\psi = 0$ then $\psi_2 = i h_1$, which yields $\psi_1= - \frac{h_1}{r}$.
Hence $R(0) = M(\bar Q)$ as in \eqref{E}. We will prove that
\begin{equation}\label{Rest}
\| R(\psi) - R(0)\|_{\partial_r \tilde X} \lesssim \|\psi\|_{LX}
\end{equation}
Suppose this is done. Then the same argument as in Section~\ref{utovw}
leads to the bound \eqref{frinX}, as well as
\[
\|u - Q\|_{\tilde X} \lesssim \gamma
\]
To upgrade the above norm to an $X$ norm we need an additional bound
for $\|r^{-1}(u-Q)\|_{L^2}$.  We first remark that the last row of
$\calO$ is a-priori known, namely $(v_3,w_3,u_3) = (\Im \psi_2,\Re
\psi_2,A_2)$; this already shows that 
\[
\| r^{-1}(v_3 - h_1)\|_{L^2} + \| r^{-1} w_3\|_{L^2} + \| r^{-1}(u_3 - h_3)\|_{L^2}  
\lesssim \gamma
\]
To transfer this information to $u_1$ and $u_2$ we use again the orthogonality
of $\calO$. For $u_1$ for instance we have 
\[
u_1 = v_2 w_3 - v_3 w_2 = v_2 w_3 - (v_3-h_1) w_2 - h_1 (w_2-1)+ h_1   
\]
which suffices.

It remains to prove the bound \eqref{Rest}. Using the second relation in 
\eqref{comp1rep} we have
\[
\begin{split}
\psi_1 = &  \psi + i\frac{\psi_2}{r} = -i A_2 \partial_r \psi_2 
+ |\psi_2|^2  \psi + i  \frac{\psi_2}{r} |\psi_2|^2
\\
= & -\frac{h_1}{r} -i A_2 \partial_r (\psi_2 - i h_1) + (A_2-h_3) \partial_r h_1
+  |\psi_2|^2  \psi + \frac1r (  i  \psi_2 |\psi_2|^2 + h_1^3)
\end{split}
\]
The first term is the value that corresponds to $\psi = 0$. 
The second is placed in $\partial_r \tilde X$ by \eqref{XdrX} and \eqref{h1X}.
The remaining terms are estimated in $\partial_r \tilde X$ just based on their size,
via \eqref{drtXsize}. The third term is pointwise bounded by $\gamma \la r \ra^{-3}$.
For the third one we use  the $L^2$ bound for $\psi$, combined with 
the $l^2 L^\infty$ bound on $\psi_2-ih_1$ for small $r$ and the pointwise 
$r^{-\frac12}$ bound on $\psi_2-ih_1$ for large $r$. The fourth one is similar, only 
the  $L^2$ bound for $\psi$ is replaced by the  $L^2$ bound for $r^{-1}(\psi_2-ih_1)$.

\subsection{ Local energy bounds} 
A key role in the study of the Schr\"odinger type equation \eqref{psieqa}
for $\psi$ is played by the dispersive estimates for $\psi$, most notably the
local energy decay, which allows us to control a norm for $\psi$ which 
is of the form
\[
\|\psi\|_{LE[\lambda]} = 
\sum_{k <  0}\frac{1}{k 2^k} \|P_k^\lambda \psi\|_{LE_k} + \left(
\sum_{k \geq  0} \|P_k^\lambda \psi\|_{LE_k}^2\right)^\frac12
\]
where $\lambda$ is a function time for which
\begin{equation} \label{laest}
|\l -1 | \ll 1
\end{equation}
This is always satisfied in the context of this paper, as the 
functions $(\psi_2,A_2)$ given by Proposition~\ref{constr} satisfy
\begin{equation} \label{laesta}
|\l -1 | + |\alpha| \lesssim |A_2(1) - h_3(1)| + |\psi_2(1) - h_1(1)|  \lesssim \| \psi
\|_{LX}
\end{equation}
A-priori this norm depends on the choice of $\lambda$. However,
using Proposition~\ref{p:transfa} and Proposition~\ref{p:bumpft}
it is easy to prove that different choices of $\lambda$ subject to
\eqref{laest} yield equivalent norms.

 In this section we study to what extent the 
local energy decay bounds for $\psi$ can be transferred to $(\psi_2,A_2)$
via the system \eqref{comp1rep}. 
At first one might attempt to prove local
energy decay bounds for $\psi_2-ih_1$ and $A_2-h_3$.  
If that were true,
it would imply square integrability for $\lambda(t)-1$ and $\alpha(t)$,
where $\lambda(t)$ and $\alpha(t)$ are the parameters defined in 
\eqref{analdef} describing $(\psi_2,A_2)$ at $r = 1$.
However, such decay estimates turn out not to hold. 

Our remedy for this difficulty is to start with $\lambda$ and $\alpha$
defined in \eqref{analdef} and to compare $(\psi_2,A_2)$ with their
value associated to the harmonic map $Q_{\alpha,\lambda}$. Precisely,
with $\lambda$ and $\alpha$ given by
\begin{equation} \label{defal}
A_2(1)=h_3^\l(1), \psi_2(1)=i e^{i\alpha} h_1^\l(1)
\end{equation}
we seek to estimate the differences
\begin{equation} \label{defdelta}
\delta^{\l,\alpha} \psi_2=\psi_2 - i e^{i\alpha}h_1^\l, \qquad\delta^\l A_2=A_2 - h_3^\l. 
\end{equation}
For $\lambda$ and $\alpha$ in \eqref{defal} we assume that 
\begin{equation} \label{allaest}
\| \alpha \|_{L^\infty} + \|\l -1 \|_{L^\infty} \ll 1
\end{equation}
In the context of Proposition~\ref{constr} this is a consequence 
of the bound
\[
\| \psi\|_{L^\infty LX} \ll 1
\]
The main  result of this section is the following

\begin{p1} \label{p:leA2}
a)  Suppose that  $\psi \in L^2$, small. Let 
$(\psi_2,A_2)$  be the solutions to \eqref{comp1rep} with initial data 
as in \eqref{defal}, \eqref{allaest}. Then we have the fixed time bound
\begin{equation}
 \| \delta^{\l,\alpha} \psi_2\|_{\dHe} +  \| \delta^\l A_2\|_{\dHe} 
\lesssim \|\psi\|_{L^2}
\label{h1fowpsi2}\end{equation}

b) Assume in addition that $\psi$ is small in $L^\infty LX$ and that
\eqref{allaest} is valid.  Then the following space-time bound holds:
\begin{equation}
\|\frac{\langle r \rangle^{ -\epsilon}}{r} \delta^{\l,\alpha} \psi_2 \|_{L^2} +
\|\frac{ \langle r \rangle^{\frac12-\epsilon}}{r} \delta^\l A_2\|_{L^2}\lesssim  
\| \psi \|_{LE[1]} 
\label{dap-eloc}\end{equation}
\end{p1}
We remark that heuristically \eqref{h1fowpsi2} 
can be viewed as a consequence of the estimate 
\eqref{goodal} and the relation \eqref{basicpsi}. 

\begin{proof} 
a) By \eqref{allaest}, $\lambda$ is close to $1$. Solving 
\eqref{comp1rep} on the time interval $[1,\lambda]$ we 
obtain the bound
\begin{equation} \label{switchr}
 | \delta^{\l,\alpha} \psi_2(\lambda)| +  | \delta^\l A_2(\lambda)| \lesssim \|\psi\|_{L^2_{comp}}
\end{equation}
Then we can use  a rotation and scaling to set  
 $\lambda = 1$ and $\alpha = 0$ in \eqref{defdelta}
at the expense of replacing \eqref{h1fowpsi2} by 
\begin{equation}
 \| \delta \psi_2\|_{\dHe} +  \| \delta A_2\|_{\dHe} 
\lesssim \|\psi\|_{L^2} +  | \delta \psi_2(1)| +  | \delta A_2(1)|
\label{h1fowpsi3}\end{equation}
under a smallness assumption on the right hand side. Here we make the 
convention that if $\l=1,\alpha=0$ then we drop the upper-scripts from $\delta$.
Using \eqref{sphere} we rewrite the equation \eqref{comp1rep} in the 
equivalent form
\begin{equation}
 \label{comp1rep1}
\left\{ \begin{array}{l}
L \delta \psi_2 = i A_2 \psi - \frac{1}r \delta A_2 \psi_2, \cr \cr
L_1 \delta  A_2= \Im{(\psi \bar{\psi}_2)}+ \frac{1}r |\delta A_2|^2
\end{array} \right.
\end{equation}
where the operators $L$ and $L_1$ are given by
\[
L = \partial_r + \frac{h_3}{r}, \qquad L_1 = \partial_r + 2\frac{h_3}{r}
\]
The functions $h_1$, respectively $h_1^2$ solve  the homogeneous
equations $Lf = 0$, respectively $L_1 f = 0$.
Then the inverses $T$, respectively $T_1$  of $L$, respectively $L_1$ 
with zero Cauchy data at $r=1$ have the form
\[
 T f (r) = h_1(r) \int_{1}^r \frac{f(s)}{h_1(s)} ds, \qquad
 T_1 f (r) = h_1^2(r) \int_{1}^r \frac{f(s)}{h_1^2(s)} ds
\]
Then we have:

\begin{l1} \label{ldecay}
The operators $T$ and $T_1$ satisfy the bounds
\[
 \| r^{\alpha} T f\|_{\dHe(1,r_0)} 
  \lesssim \| r^{\alpha} f\|_{L^2(1,r_0)} 
\]
where the range of $\alpha$ is $\alpha < 1$ if $r_0 > 1$,
respectively $\alpha > -1$ if $r_0 < 1$. 
\end{l1}
The proof of the lemma is straightforward, and is left for the reader.
To continue with the proof of the proposition we rewrite
\eqref{comp1rep1} as 
\begin{equation}
 \label{comp1rep2}
\left\{ \begin{array}{l}
 \delta \psi_2 = h_1 \delta \psi_2(1) + T (
 i A_2 \psi - \frac{1}r \delta A_2 \psi_2), \cr \cr
 \delta  A_2=h_1^2 \delta A_2(1)+ T_1 ( \Im{(\psi \bar{\psi}_2)}+ \frac{1}r |\delta A_2|^2)
\end{array} \right.
\end{equation}
and solve this equation using the contraction principle in 
$\dHe \times \dHe$. Given Lemma~\ref{ldecay} it suffices 
to show that the map
\[
(\psi, \delta \psi_2, \delta A_2) \to ( i A_2 \psi - \frac{1}r \delta A_2 \psi_2,  \Im{(\psi \bar{\psi}_2)}+ \frac{1}r |\delta A_2|^2)
\]
is locally Lipschitz from $L^2 \times \dHe \times \dHe$ to
$L^2 \times L^2$. This follows easily since $\dHe \subset L^\infty
\cap r L^2$. We note that the requisite smallness in the contraction
principle comes from the smallness of the right hand side in
\eqref{h1fowpsi3}, while the small Lipschitz constant is produced 
by unbalancing the norms
\[
\| (\psi, \delta \psi_2, \delta A_2)\|_{L^2 \times \dHe \times \dHe}
= \|\psi\|_{L^2} + \| \delta \psi_2\|_{\dHe} + M \|  \delta A_2\|_{\dHe}
\]
with large $M$ (and similarly for $L^2 \times L^2$).

b) After a time dependent rescaling and rotation we can assume 
that $\lambda=1$ and $\alpha = 0$ in \eqref{defdelta}. The price we pay is
twofold:

i)  As in part (a), the initial condition becomes 
$\delta \psi_2(\lambda)=0$, $\delta A_2(\lambda) = 0$. However,
we can shift back to $r=1$ using the time integrated form 
of \eqref{switchr}.

ii) The norm $LE[1]$ in \eqref{dap-eloc} is replaced by 
$LE[\lambda^{-1}]$. However, these two norms are equivalent.

After this reduction, it remains to prove the estimate
\begin{equation}
\|\frac{\langle r \rangle^{ -\epsilon}}{r} \delta \psi_2 \|_{L^2} +
\|\frac{ \langle r \rangle^{\frac12-\epsilon}}{r} \delta A_2\|_{L^2}\lesssim  
\| \psi \|_{LE[1]} +   \| \delta \psi_2(1)\|_{L^2} +  \| \delta A_2(1)\|_{L^2}
\end{equation}

In the interval $[0,R]$ this is obtained directly from \eqref{comp1rep2}
via Lemma~\ref{ldecay} with $\alpha = 0$. This yields
\[
\left\{
\begin{array}{l}
\| \delta \psi_2\|_{\dHe[0,1]} \lesssim | \delta \psi_2(1)|
+ \|\psi \|_{L^2[0,1]} + \| \delta A_2\|_{\dHe[0,1]}
(1+  \| \delta \psi_2\|_{\dHe[0,1]})
\cr \cr
\| \delta A_2\|_{\dHe[0,1]} \lesssim | \delta A_2(1)|
+ \|\psi \|_{L^2[0,1]} + \| \delta A_2\|_{\dHe[0,1]}^2 
\end{array} \right.
\]
From part (a) we have $\| \delta \psi_2\|_{\dHe}+ \| \delta A_2\|_{\dHe}
\lesssim \epsilon$, which allows us to close and obtain
\[
\| \delta \psi_2\|_{\dHe[0,1]}
+\| \delta A_2\|_{\dHe[0,1]}  \lesssim | \delta \psi_2(1)|
+ | \delta A_2(1)|
+ \|\psi \|_{L^2[0,1]}
\]
We square this and integrate in time.

It remains to consider the interval $[1,\infty)$. Here we apply
Lemma~\ref{ldecay} with $- \frac12 < \alpha < 0$ to obtain
\[
\begin{split}
\| r^{\alpha-1} \delta \psi_2\|_{L^2[1,\infty)} \lesssim & \ | \delta \psi_2(1)|
+ \| r^{\alpha-1}  T (A_2 \psi) \|_{L^2[1,\infty)} \\ & + \|r^{\alpha-1} \delta A_2\|_{L^2[1,\infty)}
(1+  \| \delta \psi_2\|_{\dHe})
\end{split}
\]
respectively
\[
\begin{split}
\|r^{\alpha-1} \delta A_2\|_{L^2[1,\infty)} \lesssim & \ | \delta A_2(1)|
+ \| r^{\alpha-1} T_1 (\psi_2 \psi) \|_{L^2[1,\infty) } \\ & + \|r^{\alpha-1} \delta A_2\|_{L^2[1,\infty)} 
\| \delta A_2\|_{\dHe} 
\end{split}
\]
At this point we use the assumption that $\psi$ is small in $L^\infty
LX$; by Proposition~\ref{constr} this implies $|\delta \psi_2|
\lesssim r^{-\frac12}$ and $|\delta A_2| \lesssim r^{-1}$.  These
bounds allow us to obtain a favorable bound for most of the
contributions of $\psi$, namely
\[
\| r^{\alpha-1} T ((A_2-1) \psi)\|_{L^2[1,\infty)}
+ \| r^{\alpha-1} T_1 (\psi_2 \psi) \|_{L^2[1,\infty) } 
\lesssim \| r^{\alpha -\frac12} \psi\|_{L^2}
\]
It remains to prove the linear estimate
\begin{equation}\label{tinl2}
\| r^{\alpha-1} T \psi \|_{ L^2[1,\infty)} \lesssim \|\psi\|_{LE[1]}
\end{equation}
For this we write
\[
 T \psi = L^{-1} \psi - h_1 L^{-1} \psi(1)
\]
We consider a dyadic decomposition of $\psi$, 
$\psi = \sum \psi_k$. 
For $k \geq 0$ the kernels $K_k^1$ of $L^{-1} P_k$,
described in Proposition~\ref{p-lp}, 
satisfy 
\[
|K_{k}^1(r,s)| \lesssim \frac{1}{(1+2^k|r-s|)^N(r+s)},
\]
This gives the $L^2$ bound 
\[
 \| r^{\alpha-1} L^{-1} \psi_k \|_{ L^2[1,\infty)} \lesssim 2^{-\frac{3k}2} \|\psi_k\|_{LE_k}
\]
and the pointwise bound
\[
 \| L^{-1} \psi_k (1)\|_{L^2_t} \lesssim 2^{-k} \|\psi_k\|_{LE_k}
\]
which are easy to sum up.

For $k < 0$ the kernels $K_k^1$ split into 
\[
K_k^1 = K_{k,reg}^1 + K_{k,res}^1
\]
where the regular part satisfies
\[
|K_{k,reg}^1(r,s)| \lesssim \frac{2^k \log{(1+r)}}{|k|(1+2^k|r-s|)^N(1+2^k(r+s))},
\]
while the resonant part, present only for $k < 0$, satisfies
\[
|K_{k,res}^1(r,s)| = \frac{1}{2^{k}|k|} h_1(r) \chi_{2^kr \leq 1}(r)  c_k(s)
\]
where $\chi_{2^kr \leq 1}$ is a bump which equals $1$ for $2^k r \ll 1$,
and $|c_k(s)| \lesssim (1+ 2^k s)^{-N}$.
For the regular part we have
\[
  \| r^{\alpha-1} L^{-1}_{reg} \psi_k \|_{ L^2[1,\infty)} + 
\|  L^{-1}_{reg} \psi_k(1) \|_{ L^2_t} \lesssim \frac{1}{|k| 2^{k}} \|\psi_k\|_{LE_k}
\]
which suffices due to the extra $|k| 2^{k}$ weight in the definition 
of $LE[1]$.

Finally for the resonant part we take advantage of the 
cancellation of the resonance. We have 
\[
 L^{-1}_{res} \psi_k = \frac{1}{2^{k}|k|} h_1(r) \chi_{2^k r \leq 1}(r)  f_k(t),
\qquad \|f\|_{L^2} \lesssim 2^{-k} \|\psi_k\|_{LE_k}
\]
Then for the corresponding part of $T$ we have
\[
 T_{res} \psi_k = \frac{1}{2^{k}|k|} h_1(r) (\chi_{2^kr \leq 1}(r)-1)  f_k(t)
\]
which leads to the stronger bound
\[
 \|r^{\alpha -1}  T_{res} \psi_k \|_{L^2} \lesssim 2^{-\alpha k}\frac{1}{2^{k}|k|} \|\psi_k\|_{LE_k}
\]
This concludes the argument for the part of \eqref{tinl2} concerning $\delta \psi_2$;
however we need to improve on the decay for the $\delta A_2$ term. We make the following
general observation which will be of use later too.
In some estimates we need better decay bounds for $\delta A_2$ near 
spatial infinity. For that we observe that for large $r$ the function 
$\delta A_2$ can be algebraically estimated as 
\begin{equation} \label{beterdelta}
 |\delta A_2| \lesssim h_1 |\delta \psi_2| + |\delta \psi_2|^2 
\end{equation}
For the first term on the right we have an extra order of decay.
For the second we can either use the $LX$ norm of $\psi$ to get another
half unit of decay, or we can get almost an $L^\infty L^1$ bound.  
In particular, this justifies the second part of \eqref{tinl2}. 
 \end{proof}

\section{The nonlinear equation for $\psi$}
\label{nonlin}

\subsection{A short time result}
We write the nonlinear equation for $\psi$ as
\begin{equation}
(i \partial_t - \tilde H) \psi = W \psi, 
\qquad \psi(0) = \psi_0, \qquad W = A_0 - 2 \frac{\delta A_2}{r^2}-
\frac{1}r \Im({\psi}_2 \bar{\psi})
\label{wnlin-eq1}\end{equation}
with $A_2$ and $\psi_2$   uniquely determined by $\psi$, see
Proposition \ref{constr}, and $\delta A_2=A_2 - h_3$. $A_0$ is
given by \eqref{aoef} which we recall for convenience
\[
A_0(r) =  -\frac12 |\psi|^2  + \frac{1}r \Im (\psi_2 \bar \psi)
+ [r\partial_r]^{-1}( |\psi|^2 - \frac{2}{r} \Im (\psi_2 \bar \psi))
\]
Treating the right hand side perturbatively, we prove a
local in time well-posedness result for \eqref{wnlin-eq1}:

\begin{t1} \label{t:localwp}
For each initial data $\psi_0$ satisfying 
\[
\|\psi_0\|_{LX} \leq \gamma \ll 1
\]
there is an unique solution $\psi$ for  \eqref{wnlin-eq1}  in the time interval $I=[0,1]$,
satisfying
\begin{equation} \label{psionI}
\| \psi\|_{\WSs[1](I)} \lesssim \gamma
\end{equation}
Furthermore, the solution map $ \psi_0 \to \psi$ is Lipschitz from
$LX$ to $\WSs[1](I)$.
\end{t1}
\begin{proof}
  By Proposition~\ref{p:mainS} it suffices to show that the map $\psi
  \to W\psi$ is Lipschitz from $\WSs[1](I)$ to $\WNs[1](I)$ with a small
  Lipschitz constant for $\psi$ as in \eqref{psionI}. We consider each term in $W$
and use Proposition~\ref{constr} to describe the dependence of $A_2$ and 
$\psi_2$ on $\psi$.  Some but not all of the estimates below depend 
on the size of the time interval $I$. To identify those we use the $\lesssim_I$ 
notation. For convenience, we use $\WSs, \WNs$ below instead of $\WSs[1](I), \WNs[1](I)$.

{\bf 1. The $A_2$ term}  is estimated using the bound
\begin{equation}
\| r^{-2}  f  g \|_{\WNs}  \lesssim_I \| f \|_{L^\infty X} \|g\|_{\WSs}
\end{equation}
For high frequencies in the output we use the local energy norms,
\begin{equation} \label{lwpc1}
\| r^{-2} f g\|_{LE^*}  \lesssim \| \la r \ra^\frac12  f\|_{L^\infty } \|g\|_{LE} 
\lesssim  \| f\|_{L^\infty X} \|g\|_{S}
\end{equation}
For low frequencies in the output we use \eqref{l1emb}  and an  $L^1$ bound
\[
\begin{split}
\| r^{-2} f g\|_{L^1}  \lesssim_I & \ \| r^{-2} f g\|_{L^2 L^1}  \lesssim 
\| |\log(1+r)|^{-1} f\|_{L^\infty L^2} \| r^{-2} \log(1+r) g\|_{L^2} \\
\lesssim & \   
\| f\|_{L^\infty X} \|g\|_{S}
\end{split}
\]

{\bf 2. The $\psi_2$ term} is estimated using the bounds
\begin{equation}
\begin{split}
\| r^{-1} f g h \|_{\WNs} +\| r^{-1} h_1  g h \|_{\WNs} \lesssim \ & (1+ \|f\|_{L^\infty X})
 \|g\|_{\WSs} \| h\|_{\WSs} 
\end{split}
\end{equation}
We only discuss the first bound; the second is similar but easier.
For high frequencies it suffices to write
\begin{equation}\label{lwpc2}
\| r^{-1} f g h \|_{LE^*} \lesssim   \| \la r \ra^\frac12  f\|_{L^\infty }  
\| g\|_{L^4} \|h\|_{L^4} \lesssim  \|f\|_{L^\infty X} \|g\|_{S} \|h\|_{S}
\end{equation}
For low frequencies we use \eqref{l1emb}  and an  $L^1$ bound derived 
using \eqref{betterle}:
\begin{equation} \label{lwpc5}
\begin{split}
\| r^{-1} f g h \|_{L^1} \lesssim & \
 \| \la r \ra^\frac12  f\|_{L^\infty }  
 \|r^{-\frac12} \la r \ra^{-\frac14}  g \|_{L^2} 
\|r^{-\frac12} \la r \ra^{-\frac14} h\|_{L^2}
\\
 \lesssim  & \ \|f\|_{L^\infty X} \|g\|_{\WSs} \|h\|_{\WSs}
\end{split}
\end{equation}

{\bf 3. The $|\psi|^2$ part of the $A_0$ term} 
 is estimated using the bounds
\begin{equation} 
\|  f g h \|_{\WNs} \lesssim_I  \|f\|_{\WSs} \|g\|_{\WSs} \| h\|_{\WSs}
\end{equation}
Indeed, for the high frequency part
we have 
\begin{equation} \label{lwpc3}
\|  f g h \|_{N} \lesssim \|  f g h \|_{L^\frac43} \lesssim   \|f\|_{L^4} \|g\|_{L^4} \| h\|_{L^4}
\lesssim  \|f\|_{S} \|g\|_{S} \| h\|_{S}
\end{equation}
while for the low frequency part we write
\begin{equation}
\|  f g h \|_{L^1} \lesssim_I \|  f g h \|_{L^2 L^1} \lesssim   \|f\|_{L^4} \|g\|_{L^4} \| h\|_{L^\infty L^2}
\lesssim  \|f\|_{S} \|g\|_{S} \| h\|_{S}
\end{equation}

{\bf  4. The $[r \partial_r]^{-1} |\psi|^2$ part of the $A_0$ term}  is estimated 
as in Case 3 using in addition the Hardy type inequality 
\eqref{rdrm} for $[r \partial_r]^{-1}$.

{\bf 5. The $[r \partial_r]^{-1} ( r^{-1} \Im (\psi_2 \bar \psi))$
  part of the $A_0$ term} is estimated as in Case 2: using \eqref{rdrm}
  with $p=2$ for the corresponding case to \eqref{lwpc2} and using \eqref{rdrmw}
  with $p=1$ and $w=\la r \ra^{-\frac12}$ for the corresponding case to 
  \eqref{lwpc5}.
\end{proof}

\subsection{The long time result} 
We rewrite the equation for $\psi$ in the form
\begin{equation}
(i \partial_t - \tilde H_\lambda) \psi = W_\l \psi, 
\quad \psi(0) = \psi_0, \quad  
W_\l = A_0 - 2 \frac{\delta^\l A_2}{r^2}-\frac{1}r \Im({\psi}_2 \bar{\psi})
\label{wnlin-eq}\end{equation}
with $A_0$, $A_2$ and $\psi_2$ as well as the time dependent parameter
$\lambda$ uniquely determined by $\psi$ (see \eqref{defal} and
\eqref{defdelta}).  Our main long time bootstrap result is as follows:

\begin{t1} \label{thflow}
 Let  $T \in (0,\infty]$, $\epsilon \leq 1$ and $\gamma \leq \gamma_0 \ll 1$.
Suppose that the initial data for $\psi$ satisfies
\begin{equation}
\| \psi(0)\|_{L^2} \leq \epsilon \gamma, \qquad
\| \psi(0)\|_{LX} \leq \gamma
\label{psidata}\end{equation}
Assume that the parameter 
$\lambda$ and the function $\psi$ satisfy the following bootstrap assumptions:
\begin{equation}
\qquad \|\lambda-1\|_{Z_0[0,T]} \leq \gamma_0.
\label{bootlambda}\end{equation}
respectively
\begin{equation}
\| \psi\|_{\ldSs[0,T]} \leq \epsilon \gamma_0, \qquad
\| \psi\|_{\WSs[\tl][0,T]} \leq \gamma_0, 
\label{bootpsi}\end{equation}
where $\tl$ is any function with the following properties:
\begin{equation}
\qquad \|\tl-1\|_{\Zu[0,T]}+ \| \l - \tl\|_{(L^2 \cap L^\infty)[0,T]} \lesssim \gamma_0
\label{tlambda}\end{equation}

Then the functions $\psi$ and $\lambda$ must satisfy 
the stronger bounds
\begin{equation}
\| \psi\|_{\ldSs[0,T]} \lesssim \epsilon (\gamma + \gamma_0^2), \qquad
\| \psi\|_{\WSs[\tl][0,T]} \lesssim \gamma + \gamma_0^2, 
\label{rebootpsi}\end{equation}
respectively
\begin{equation}
\qquad \|\lambda-1\|_{Z_0[0,T]} \lesssim \gamma+\gamma_0^2.
\label{rebootlambda}\end{equation}
\end{t1}

To close the bootstrap it suffices to choose $\gamma_0 = C \gamma$ for
a fixed large universal constant $C$. 

We remark that for the global well-posedness result it suffices to
take $\epsilon = 1$. However, the parameter $\epsilon$, along with the
stronger bounds in \eqref{betterbootx}, is needed for the proof of the
instability result.

The additional parameter $\tl$ is needed because
the spaces $\WSs[\lambda]$, $\WNs[\lambda]$ and the linear result
in Proposition~\ref{p:lin-x} require $\lambda -1 \in \Zu$, while 
above we only have $\lambda -1 \in Z_0$.
There is some flexibility in the choice of $\tl$. An acceptable choice 
would be for instance $\tl = Q_{\leq 1} \lambda^{ext}$ where 
$\lambda^{ext}$ is any suitable extension of $\lambda$ in $Z_0$.

For brevity we omit the time interval $[0,T]$ in the notations 
in this section. For the most part this plays no role.  At one point 
in the proof this  requires an additional discussion.

For the first bound in \eqref{rebootpsi} we  use Theorem~\ref{p:lin-l2}.
Hence it suffices to estimate the nonlinear expression $W\psi$, for which 
we will prove 
\begin{equation}
 \| W \psi\|_{\ldNs} \lesssim \epsilon \gamma_0^2.
\label{betterboot}\end{equation}

For the second bound \eqref{rebootpsi} for $\psi$ we 
rewrite the equation in the form 
\begin{equation}
(i \partial_t - \tilde H_{\tl}) \psi = (V_\l-V_{\tl})\psi + W_\l \psi, 
\quad \psi(0) = \psi_0
\label{wnlin-eq-tl}\end{equation}
and use Theorem~\ref{p:lin-x}.    Hence it suffices 
to estimate the linear term 
$(V_\l-V_{\tl})\psi$ and the nonlinear expression $W\psi$. Precisely, 
we will prove the bounds
\begin{equation}
\|(V_\l-V_{\tl})\psi\|_{\WNs[\tl]} \lesssim \epsilon \gamma_0^2
\label{betterbootV}\end{equation}
\begin{equation}
 \| W_\l \psi\|_{\WNs[\tl]} \lesssim  \epsilon \la \log \e \ra^2 \gamma_0^2.
\label{betterbootx}\end{equation}
The three bounds above follow from Propositions~\ref{p:nonlin-l2}, \ref{p:nonlin-V},
\ref{p:nonlin-x} below. The last part of the theorem, namely 
the $\lambda$ bound \eqref{rebootlambda}, is proved in the next section.
We begin with the nonlinear bound in $L^2$:
\begin{p1}\label{p:nonlin-l2}
Suppose that $\psi$  satisfies \eqref{bootpsi}.
Then
 \begin{equation}
\| W u\|_{N} \lesssim  \gamma_0 \|u\|_{S} 
\label{Wl2}\end{equation}
\end{p1}

\begin{proof}
  We consider each of the terms in $W$ as in the five cases in the
  proof of Theorem~\ref{t:localwp}.  The bound \eqref{Wl2} follows by
  applying \eqref{lwpc1}, \eqref{lwpc2}, \eqref{lwpc3} and their
  counterparts in the last two cases. It is essential that none of
  these bounds depend on the size of the time interval.

\end{proof}

Next we consider the linear bound \eqref{betterbootV}:

\begin{p1} \label{p:nonlin-V}
Suppose that  $\lambda$ and $\tl$  satisfy \eqref{bootlambda} and 
\eqref{tlambda}. Then
\begin{equation} \label{Vest}
\|(V_\l- V_{\tl}) u\|_{\WNs[\tl]} \lesssim  \gamma_0 \|u\|_{l^2\Ss}  
\end{equation}
\end{p1}

\begin{proof}
Since both $\l$ and $\tl$ are close to $1$, it follows that
\[
 |V_{\lambda_1} - V_{\lambda_2}| \lesssim |\lambda_1-\lambda_2|
(1+r^2)^{-2}
\]
Hence, using the embedding \eqref{LXemb} and the $LE$ norm for $u$, we 
obtain
\[
\begin{split}
 \|( V_{\lambda_1} - V_{\lambda_2})u\|_{L^1 LX} \lesssim &
\|( V_{\lambda_1} - V_{\lambda_2})u\|_{L^1 (L^1 \cap L^2)}\lesssim 
\| \la r \ra^{3-} (V_{\lambda_1} - V_{\lambda_2}) u \|_{L^1 L^2}
\\
\lesssim & \| u \|_{LE} \|\lambda_1-\lambda_2\|_{ L^2}
\end{split}
\]
Thus \eqref{Vest} follows.
\end{proof}

Next we consider the bound \eqref{betterbootx},
which corresponds to initial data in the smaller space $LX \subset L^2$.

\begin{p1} \label{p:nonlin-x}
Suppose that $\psi$  satisfies \eqref{bootpsi} and $\lambda$ satisfies
\eqref{bootlambda}. Then for $\tl$ as in \eqref{tlambda} we have
\begin{equation} \label{nonest}
\| W_\l u\|_{\WNs[\tl]} \lesssim \epsilon \la \log \e \ra^2 \gamma_0 \|u\|_{\WSs[\tl]} + 
\g \| u \|_{\ldSs}
\end{equation}
\end{p1}

\begin{proof} 
  The difference in the potentials $W_\l - W = 2
  \dfrac{h_3^\l-h_3}{r^2}$ decays rapidly enough at $\infty$ so that,
  in a similar manner to \eqref{Wl2}, one easily derives
\[
\| (W_\l - W) u \|_{N} \lesssim |\l-1| \| u \|_{S} \lesssim \g \| u \|_{S}
\]
Combining this with \eqref{Wl2} gives us
\begin{equation} \label{weasy}
\| W_\l  u \|_{N} \lesssim \g \|u\|_{S}
\end{equation} 
By the inclusion $N \subset l^2 N $ and the first part of
\eqref{wsl=ws}, we can use the above results to  estimate 
the high frequency output $\|P^{\tl}_{\geq 0}(W_\lambda
u)\|_{\WNs[\tl]}$.

It remains to estimate the low frequency output 
$\|P^{\tl}_{< 0}(W_\lambda u)\|_{\WNs[\tl]}$. We divide the potential
$W_\lambda$ into three parts, $W_\lambda = W_0 + W_1+W_2$
where
\[
W_0 = \frac{\delta^\l A_2}{r^2}, \qquad 
 W_1 = - \int \frac{1}{r^2} \Im (\psi_2 \bar \psi), 
\qquad
W_2 = -\frac12 |\psi|^2 - \int \frac{1}{r}|\psi|^2 dr  
\]
The contribution of $W_0$ can be estimated directly
by using the local energy decay for 
$u$ and $\delta^\l A_2$, see \eqref{le=lek} and \eqref{dap-eloc}, 
 to write
\[
\|W_0 u\|_{L^1 LX} \lesssim 
\| W_0  u\|_{L^1} 
\lesssim 
\left\| \frac{\log(2+r)}{r} \delta^\l A_2 \right\|_{L^2} 
\left\| \frac{1}{r \log(2+r)} u \right\|_{L^2} 
\lesssim  \gamma_0 \|u \|_{S}
\]
A similar argument applies for the contribution of $W_1$. Indeed,
consider first the simpler potential $\tilde W_1 = \dfrac{\psi_2}{r}
\bar\psi$ and use the pointwise bound for $\psi_2$, $|\psi_2| \lesssim
\la r \ra^{-\frac12}$.  Using \eqref{betterle} and \eqref{bootpsi} we obtain
\[
\begin{split}
\| \tilde W_1 u\|_{L^1} & \lesssim \| \la r \ra^\frac12 \psi_2 \|_{L^\infty} \left\|\frac{\ln(2+r)}r  \psi\right\|_{L^2} 
\left\| \frac{1}{\la r \ra^\frac12 \ln(2+r)}u\right\|_{L^2} \\
& \lesssim \epsilon \la \log \epsilon \ra^2 \g \|u\|_{\WSs[\tl]}
\end{split}
\]
A similar bound is obtained for $W_1$ by using \eqref{rdrmw} with $p=2$ and $w=\la r \ra^\frac12 \ln(2+r)$.

The bulk of the proof is devoted to the low frequency estimate for
$W_2 u$, which is independent of $\l$. It suffices to show that
\[
\| P_{<0}^{\tl}  (W_{2}\ u)\|_{\WNs[\tl]}
\lesssim \epsilon \gamma_0^2 \| u\|_{\WSs[\tl]}+ \gamma_0^2 \|u\|_{\ldSs}
\] 
The expression $W_2 u$  is a trilinear expression in  $(\psi,\psi,u)$,
\[
W_2 u = N(\psi,\psi,u) + \tilde N(\psi,\psi,u),
\]
where
\[
N(\psi_1,\psi_2,\psi_3)=-\frac12 \psi_1 \overline{\psi}_2 \psi_3  \qquad
\tilde{N}(\psi_1,\psi_2,\psi_3)=-\psi_3 [r \partial_r]^{-1} ( \psi_1 \overline{\psi}_2)
\]
Given the bounds \eqref{bootpsi}
on $\psi$, the above inequality can be rewritten in a more symmetric
way as
\begin{equation}
\begin{split}\!\!
\| P_{<0}^{\tl} N(\psi^1\!,\psi^2\!,\psi^3) \|_{\WNs[{\tl}]} 
% + \| P_{<0}^{\tl} \tilde N(\psi^1,\psi^2,\psi^3) \|_{\WNs[{\tl}]}\\
\!\! \lesssim \!\! \sum_{\sigma \in S_3} \!\!
\| \psi^{\sigma(1)} \|_{\ldSs} \|\psi^{\sigma(2)}\|_{\WSs[{\tl}]}  
\|\psi^{\sigma(3)}\|_{\WSs[{\tl}]} 
\end{split}
\end{equation}
and the similar bound for $\tilde N$.  To avoid repetition, we focus
on establishing this inequality for $N$. Every step in the analysis
for $N$ has its counterpart for $\tilde N$.  As a general rule, the
estimate for $\tilde N$ is similar to the one for $N$ by the use of
\eqref{rdrm} and \eqref{rdrmw}, with one exception which require 
separate analysis.

We decompose in the time dependent frame
$\psi^i  = \sum \psi^i_k$, i.e. $\psi^i_k=P^{\tl}_k \psi^i$
For a large constant $n_0$ we expand
\[
\begin{split}
& P_{<0}^{\tl} N(\psi^1,\psi^2,\psi^3)
= \sum_{\sigma \in S_3}
\sum_{k_1 \geq k_2 \geq k_3} P^{\tl}_{<{k_3^--2n_0} }
 N^{\sigma}(\psi^{\sigma(1)}_{k_1} ,\psi^{\sigma(2)}_{k_2},\psi^{\sigma(3)}_{k_3}) \\
 & + \sum_{j < 0}  P_{[j-2n_0,0]}^{\tl} \left( N(\psi^1_{j} ,\psi^{2},\psi^{3}) + N(\psi^1_{\geq j} ,\psi^{2}_{j},\psi^{3})+ N(\psi^1_{\geq j} ,\psi^{2}_{\geq j},\psi^{3}_{j})\right)
\end{split}
\]
where $N^{\sigma}(f^1, f^2, f^3):=N(f^{\sigma^{-1}(1)},
f^{\sigma^{-1}(2)}, f^{\sigma^{-1}(3)})$.  We observe that the second
term has a favorable frequency balance and is estimated directly using
only the Strichartz norms,
\[
\begin{split}
\|  P_{[j-2n_0,0]}^{\tl} N(\psi^1_{j} ,\psi^{2},\psi^{3}) \|_{\WNs[{\tl}]}
\lesssim & \ \frac{1}{| j | 2^j} \| N(\psi^1_{j} ,\psi^{2},\psi^{3})\|_{N}\\
\lesssim & \ \sum_{k < j} \frac{1}{| j | 2^j} \| \psi^{1}_{j}\|_{S}  \|\psi^2\|_S \|\psi^3\|_{S}
\end{split}
\]
followed by a straightforward summation with respect to $j$. The other
two nonlinear factors in the second term are treated similarly. Using
\eqref{rdrm} the same argument works for $\tilde N$.

The first term requires considerably more work; in what follows, 
we will prove that for all $k_1 \geq k_2 \geq k_3$ and $\sigma \in S^3$
the following bound holds:
\begin{equation}
  \| P_{< k_3^--2n_0} ^{\tl}
N^\sigma(\psi_{k_1}^1,
\psi_{k_2}^{2},\psi_{k_3}^{3})\|_{\ldNs} \lesssim 
\frac{\|\psi_{k_1}^{1} \|_{\ldSs} \|\psi_{k_2}^{2} \|_{\ldSs} \|\psi_{k_3}^{3} \|_{\ldSs}}
{ 2^{{k_1^+}/8} \la k_2^- \ra  2^{k_2^-} \la k_3^-\ra  2^{k_3^-}} 
\label{w2}\end{equation}
as well as  the similar one for $\tilde N$.

One difficulty is that the functions $\psi^{\sigma(j)}_{k_j}$ are 
only localized in the frequency dependent frame. To deal with this we 
relocalize them in the fixed frame,
\[
 \psi_{k_j}^{\sigma(j)} = \tilde P_{k_j} \psi_{k_j}^{\sigma(j)}  + \psi_{k_j}^{\sigma(j),err}
\]
and estimate pointwise the error $\psi_{k_j}^{j,err}$ using
Corollary~\ref{reloc}, see \eqref{relocest},
\[
| \psi_{k_j}^{j,err}(r)| \lesssim \frac{2^{-k_j^+/2}}{r \log^2(2+r)} 
\|\psi_{k_j}^{j} \|_{L^\infty L^2}
\]
Then for the part of \eqref{w2} containing at least one
error term we combine this with the local energy decay estimates and \eqref{point-le},
\eqref{le=lek} and the low frequency part of
\eqref{LXemb}. One such term is 
\[
\begin{split}
\| N^\sigma(\psi_{k_1}^{1}, \psi_{k_2}^{2}, \psi_{k_3}^{3}) \|_{L^1} 
 \lesssim & \ \|r \log^2(2+r) \psi_{k_1}^{1,err}\|_{L^\infty} \prod_{j=2,3}
\left\| \frac{\psi_{k_j}^{j}}{r^\frac12 \log(2+r)} \right\|_{L^2} 
 \\
\lesssim & \ \frac{\|\psi_{k_1}^{1}\|_{\ldSs} \|\psi_{k_2}^{2} \|_{\ldSs} 
\|\psi_{k_3}^{3} \|_{\ldSs}}{2^{k_1^+/2} 2^{\frac{k_2}2}  2^{\frac{k_3}{2}}} 
\end{split}
\]
The other terms are treated in a similar way.  In the case of $\tilde N$,
the same argument works when combined with \eqref{rdrmw} as follows:
$p=2$ and $w=r^\frac12 \log(2+r)$ if the $err$ term is inside $[r\partial_r]^{-1}$
or $p=1$ and $w=\frac{1}{r \log^2(2+r)}$ otherwise.

After the above reduction we can assume that the functions
$\psi_{k_j}^{j}$ in \eqref{w2} are localized in the fixed frame. 
Thus the inputs $\psi_{k_j}^{j}$ no longer bear any relation 
to $\tl$.

So far the time interval has played no role. At this point
we consider suitable frequency localized\footnote{with respect to the fixed
$\lambda = 1$ frame} extensions for both 
$\psi_{k_j}^{j}$ and $\tl$, and prove that \eqref{w2} holds 
over the entire real line. This directly implies the similar 
bound over each subinterval.

A second difficulty is that the output $N$ is also localized in the 
time dependent frame, while it is more natural at this point to 
 project  the output in a fixed frame. We  denote 
\[
 g = N^\sigma(\psi_{k_1}^{1}, \psi_{k_2}^{2},\psi_{k_3}^{3})
\]
and normalize 
\begin{equation} \label{renorm}
 \| \psi_{k_i}^{\sigma(i)}\|_{\ldSs} = 1,
\end{equation}
We partition $g$ in frequency/modulation and estimate each 
piece as in \eqref{w2}.

\medskip 
{\bf Case 1: The high frequency part $P_{\geq k_3^--n_0} g$.}
Using \eqref{point-le}, we obtain
\begin{equation}\label{ginl1}
\|(1+2^{k_3} r)^{-\frac12} g\|_{L^1} \lesssim 2^{-\frac{k_1+k_2+6k_3}8}
 \|\psi_{k_1}^{\sigma(1)} \|_{L^4_{k_1}}
 \|\psi_{k_2}^{\sigma(2)}\|_{L^4_{k_2}}\|\psi_{k_3}^{\sigma(3)}\|_{LE_{k_3}}
\end{equation}
which is satisfactory for $r \lesssim 2^{-k_3}$ but not for larger
$r$.  We use this bound to estimate $P_{j}^{\tl} P_{\geq k_3^--n_0}$ in
$L^1 L^2$ for $j < k_3^--2n_0$.  From Proposition \ref{p:bumpft} and
Proposition \ref{p:transfa} it follows that for $j < k_3^--2n_0$ the kernel
of $P_j^{\tl} P_{\geq k_3^--n_0}$ satisfies
\[
\begin{split}
 |K(r,s)| \lesssim & \frac{2^{2j}}{j^2} \frac{\log (1+r)}{ (1+2^j r)^{N}} 
\left(
\sum_{k=k_3^--n_0 }^0   \frac{1}{\la k \ra^2}\frac{\log (1+s)} {(1+2^k s)^{N}} 
+  \sum_{k =1 }^\infty    \frac{ 2^{-Nk}s^2}{  (1+ s)^{N}}\right)  \\
\lesssim &\ \frac{2^{2j}}{j^2} \frac{\log (1+r)}{(1+2^j r)^{N}} 
\frac{1} {(1+2^{k_3^-} s)^{N}}
\end{split}
\]
Thus by \eqref{ginl1} we obtain
\[
\begin{split}
 \| P_j^{\tl} P_{\geq k_3^--n_0} g\|_{L^1L^2} & \lesssim \frac{2^j}{|j|}
2^{-\frac{k_1+k_2+6k_3}8} 2^{\frac{k_3^+}2}
\end{split}
\]
Summing in $j$ with weights $(2^j |j|)^{-1}$ we estimate 
$  P_{\geq {k_3^--n_0}} g$ as in \eqref{w2}. The argument 
for $\tilde N$ is obtained by using \eqref{rdrmw} with $p=2$
and $w=(1+2^{k_3}r)^{\frac18}$
if the $k_3$ frequency term is outside $[r\partial_r]^{-1}$
and with $p=\frac43$ and $w=(1+2^{k_3} r)^{-\frac12}$
if the $k_3$ frequency term is inside $[r\partial_r]^{-1}$.

\medskip

{\bf Case 2: Low frequency, high modulations:  $Q_{\geq k_1+k_3 - n_0} 
P_{< k_3^--n_0} g$.}
For this term we prove an $L^2$ bound
\begin{equation}
\sum_{j < k_3^--n_0} \frac{1}{2^j |j|}\| Q_{\geq k_1+k_3-n_0} P_j g \|_{L^2} 
\lesssim \la k_2^-\ra
\label{g3}\end{equation}
which leads to an estimate as in \eqref{w2}  by using \eqref{wnemb} 
(precisely, its high modulation part from \eqref{inNl}).  We have
\[
\| (1+ r 2^{k_2})^{\frac18}  g \|_{L^2 L^1} \lesssim 
\| \psi_{k_1}^{1} \|_{L^\infty L^2} \| \psi_{k_2}^{2} \|_{L^4_{k_2}}
 \| \psi_{k_3}^{3} \|_{L^4} 
\lesssim 1.
\] 
Hence, arguing as in the proof of \eqref{l1emb}, it follows that for $j < 0$ 
we have
\[
\frac{1}{2^j |j|} \|P_j g\|_{L^2} \lesssim \sum_{m} \frac{2^{m^-}\la m^+ \ra}{j^2} 
\|  g \|_{L^2 L^1(A_m)} \lesssim  \frac{\la k_2^- \ra}{j^2} 
\]
and \eqref{g3} follows after a $j$ summation. The bound
for $\tilde N$ is obtained by using \eqref{rdrm} with $p=\frac43$ if the $k_2$
frequency term is outside $[r\partial_r]^{-1}$ and \eqref{rdrmw} when the $k_2$
frequency term is inside $[r\partial_r]^{-1}$ as follows:
with $w=(1+ r 2^{k_2})^{\frac18}$, $p=2$ when the $k_3$ frequency term is also 
inside $[r\partial_r]^{-1}$
and with $w=(1+ r 2^{k_2})^{\frac18}$, $p=\frac43$ when 
the $k_3$ frequency term is outside $[r\partial_r]^{-1}$.
\medskip

{\bf Case 3: The low frequency, low modulations  $Q_{< k_1+k_3-n_0} 
P_{< k_3^--n_0} g$.} 
We will prove that 
\begin{equation}
\|Q_{< k_1+k_3-n_0} 
P_{< k_3^-} g \|_{L^1 LX} \lesssim 
\frac{1}{ 2^{{k_1^+}/8} \la k_2^- \ra  2^{k_2^-} \la k_3^-\ra  2^{k_3^-}}
\label{g2}\end{equation}
which, by \eqref{wnemb}, suffices for \eqref{w2}. We separate this into two cases:

{\bf Case 3A: One input has high modulation $\geq  2^{k_1+k_3-n_0}$}.
Say that factor is $\psi^1_{k_1}$; the other cases are similar. Then we bound
\[
\begin{split}
& \| N^\sigma(Q_{\geq k_1+k_3-n_0} \psi_{k_1}^{1},\psi_{k_2}^{2},\psi_{k_3}^{3})\|_{L^1} \\
\lesssim & \| Q_{\geq k_1+k_3-n_0} \psi_{k_1}^{1}\|_{L^2} 
\| \psi_{k_2}^{2}\|_{L^4}\|\psi_{k_3}^{3}\|_{L^4} \lesssim 2^{-\frac{k_3+k_1}2}
\end{split}
\]
and conclude with \eqref{l1emb}.  We note that the $L^4$ bound is
stable under cut-offs in modulation ($\leq 2^{k_1+k_3-n_0}$) due to the
$V^2_{\tilde{H}}L^2$ structure, see \eqref{modtr}. This  works 
for $\tilde N$ as well, by using \eqref{rdrm}.

{\bf Case 3B: All factors $\psi_{k_i}^i$ have low modulation.}
By duality, it suffices to estimate the quadrilinear integral
\[
 I_0 = \int N^\sigma( Q_{<k_1+k_3-n_0} \psi_{k_1}^{1},  Q_{<k_1+k_3-n_0} \psi_{k_2}^{2},
Q_{<k_1+k_3-n_0} \psi_{k_3}^{3}) \overline{Q_{<k_1+k_3-n_0} \psi_j}\ rdr dt
\]
with frequency localized inputs and show that
\begin{equation} \label{io}
|I_0| \lesssim \frac{2^j }{|j|}
\frac{1}{ 2^{{k_1^+}/8} \la k_2^- \ra  2^{k_2^-} \la k_3^-\ra  2^{k_3^-}}
\|\psi_{j}\|_{L^\infty L^2}
\end{equation}

We begin with a short frequency/modulation analysis. 
If the frequencies in the four factors are $\xi_1$, $\xi_2$, $\xi_3$ and $\xi$ 
and all modulations are $\ll 2^{k_1+k_3}$ then the time frequency
in the $I$ integrand is  
\[
 \phi = \xi_1^2 -\xi_2^2 + \xi_3^2 - \xi^2 + m
\]
where $m$ is the sum of the four modulations involved, hence $m \ll 2^{k_1+k_3}$.
Hence the time integral vanishes unless
\begin{equation}\label{defd}
(\xi_1,\xi_2,\xi_3,\xi) \in D = 
\{\xi_i \approx 2^{k_i}, \xi \approx 2^j; |\xi_1^2 -\xi_2^2 + \xi_3^2 - \xi^2| \ll 2^{k_1+k_3}\}
\end{equation}
Given the dyadic localization of $\xi_1$, $\xi_2$, $\xi_3$ and $\xi$,
this leads to one of the following two scenarios (recall that $j+n_0 \leq k_3 \leq k_2 \leq k_1$):
\begin{itemize}
\item[(i)] Equal frequency inputs,
\[
|k_1-k_2| \lesssim 1,  \quad |k_2 - k_3| \lesssim 1, \qquad  
|\xi_1^2+  \xi_2^2 - \xi_3^2 | \ll 2^{2k_3}.
\]
\item[(ii)] Unbalanced frequency inputs,
\[
|k_1 - k_2| \lesssim 1, \quad k_3 \ll k_2, \qquad  |\xi_1 -  \xi_2 | \ll  2^{k_3}.
\]
\end{itemize}
On the other hand, resonant interactions can only occur when
\begin{equation}
 \xi_1 \pm \xi_2 \pm \xi_3 \pm \xi = 0
\label{null_freq}\end{equation}
where the $\pm$ signs correspond to outgoing/incoming waves.  But this
is precluded in both cases (i) (ii). We will strongly exploit this
fact in our analysis.

In proving \eqref{io} we only use the $S_{k_3}$ norm for
$\psi_{k_3}^{3}$, together with boundedness of $Q_{<k_1+k_3-n_0}$ in
$S_{k_3}$.  For $ Q_{<k_1+k_3-n_0} \psi_{k_1}^{1}$ and
$Q_{<k_1+k_3-n_0} \psi_{k_2}^{2}$ we would also like to use the
norms $S_{k_1}$ respectively $S_{k_2}$.  Unfortunately, the operator
$Q_{<k_1+k_3-n_0}$ is not always bounded in these spaces; 
instead, from \eqref{modcut}, we have
\begin{equation} \label{lemod}
\| Q_{\leq k_1+k_3-n_0} \psi_{k_i}^{i} \|_{S_{k_i}} \lesssim 
\left(1+k_1-k_3\right)  \| \psi_{k_i}^{i} \|_{\Ss_{k_i}}, \qquad i= 1,2
\end{equation}
In our case this leads to losses 
of at most a factor of $1+(k_3 - k_1)^2$. Fortunately we are able to
 prove a stronger bound
\begin{equation} \label{ioa}
|I_0| \lesssim \frac{2^j }{|j|}
\frac{ \la k_3^-\ra}{2^{\frac{k_1+k_3}2}}
\|\psi_{k_1}^{1}\|_{S_{k_1}} \|\psi_{k_2}^{2}\|_{S_{k_2}} \|\psi_{k_3}^{3}\|_{S_{k_3}}
\|\psi_{j}\|_{L^\infty L^2}
\end{equation}
which can absorb these losses and still lead to \eqref{io}.
To keep the size of formulae below manageable, we normalize 
all four norms on the right to $1$ in the sequel.

Without restricting the generality of the argument, we restrict our attention to
$\sigma(i)=i$, in which case 
\[
I_0 = \int_0^\infty \psi^{1}_{k_1}(r) \overline{\psi}^{2}_{k_2}(r) 
\psi^{3}_{k_3}(r) \overline{\psi}_{j}(r) rdr dt
\]
In order to treat $\tilde N$, we need to consider also (the others are similar)
\[
I_1 = \int_0^\infty \overline{\psi}_j(r) \psi_{k_3}^{3}(r) \int_r^\infty \frac1s
\psi_{k_1}^{1}(s) \overline{\psi}_{k_2}^{2}(s) ds rdr dt
\]
\[
 I_2 = \int_0^\infty \overline{\psi}_{j}(r) \psi_{k_1}^{1}(r) \int_r^\infty \frac1s
\psi_{k_2}^{2}(s) \overline{\psi}_{k_3}^{3}(s) ds rdr dt
\]
Here all factors have the appropriate frequency and modulation
localization.

 The analysis is similar for all these quadrilinear forms, so we
will work with $I_0$. We switch $I_0$ to the Fourier space, where it
becomes
\[
\begin{split}
I_0\!= &\int \! \psi_{\xi_1}(r) \psi_{\xi_2}(r) \psi_{\xi_3}(r) \psi_\xi(r)r dr
\hat{\psi}_{k_{1}}^{1}(t,\xi_1) \overline{\hat{\psi}}_{k_{2}}^{2}(t,\xi_2) 
\hat{\psi}_{k_{3}}^{3} (t,\xi_3) \overline{\hat{\psi}}_{j}(t,\xi) d\xi_i d \xi dt
\\ =  & \int (G_0\chi_D)(\xi,\xi_1,\xi_2,\xi_3)
\hat{\psi}_{k_{1}}^1(t,\xi_1)
\overline{\hat{\psi}}_{k_{2}}^2(t,\xi_2) 
\hat{\psi}_{k_{3}}^3(t,\xi_3) \overline{\hat{\psi}}_{j}(t,\xi) d\xi_i d \xi dt
\end{split}
\]
where $G_0$ is the quadrilinear form on $\tilde H$-generalized 
eigenfunctions,  introduced in Section~\ref{sec:g}, and
$\chi_D$ is any function which equals $1$ in the set $D$ 
defined in \eqref{defd}. Similarly we can write $I_1$ and $I_2$ 
in terms of $G_1$ and $G_2$.

In practice we always restrict the support 
of $\chi_D$ to the sets described in cases (i), (ii) above,
and we assume it has good regularity. Hence in the support 
of $\chi_D$ we can use the bounds in Proposition~\ref{p:g}
for $G_0$, $G_1$ and $G_2$.  Thus, from \eqref{cfour}, $G_0$ satisfies
\[
|G_0| \lesssim g_{jk_1k_2k_3}=\frac{2^{\frac{j}2}}{|j|} \frac{\la k_3^-\ra 2^{-2k_3^+}}{2^{\frac{k_3}2}}
\]
and is smooth in $(\xi_1,\xi_2,\xi_3)$ on the $2^{k_3}$ scale,
and in $\xi$ on the $2^j$ scale. We consider the two cases
(i) and (ii) described above.

{\bf Case 3B(i)}. Here $|k_1 - k_2|, |k_2 - k_3| \lesssim 1$ therefore 
$G_0\chi_D$ has symbol regularity in all variables. Hence separating
the variables it suffices to look at $G_0$ of the form
\[
G_0 = g_{j k_1 k_2 k_3}
\chi_{k_1}(\xi_1)\chi_{k_2}(\xi_2) \chi_{k_3}(\xi_3)
\chi_{j}(\xi)
\]
where the $\chi_k$'s are smooth normalized dyadic bump functions.
Then the quadrilinear integral becomes
\[
I_0 =  g_{j k_1 k_2 k_3} \int \la \hat \chi_{k_{1}},\psi_{k_{1}}^1 \ra 
\la \hat \chi_{k_{2}}, \overline{\psi}_{k_{2}}^2 \ra \la \hat \chi_{k_{3}},\psi_{k_{3}}^3 \ra
\la \hat \chi_{j}, \overline{\psi}_{j} \ra dt
\]
Hence using the local energy norm for the first two $\psi$'s, 
the energy for the last two and \eqref{bumpft} for the inverse 
FT of bump functions we obtain
\[
|I_0| \lesssim g_{j k_1 k_2 k_3} 2^{\frac{k_3+j-k_1-k_2}2} \lesssim
 \frac{2^{j}}{|j|} \frac{\la k_3^-\ra}{2^{k_1}}
\]
which easily implies \eqref{ioa}. Since in this case $G_1$ and $G_2$
satisfy similar bounds, we also obtain \eqref{ioa} for $I_1$ and $I_2$.

{\bf Case 3B(ii)}. Here $|k_1-k_2| \lesssim 1, |k_3-k_2| \gg 1$, and the function $G_0\chi_D$
has symbol regularity in $\xi_3$ and $\xi$, but only on the $2^{k_3}$
scale in $\xi_1$ and $\xi_2$. Furthermore we can use $\chi_D$ to
localize it in the region $|\xi_1 -\xi_2| \ll 2^{k_3}$.  Thus we can
separate variables and obtain a decomposition of $G_0$ of the form
\[
G_0(\xi_1,\xi_2,\xi_3,\xi) = g_{j k_1 k_2 k_3}
\sum_{l=1}^{2^{k_1-k_3}} \chi_{k_1}^l(\xi_1)\chi_{k_2}^l(\xi_2)
\chi_{k_3}(\xi_3)
\chi_{j}(\xi)
\]
where $\chi_{k_1}^l,\chi_{k_2}^l$ have similar $2^{k_3}$ sized supports.
Using an argument similar to the one in Proposition \ref{p:bumpft} we obtain
\[
 |\hat \chi_{k_{1,2}}^l(r)| \lesssim 2^{k_3+\frac{k_1}2} (1+2^{k_1 }r )^{-\frac12} 
(1+2^{k_3} r)^{-N}
\]
At this point a  direct computation as in case (i) 
for the corresponding part $I_0^l$ of $I_0$ gives 
\[
|I_0^l| \lesssim \frac{2^{j}}{|j|} \frac{\la k_3^-\ra}{2^{k_1}}
\]
and, after summation with respect to $l$,
\begin{equation}\label{nogood}
|I_0| \lesssim \frac{2^{j}}{|j|} \frac{\la k_3^-\ra}{2^{k_3}}.
\end{equation}
Unfortunately this bound is not strong enough for \eqref{ioa} for 
either $I_0$ or $I_1$. The
failure of this argument is that no orthogonality with respect to $l$
is exploited. However, we remark that the bounds for $G_2$ have an extra factor
of $2^{k_3-k_1}$, which is more than enough to prove
\eqref{ioa} for $I_2$.

To remedy the above difficulty for $I_0$ and $I_1$ we separate the
nonlinear expression into two parts: one where $r$ is small, which we
can estimate directly, and the other one where $r$ is large, for which
we apply the above computation.  Given a threshold $m > -k_3^-$ we
split $I_0 = I_0^m + J_0^m$ where
\[
 I_0^m = \int_0^\infty \chi_{\geq m}(r) \psi^1_{k_1}(r) \psi^2_{k_2}(r) 
\psi^3_{k_3}(r) \psi_{j}(r) rdr dt
\]
The contribution $J_0^m$ of the region  $A_{<m}$ is estimated directly
via the local energy for the first two factors combined with the pointwise bound
\eqref{ps0} derived from the energy for the last two factors:
\begin{equation}\label{rlm}
|J_0^m| \lesssim  \| \chi_{<m}(r) \psi_{k_1}^1 {\psi}_{k_2}^2 
\psi_{k_3}^3 {\psi}_j \|_{L^1} 
\lesssim \frac{2^{j}}{|j|} \frac{ \la m \ra 2^{(m+k_3)/2}}{2^{k_1}}
\end{equation}
For $I_0^m$ we proceed as in the derivation of \eqref{nogood} 
but with $G_0$ replaced by its corresponding truncated version $G_0^m$ which
by Proposition~\ref{p:g} satisfies a better bound than $G_0$, namely
\[
|G_0^m| \lesssim \frac{2^{j/2}}{|j|} \frac{\la k_3^-\ra}{2^{k_3/2}}
2^{-N(m+k_3)}
\]
This leads to a similar improvement over \eqref{nogood}, namely
\begin{equation}\label{nowgood}
|I_0^m| \lesssim \frac{2^{j}}{|j|} \frac{\la k_3^-\ra}{2^{k_3}}2^{-N(m+k_3)}
\end{equation}
Adding \eqref{rlm} and \eqref{nowgood} gives
\[
|I_0| \lesssim \frac{2^{j}}{|j|}\frac{\la m \ra}{2^{k_1}} 
\left(  2^{(m+k_3)/2}
+ 2^{k_1-k_3}
2^{-N(m+k_3)}\right)
\]
Optimizing with respect to $m$ we obtain \eqref{ioa} for $I_0$.

We still need to consider $I_1$.  The bound for $I_1^m$ is identical 
to the one for $I_0^m$ since
 $\eqref{gg1R}$ gives the same bounds for $G_0^m$ and
$G_1^m$. For $J_1^m$ we apply a similar argument but some extra care 
is required due to the presence of the $[r \partial_r]^{-1}$ operator.
Precisely, using the $L^4_{k_1}$ and $LE_{k_2}$ norms for the two factors
we obtain
\[
2^{k_1} \| \psi_{k_1}^1 \psi_{k_2}^2\|_{L^{\frac43}(A_{<-k_1})}
+ 2^{\frac{5k_1-3m}8} \sum_{m \geq k_1}\| \psi_{k_1}^1 \psi_{k_2}^2\|_{L^\frac43(A_{m})}
\lesssim 1
\]
and the norm in the LHS above is preserved by the operator 
$[r \partial_r]^{-1}$ by \eqref{rdrmw}. Then using the pointwise bound derived from the energy
norm for $\psi_j$ respectively the $L^4_{k_3}$ norm for $\psi_{k_3}^3$
we obtain
\[
|J_1^m| \lesssim  \frac{2^{j}}{|j|} \frac{\la m \ra 2^{(m+k_3)/4}}
{2^{(5k_1+3k_3)/8}} 
\]
Though slightly weaker than \eqref{rlm}, this still leads to
\eqref{ioa} for $I_1$ when combined with \eqref{nowgood} for
$I_1^m$. The proof of the proposition is complete.

\end{proof}

\section{The bootstrap estimate for the 
$\l$ parameter.}
\label{seclambda}

In this section we show that the $Z_0$ regularity of the parameter
$\lambda$ can be bootstrapped, completing the proof of
Theorem~\ref{thflow}. The result in \eqref{rebootlambda} is 
obtained by replacing $\gamma$ with $\gamma + \gamma_0^2$ in the the
following Proposition.

\begin{p1}
Let $T \in (0,\infty]$. Consider $\tilde \lambda$ which satisfies the bound
\begin{equation}
\|\tilde \lambda-1\|_{\Zu[0,T]} \ll 1.
\end{equation}
and a  function $\psi$ satisfying
\begin{equation} \label{psibb}
\| \psi\|_{\WSs[\tilde \lambda][0,T]} \leq \gamma \ll 1 
\end{equation}
Then the functions $\lambda$, $ \alpha$ defined by
\begin{equation} \label{ltdef}
\psi_{2}(t,1) = i e^{i \alpha(t)} h_1^{\lambda(t)}(1)
\end{equation}
satisfy the bounds
\begin{equation}
\| \lambda -1 \|_{Z_0[0,T]} + \| \alpha \|_{Z_0[0,T]} \lesssim \gamma
\label{lzb}\end{equation}

\end{p1}

\begin{proof}
 From the fixed time analysis in Proposition~\ref{constr} we have 
the uniform bound in $[0,T]$
\[
\|\psi_2(1)-1\|_{L^\infty_t L^\infty_r[\frac12,2]} \lesssim \gamma
\]
which leads to the uniform bound
\[
 \| \lambda -1\|_{L^\infty} + \| \alpha\|_{L^\infty}
\lesssim \gamma
\]
To continue we recall the equation for $\psi_2$ from \eqref{comp1}: 
\[
\partial_r \psi_2  =  i A_2 \psi  - \frac{A_2 \psi_2}{r} 
\]

Approximating the second $A_2$ with $h_3^{\lambda}$
we rewrite this in the form
\[
 L_{\lambda} \psi_2 = i A_2 \psi  - 
\frac{(A_2 - h_3^{\lambda}) \psi_2}{r} := i \psi + N
\]
A solution to the homogeneous equation is $h_1^{\lambda}$.
Then by the same reasoning as in  Proposition~\ref{constr}
we must have 
\[
\psi_2(r)= i \lambda h_1^{\lambda}(r) +
 L_{\lambda}^{-1}(i \psi + N)
\]
where the factor $ \lambda$ in the first term on the right 
is dictated by the requirement that $\psi_2-i h_1 \in X$
which plays the role of the boundary condition at infinity.
Also $ L_{ \lambda}^{-1}$ is the $LX \to X$ inverse of 
$ L_{\lambda}$.  For functions $f \in L^1(rdr)$ we have the integral
formula
\[
 L_{\lambda}^{-1} f(r) = - h_1^{\lambda}(r)\int_{r}^\infty
\frac{f(s)}{ h_1^{ \lambda}(s)} ds
\]
We will be able to apply this formula for $N$ above, but not for $i\psi$.
However, in the case of $i\psi$ there is another computation that allows us 
to replace $L_{\lambda}^{-1}$ by  $L^{-1}$. Precisely,
\[
 L_{\lambda}^{-1} f(r) =  \frac{\lambda h_1^{\lambda}(r)}{
  h_1(r)}   L^{-1} f(r) -
 \lambda h_1^{ \lambda}(r) \int_{r}^\infty 
\left( \frac{1}{ \lambda h_1^{ \lambda}(s)} 
-  \frac{1}{h_1(s)}\right)  f(s) ds
\] 
which holds for all $f \in LX$. The integral converges even for all 
$f \in L^2$.

Then we rewrite $\psi_2$ in the form 
\[
 \psi_2(r)= i \lambda h_1^{ \lambda}(r) 
(1+(h_1(r))^{-1}A(r)-B(r)) 
\]
where
\[
A(r) =  L^{-1} \psi, \qquad B(r)= \int_{r}^\infty \psi(s)
\left( \frac1{ \lambda h_1^{\lambda}(s)}
  -\frac{1}{ h_1(s)} \right)-i \frac{N(s)}{
  \lambda h_1^{ \lambda}(s)}ds
\]
Set $r=1$ and recall that $\psi_2(1) = i e^{i\theta} h_1^{ \lambda}(1)$. Then
we obtain
\[
 e^{i\theta} = \lambda(1-A(1)-B(1))
\]
Since $Z_0$ is an algebra, it suffices to estimate $A(1)$ and $B(1)$
in $Z_0$.  For $A(1)$ we will establish a linear $Z$ bound,
\begin{equation}
\| A(1) \|_{Z[0,T]} \lesssim \gamma
\label{Abound}\end{equation}
which is facilitated by the fact that $ \lambda$
does not appear in the expression for $A$. On the 
other hand $ \lambda$ does appear in the $B$ expression;
however, there this does not matter since 
for $B(1)$ we establish a stronger bound
\begin{equation}
\|B(1)\|_{L^2_t[0,T]} \lesssim \gamma
\label{Bbound}\end{equation}
using only the pointwise bound for $\lambda$.
Together \eqref{Abound} and \eqref{Bbound} imply 
\eqref{lzb}. We consider these two bounds separately:

\bigskip

{\bf 1. The estimate \eqref{Abound} for the linear term.}

By the local energy decay estimate for $\psi$ we can replace 
$A(1)$ with a local average 
\[
\tilde A = \int A(r) \chi( r) dr 
\]
where $\chi$ is a bump function supported near $1$ 
with the normalization 
\[
\int \chi(r) h_1(r) dr = 1
\]
The difference admits good $L^2$ and $L^\infty$ bounds,
\[
\| A(1) -  \tilde A\|_{L^\infty} 
\lesssim \|\psi \|_{L^\infty L^2},
\qquad
\| A(1) -  \tilde A\|_{L^2} 
\lesssim \|\psi \|_{LE}
\]
Hence it remains to show that 
\begin{equation}
\| \tilde A \|_{Z_0[0,T]} \lesssim \| \psi\|_{\WSs[\tl][0,T]}
\label{Aboundz}\end{equation}

Using the Fourier expansion for $\psi$ in the $\tilde H$ frame,
respectively for $ L^{-1} \psi$ in the $ H$ frame,
we can represent $L^{-1} \psi$ in the form
\[
L^{-1} \psi(r) = \int \xi^{-1} \phi_\xi (r) \hat\psi(\xi) d\xi 
\]
Then for $\tilde A$ we have
\[
\tilde A =  \int h(\xi)  \hat\psi(\xi) d\xi,
\qquad
h(\xi) =  \xi^{-1} \int  \phi_\xi (r) \chi(r) dr.
\]
Given the representation of $\phi_\xi$ in Section~\ref{spectral},
it follows that $h$ is a smooth function in $(0,\infty)$ 
which has symbol type regularity, rapid decay at infinity
and whose size near $\xi = 0$ is given by
\[
h(\xi) \approx \frac{\xi^{-\frac32}}{|\log \xi|},\qquad \xi \ll 1.
\]

The bound \eqref{psibb} for $\psi$ in the proposition is given in
terms of the $\lambda$ frame, which is inconvenient as it makes it
difficult to track the different modulations. Fortunately, we are able
to use the  bound  \eqref{wsembI} to transfer enough of the $\WSs[\tl]$ 
norm to the fixed frame. It remains to show that
\begin{equation}
\| \tilde A \|_{Z[0,T]} \lesssim \| \psi\|_{\WS^r[1][0,T]}
\label{tAboundz}\end{equation} 
At this stage we can replace $\psi$ by any admissible $\WS^r[1]$ 
extension to the real line, and show that the above bound 
holds globally in time.

We decompose $\psi$ in frequency with respect to the  fixed $\lambda=1$ frame,
\[
\psi = \sum \psi_k, \qquad \psi_k = P_k \psi
\]
Correspondingly $\tilde A = \sum A_k$, where
\[
A_k =  \la g_k, \psi_k \ra, \qquad  \hat g_k(\xi) = 
\tilde\chi_k(\xi) h(\xi) 
\]
By Proposition~\ref{p:bumpft}, the functions $g_k$ satisfy the pointwise
bounds
\[
|g_k(r)| \lesssim  \frac{\log(1+r^2)}{\la k^- \ra^2} (1+ 2^k r)^{-N}
2^{-N k^+}
\]
The contribution of high frequencies $k \geq 0$ is easily estimated in
$L^\infty \cap L^2$,
\[
\|A_k\|_{L^\infty} \lesssim 2^{-Nk} \|\psi_k\|_{L^\infty L^2},
\qquad 
\|A_k\|_{L^2} \lesssim 2^{-Nk} \|\psi_k\|_{LE_k}
\]
It suffices to show that for the low frequencies we have 
\begin{equation} \label{akpsik}
\| A_k\|_{Z} \lesssim  \frac{2^{-k}}{|k|} \|\psi_k\|_{S^r_k}, \qquad k < 0
\end{equation}

On one hand we can use local energy decay to obtain an $L^2$ bound,
\[
\| A_k \|_{L^2} \lesssim \frac{2^{-2k}}{|k|} \| \psi_k\|_{LE_k}
\]
which suffices at low modulations,
\[
\|Q_{\leq 2k}  A_k\|_{Z} \lesssim \|Q_{\leq 2k}  A_\mu\|_{\dot H^\frac12}
\lesssim 2^{k} \| A_k \|_{L^2} 
\lesssim   \frac{2^{-k}}{|k|} \| \psi_k\|_{LE_k}.
\]

On the other hand for high modulations we can use the high modulation 
component of the $S^r_k$ norm:
\[
\begin{split}
\|Q_{> 2k}  A_k\|_{Z} \lesssim  \| Q_{> 2k}  A_k\|_{Z}
\lesssim \| g_k\|_{L^2} \| Q_{> 2k} \psi_k\|_{ ZL^2} 
\lesssim   \frac{2^{-k}}{|k|} \|\psi_k\|_{S^r_k}
\end{split}
\]
This concludes the proof of \eqref{akpsik} and thus the 
estimate for $A(1)$.

\bigskip

{\bf 2. The estimate \eqref{Bbound} for the nonlinear term.}

The analysis in this case is identical whether we work in a compact
interval or on the real line.  It suffices to place the integrand in
the formula of $B(r)$ in $L^2_t L^1_r[1,\infty)$.  For the first term
we note that
\[
\left|\frac1{\lambda h_1^{\lambda}} -\frac{1}{h_1} \right|
\lesssim r^{-1}.
\]
Then we can use the local energy bound \eqref{betterle}:
\[
 \|r^{-1} \psi\|_{L^2_t L^1_r([1,\infty];dr)}
= \|r^{-2} \psi\|_{L^2_t L^1_r ( A_{>0})} \lesssim \| \psi \|_{LE} \lesssim \|\psi\|_{\WSs[\lambda]} \lesssim
\gamma
\]
For the second term in $B$ we need to estimate
\[
\| r N\|_{L^2_t L^1_r([1,\infty];dr)} = \| N\|_{L^2_t L^1(A_{>0})}
\]
where
\[
 N = (A_2-1)\psi + \frac{1}r (A_2 - h_3^{\tilde \lambda}) \psi_2
\]
Since for $r > 1$ we have $|A_2-1| \lesssim |\psi_2|^2$, for the first term 
we use the local energy decay \eqref{betterle}
and the uniform in time bounds \eqref{linX} and \eqref{pointX} to write
\[
 \| (A_2-1)\psi\|_{L^2 L^1} \lesssim
\|  \la r \ra^{-\frac12-\delta} \psi \|_{L^2}
\|  \la r \ra^{-\delta} \psi_2 \|_{L^\infty L^2} \|  \la r
\ra^{2\delta}  \psi_2 \|_{L^\infty} \lesssim \gamma
\]
for some $0 < \delta < \frac14$.

Finally, for the second term in $N$ we use the local energy decay
\eqref{dap-eloc} and the uniform in time bound \eqref{linX}: 
\[
\|\la r\ra^{-1} (A_2 - h_3^{\tilde \lambda}) \psi_2\|_{L^2 L^1} \lesssim
\| \la r \ra^{-\frac12-\delta} (A_2 - h_3^{\tilde \lambda})\|_{L^2} \| \la r
\ra^{-\frac12+\delta}  \psi_2\|_{L^\infty L^2} \lesssim 
\gamma
\]
\end{proof}

\begin{r1}
  It is likely that $\lambda$ actually belongs to $Z$. However this is
  not needed in the present paper, and proving it would require a much
  more involved analysis of the nonlinear contribution $B$ above.
\end{r1}

\section{The bootstrap argument} \label{secboot}

In this section we prove Theorem~\ref{tmain-X}, i.e. our main global
well-posedness result. It is convenient to state the result in a more precise form:

\begin{t1}
a) Let $m = 1$ and $\e \leq 1, \gamma \ll 1$. Then for each 
$1$-equivariant initial data $u_0$ satisfying 
\begin{equation}
\| u_0 - \mathcal Q^1 \|_{\dot H^1} \leq \e \g , \qquad \| \bar{u}_0 - \bar Q\|_X \leq \gamma 
\label{stdata}\end{equation}
there exists a unique global solution $u$ so that
$\bar{u} - \bar Q \in C(\R;X)$ and 
\begin{equation}
 \| \bar u - \bar Q\|_{C(\R;X)} \lesssim \gamma
\label{stsolution}\end{equation}
Furthermore, this solution has a Lipschitz dependence
on the initial data in $X$, uniformly on compact time intervals.

b) Let $\psi$ be the reduced field associated to the solution $u$
above and $(\alpha,\l)$ defined by \eqref{defal}.
Then the following estimates 
are valid:
\begin{equation}\label{tbig:ld}
\|\psi \|_{\ldSs} \lesssim \e \g , \qquad \| W \psi \|_{\ldNs} \lesssim  \e \g^2,
\qquad \|\l-1\|_{Z_0} + \| \alpha\|_{Z_0} \lesssim \g 
\end{equation}
In addition, for any function $\tilde \l$  satisfying \eqref{tlambda}
we have
\begin{equation}\label{tbig:w}
 \| \psi \|_{\WSs[\tl]} \lesssim \g, \qquad
 \| W_{\tl} \psi \|_{\WNs[\tilde \l]} \lesssim \e \la \log \e \ra^2 \g^2 
\end{equation}

\label{tmain-XS}\end{t1}

We prove the theorem in two stages. First we use a
bootstrap argument establish global well-posedness and bounds for
regular solutions, i.e. solutions with initial data in $Q+H^2$. Then we
use a density argument to extend this result to initial data in $Q+X$.

\subsection{Regular solutions}

a) Given an initial data $u_0$ so that $u_0 - Q \in H^2$, by
Theorem~\ref{th2} we know that there exists a unique short time
solution $u$ so that $u(t) - Q \in C(0,T;H^2)$ for some small $T$.
In addition, we assume that $u_0$ satisfies \eqref{stdata}.

We will use a bootstrap argument to extend the time interval $T$ for
which the solution exists, with some suitable bounds. We begin by
describing the bounds we will bootstrap.  These are all expressed in terms 
of the reduced field $\psi$ and the soliton parameter $\lambda$.
We remark that in view of \eqref{stab}, \eqref{basicpsi}, 
and Theorem~\ref{thex}(a), the bounds \eqref{stdata} imply that 
the initial data $\psi(0)$ satisfies
\begin{equation}\label{psio}
\|\psi(0)\|_{L^2} \lesssim \epsilon \gamma, \qquad \|\psi(0)\|_{LX}
 \lesssim  \gamma
\end{equation}

We choose 
\[
\gamma < \gamma_0 = C \gamma  \ll 1
\]
where $C$ is chosen larger than the constants used in defining $\les$ in all
estimates in Theorem \ref{thflow}. 
Suppose we have a solution $u$ as
above in some interval $[0,T]$. With $\lambda$ defined as in
\eqref{analdef}, the first bootstrap bound will be
\begin{equation} \label{bootl}
\| \lambda -1 \|_{Z_0[0,T]} \leq \gamma_0
\end{equation}
The second bootstrap bound is concerned with the size of $\psi$ 
as an $L^2$ solution for a Schr\"odinger equation,
\begin{equation} \label{boot-psi-l2}
\| \psi \|_{l^2\Ss[0,T]}  \leq \epsilon  \gamma_0
\end{equation}
while the third bootstrap bound keeps track of the norm of $\psi$ in
$\WSs[\l]$ type spaces,
\begin{equation} \label{boot-psi-lx}
\| \psi \|_{\WSs[\tl][0,T]}  \leq \gamma_0
\end{equation}
for some  function $\tl$ which has the property that 
\begin{equation} \label{boottl}
\| \tl -1 \|_{\Zu[0,T]} \lesssim  \gamma_0, 
\qquad \| \tl -\l \|_{(L^2 \cap L^\infty)[0,T]} \lesssim  \gamma_0
\end{equation}

We denote by 
\[
\mathcal A = \{  \begin{array}[t]{l}T_0 \geq 0;\  \text{the Schr\"odinger map 
equation \eqref{SM} admits a solution $u \in C(0,T_0;H^2(\R^2))$} \\
\text{ so that its reduced field $\psi$ and soliton parameter $\lambda$  satisfy
 \eqref{bootl}, \eqref{boot-psi-l2}
and \eqref{boot-psi-lx}}\\ \text{ for all $T \leq T_0$}\}
\end{array}
\]
Our goal is to prove that $\mathcal A = [0,\infty)$. Once this is done, it follows
that we have a global solution satisfying \eqref{bootl}, \eqref{boot-psi-l2}
and \eqref{boot-psi-lx}. The bounds on $W\psi$  and $W_{\tl} \psi$ 
in \eqref{tbig:ld}, respectively \eqref{tbig:w} follow from Propositions~\ref{p:nonlin-l2}, \ref{p:nonlin-V},\ref{p:nonlin-x}. The estimate of \eqref{stsolution} follows 
from Theorem~\ref{thex}(b); this is where the qualitative property $u-Q \in L^2$ is 
is used in order to uniquely identify $u$ as the map associated to $\psi$ via 
Theorem~\ref{thex}(b).

By definition $\mathcal A$ is an interval containing $0$. Thus it suffices 
to prove the following two properties:
\begin{itemize}
\item[(i)] $\mathcal A$ is open in $[0,\infty)$.
\item[(ii)] $\mathcal A$ is closed in $[0,\infty)$.
\end{itemize}

\bigskip

{\bf (i)  $\mathcal A$ is open.}
Let $T_0 \in \mathcal A$. Then $u(T_0)-Q \in H^2$, therefore by
Theorem~\ref{th2} we have a local solution $u-Q \in
C([T_0,T_0+\delta_0];H^2)$.  From \eqref{bootl}, \eqref{boot-psi-l2} and
\eqref{boot-psi-lx} at $T_0$, applying Theorem~\ref{thflow}, we obtain
the  bounds 
\begin{equation}\label{better-boot}
 \| \lambda -1 \|_{Z_0[0,T_0]} \lesssim \gamma+ \gamma_0^2,
\ \   \| \psi \|_{l^2\Ss[0,T_0]}  \lesssim \epsilon(\gamma+ \gamma_0^2) ,
\ \  \| \psi \|_{\WSs[\tl][0,T_0]}  \lesssim \gamma+ \gamma_0^2
\end{equation}
which are stronger since $\gamma \ll \gamma_0 \ll 1$.
It remains to show that the above norms cannot grow much 
when replacing $T_0$ by $T_0+\delta$ with small $\delta$.
From the first norm we obtain $|\lambda(T_0)-1| \lesssim \gamma+\gamma_0^2$.
Since $\lambda$ is a continuous function of time, it follows 
that for small $\delta$ we have 
\[
\| 1_{[T_0,T_0+\delta]} (\lambda -1) \|_{L^2 \cap L^\infty} \lesssim 
\gamma + \gamma_0^2 
\]
This concludes the bootstrap for \eqref{bootl}.

Next we consider the second norm in \eqref{better-boot}.
For this we use the equation \eqref{wnlin-eq1} for $\psi$.
The linear part is well-posed in $L^2$, therefore it suffices 
to obtain a good bound for $W\psi$ in $[T_0,T_0+\delta]$,
\[
\| W\psi\|_{L^1([T_0,T_0+\delta]; L^2)} \lesssim \delta
\]
where, here and below, the implicit constant is allowed to depend on
the uniform $H^2$ bound for $u-Q$ in $[T_0,T_0+\delta_0]$. Indeed,
from this and the result of Theorem \ref{thflow}, we obtain
\[
\| \psi \|_{\ldSs [0,T_0+\delta]} \lesssim \| \psi(0) \|_{L^2} + \| 1_{[0,T_0]} W \psi \|_{\ldNs} +  \| 1_{[T_0,T_0+\delta]} W \psi \|_{\ldNs} 
\lesssim (\gamma +\g_0^2) \e + \delta 
\]
Similarly, for the third norm in \eqref{better-boot} we use the equation
\eqref{wnlin-eq} for $\psi$. The linear part is well-posed in $LX$,
therefore it suffices to obtain a good bound for $W_{\tl} \psi$ in
$[T_0,T_0+\delta]$,
\[
\| W_{\tl} \psi\|_{L^1([T_0,T_0+\delta]; LX)} \lesssim \delta
\]
with respect to a new\footnote{The exact choice of $\tl$ does not
  matter here; to fix things one could simply take $\tl = 1 + Q_{\leq
    1} (1_{[0,T+\delta]} (\lambda-1))$.}
function $\tl$ associated to $\l$ on the  interval $[0,T+\delta]$.
Indeed, in view of the result in Theorem \ref{thflow} and \eqref{linl}, 
\[
\begin{split}
\| \psi \|_{\WSs [\tl][0,T_0+\delta]} 
& \lesssim \|\psi_0 \|_{LX} +  \| 1_{[0,T_0]} W_{\tl} \psi \|_{\WNs[\tl]}
+ \| 1_{[T_0,T_0+\delta]} W_{\tl} \psi \|_{\WNs[\tl]} \\
& \lesssim (\g+\g_0^2)  + \delta 
\end{split}
\]
By choosing $\delta$ small enough, we can then use the result of part
b) in Theorem \ref{thflow} to bootstrap the bounds above and claim the
last two bounds in \eqref{better-boot} for $T_0+\delta$.

Using also the embedding \eqref{LXemb}, it suffices to show the fixed
time bound
\begin{equation}\label{wesr}
\|  W\psi\|_{L^2} \lesssim 1, \qquad \|  W_{\tl} \psi\|_{L^1 \cap L^2} \lesssim 1
\end{equation}
For this we use the $H^2$ regularity for $u-Q$ and its consequences in
Corollary~\ref{c:coulomb}. By Theorem~\ref{th2}, this regularity
persists up to time $T_0$ and also for a short time past $T_0$. We
consider each term in $W$ or $W_{\tl}$. For the cubic term we
simply use Sobolev embeddings to get $\|\psi\|_{L^2 \cap L^\infty}
\lesssim 1$. For the $\psi_2$ term we use the same, plus the $L^2$
bound for $\psi_2/r$. Finally, we split the $\delta A_2$ term into
two. For large $r$ we have the $r^{-2}$ decay factor so we only need
to use the pointwise boundedness of $\delta A_2$.  However, for small
$r$ we need to cancel that factor. This is easily done since for $r
\ll 1$ we have
\[
|\delta A_2| \lesssim r^2 + |\psi_2|^2 \lesssim  r^2 + |u_1|^2+|u_2|^2
\lesssim r^2
\]
where in the last step we have used Sobolev embeddings 
for $u_1$ and $u_2$, which vanish at the origin.
\bigskip

{\bf (ii) $\mathcal A$ is closed.}  Suppose that \eqref{bootl},
\eqref{boot-psi-l2} and \eqref{boot-psi-lx} hold for all $T <
T_0$. Then we have a Schr\"odinger map $u-Q \in C([0,T_0);H^2)$.
Passing to the limit in \eqref{bootl} we obtain \eqref{bootl}
for $T = T_0$. In particular this shows that $\lambda$ stays close 
to $1$ up to $T=T_0$. Then, by Theorem~\ref{th2}, it follows that the 
$H^2$ bounds persist up to (and beyond) $T=T_0$. Once we have 
$u-Q \in C([0,T_0];H^2)$ we repeat the above argument 
using  Theorem~\ref{thflow} in $[0,T-\delta]$ and then the bounds
\eqref{wesr} in $[T-\delta,T]$ with small $\delta > 0$ 
in order to prove \eqref{boot-psi-l2} and \eqref{boot-psi-lx}
for $T=T_0$.

b) The linear bounds in \eqref{tbig:ld} and \eqref{tbig:w} have been 
established above and the nonlinear bounds follow from \eqref{betterboot},
\eqref{betterbootV} and \eqref{betterbootx}.  

\subsection{ Rough solutions}

Given an initial data $u_0$ which satisfies 
\[
\| \bar u_0 - \bar Q\|_{X} \lesssim \gamma,
\]
we approximate it in the above topologies with a sequence of more
regular initial data $u_0^{(n)} \in Q + H^2$. Such approximations can
be obtained for instance by removing both the low and the high
frequencies in the $H$ frame,
\[
\bar u_0^{(n)} = \Pi_{\S^2} (\bar Q+ P^H_{[-n,n]} (\bar u_0- \bar Q)) 
\]
where $\Pi_{\S^2}$ represents the radial projection onto the sphere.
Since the projections $ P^H_{[-n,n]} (\bar u_0-\bar Q)$ stay pointwise small,
the convergence of $\bar u_0^{(n)}$ to $\bar u_0$ in $\bar Q+X$ follows from the
algebra property of $X$ and the bound \eqref{XLX} in Proposition \ref{Xalg}.

According to the previous step in the proof, for the initial data 
$u_0^{(n)}$ we obtain global solutions $u^{(n)}$ with $\bar u^{(n)}-\bar Q \in C(\R;X)$,
so that the bounds \eqref{better-boot} hold uniformly for 
the corresponding functions $\lambda^{(n)}$ and $\psi^{(n)}$. In particular
we have a uniform bound
\begin{equation}\label{psin}
\| \psi^{(n)}\|_{L^\infty LX} \lesssim \gamma
\end{equation}
By the first part of Theorem~\ref{thex}, the $X$ convergence of $\bar
u_0^{(n)}$ to $\bar u_0$ implies that $ \psi_0^{(n)}$ converges to $
\psi_0$ in $LX$.  By the short time result in Theorem~\ref{t:localwp},
it follows that the sequence $\psi^{(n)}$ converges in $\WSs[1][0,1]$
to some solution $\psi$ to \eqref{wnlin-eq1} with initial data
$\psi_0$.  In view of the uniform bound \eqref{psin} we can reiterate and
obtain a global solution $\psi$ to \eqref{wnlin-eq1}, so that for all $T > 0$
\[
\psi^{(n)} \to \psi \qquad \text{in } \WSs[1][0,T]
\]
 Furthermore, $\psi$ satisfies \eqref{better-boot}
globally in time. We note however that above we do not obtain
uniform convergence with respect to $T$.

Finally, given $\psi$ we apply the second part  of Theorem~\ref{thex} to 
construct a global Schr\"odinger map $u$ so that 
\[
\bar u^{(n)}-\bar Q \to \bar u-\bar Q \qquad \text{in } L^\infty X
\]
The local in time Lipschitz dependence of the solution $\bar u$ in $\bar Q+L^\infty X$
on the initial data $\bar u_0 \in \bar Q+X$ follows also by iterated application
of Theorem~\ref{t:localwp}, with the transition back and forth 
between $u$ and $\psi$ done via Theorem~\ref{thex}.

\section{ The $\dot H^1$ instability result}
\label{instabil}

In this section we prove the instability result
in  Theorem~\ref{tmain-H}.
% \begin{t1}
% Let $\gamma > 0$ small enough and $\lambda_0$ so that $|\lambda_0-1| \leq
% \gamma$. 
% Then for each $\epsilon > 0$ there exists
% a Schr\"odinger map $u$ and a time $t_0$ large enough with the following properties:
% \begin{equation}
% \|u(t)-Q\|_{X} \approx \gamma, \qquad t \in \R
% \label{inst1}\end{equation}
% \begin{equation}
% \|u(0)-Q_{\lambda_0}\|_{\dot H^1} \leq  \epsilon
% \label{inst2}\end{equation}
% \begin{equation}
% \|u(t)-Q\|_{\dot H^1} \leq  |\log \epsilon|^{-1}, \quad |t| > t_0
% \label{inst3}\end{equation}
% \end{t1}
% In other words, we can find solutions which start arbitrarily $H^1$-close
% to a soliton, yet as $t \to \pm \infty$ they travel for a fixed distance
% along the soliton family.
Let $\epsilon, \gamma \ll 1$ and $\alpha_0,\lambda_0$ as in
\eqref{alinst}, i.e.  so that $|\alpha_0|+|\lambda_0-1| \approx
\gamma$. We interpret $\e$ as a frequency parameter, and choose the
initial data $u_0$ so that $\bar u_0$ takes values in the $(\vec
i,\vec k)$ plane and
\[
\bar u(r) = \left\{ \begin{array}{ll} \bar Q_{\alpha_0,\lambda_0}(r),
    & r \ll \epsilon^{-1} \cr\cr \bar Q(r), & r \gg
    \epsilon^{-1} \end{array} \right.
\]
with a smooth transition on the $\epsilon^{-1}$ scale between the two
regions, so that
\begin{equation} \label{utran} |r (r \partial_r)^\alpha (\bar u_0-\bar Q)|
  \lesssim_\alpha \g, \qquad r \approx \e^{-1} \ \mbox{and} \ \alpha \geq 1.
\end{equation}
Using this and the form of the energy from \eqref{energy}, a direct 
computation shows that the bound \eqref{tdatai} holds,
\[
\| u_0 -Q_{\alpha_0,\lambda_0}\|_{\dot H^1} \lesssim \epsilon \gamma.
\]

To study the evolution of $u$ we switch to the $\psi$ variable.  From
the above information we characterize $\psi$ at time $t=0$. The
construction of the Coulomb gauge associated to $u(0)$ is trivial in
this case.  Since $\bar u$ stays in the $(\vec i,\vec k)$ plane and
$\bar v(\infty) = \vec k$, from the form of the ODE \eqref{cgeq} it
follows that $\bar v$ stays in the same plane. This implies that 
\[
\bar w(r) = \vec j, \qquad \bar v(r) = \bar u(r) \times \vec j
\]
 Recalling that
\[
\psi(0) = \partial_r u \cdot v + i \partial_r u \cdot w
\]
and using the above characterizations for $u,v,w$ we obtain the
following characterization for $\psi$:
\[
\psi(0)=0, \qquad r \ll \e^{-1} \ \mbox{and} \ r \gg \e^{-1}
\]
\[
|r^2(r \partial_r)^\alpha \psi| \lesssim_\alpha \g, \qquad r \approx \epsilon^{-1}
\]
In other words, $\psi(0)$ is a bump of size $\gamma \epsilon^2$
localized in the annulus $r \approx \epsilon^{-1}$. Hence it satisfies
the $L^2$ bound
\begin{equation}
\| \psi(0)\|_{L^2} \lesssim \gamma \epsilon 
\label{epsil2}\end{equation}

Using the above estimates on $\psi$ and the characterization of the
$\psi_\xi(r)$ from section \ref{spectralth}, it also follows that the
Fourier transform of $\psi$ satisfies
\begin{equation} \label{psi0est}
|(\xi \partial_\xi)^\alpha \mathcal{F}_{\tilde{H}} \psi(0,\xi)| 
\lesssim_{\alpha,N} \gamma \frac{\langle \ln \epsilon \rangle}{\langle \ln \xi \rangle}
 \xi^{\frac12} 
\langle \xi \epsilon^{-1} \rangle^{-N}, \qquad \alpha,N \in \N
\end{equation}
which directly leads to an $LX$ bound
\begin{equation}
\| \psi(0)\|_{LX} \lesssim \gamma
\label{epsix}\end{equation}
 This places us in the framework of the rest of the paper. Precisely, we obtain 
a global solution $u$ as in Theorem~\ref{tmain-XS},  which also satisfies 
the bounds \eqref{tbig:ld}, \eqref{tbig:w}.

By \eqref{alclose}, the desired bound \eqref{tsolutioni} would follow from the 
 estimate
\begin{equation}
|\psi_2(1,t)-i| \lesssim \g  |\log \epsilon|^{-1} \quad |t| > t_0
\end{equation}
We will in effect prove a stronger bound
\begin{equation} \label{stronger}
\|\psi_2(t)-ih_1\|_{\dot H^1} \lesssim \g |\log \epsilon|^{-1} \quad |t| > t_0
\end{equation}

By  Theorem \ref{tmain-XS}, we propagate the  bounds \eqref{epsil2} and
\eqref{epsix} along the flow:
\begin{equation} \label{phiflow}
\| \psi(t)\|_{L^2} \lesssim \gamma \epsilon, \qquad \| \psi(t)\|_{LX} \lesssim \gamma
\end{equation}
By \eqref{tbig:w} we have a
good\footnote{We actually get a stronger $\gamma \epsilon |\log
  \epsilon|^2$ bound here.}  bound for the nonlinearity in the $\psi$
equation,
\begin{equation}
\| (i \partial_t - \tilde H_{\tl}) \psi\|_{\WSs[\tl]} \lesssim
\g |\log \epsilon |^{-1}
\end{equation}
where $\tl$ is as in \eqref{tlambda}.
This shows that we can approximate $\psi$ in $LX$ by the 
solution to the corresponding linear equation, with $\gamma |\log \epsilon |^{-1}$
errors. Then, by Proposition~\ref{p:comp}, we can  compare
solutions to the linear $\tilde H$ equation with solutions to the 
 linear $\tilde H_{\tl}$ equation, 
\begin{equation}
\| \psi(t) - e^{it\tilde H} \psi(0)\|_{LX} \lesssim \gamma |\log \epsilon|^{-1}
\end{equation}
Thus by Proposition \ref{constr} it suffices to look at $\tpsi(t) =
e^{it\tilde H}\psi(0)$ and show that the corresponding $\tpsi_2$
associated to it satisfies \eqref{stronger}.

From \eqref{compsolh1} it follows that
\begin{equation} \label{auxpsib}
\|\tilde{\psi}_2(t) - i h_1\|_{\dot{H}^1} \lesssim \| L^{-1} \tilde{\psi}(t) \|_{\dot{H}^1}
\end{equation}
Denoting  $g(t)=L^{-1} \tilde \psi(t)$, we will prove that
\begin{equation}
\| g(t)\|_{\dot H^1} \lesssim \g (\epsilon +  \frac{|\log \epsilon|}{|\log t|}).
\label{winh1}\end{equation}
Together with \eqref{auxpsib} this establishes \eqref{stronger} and
concludes the proof of the theorem. Therefore we are left with proving
\eqref{winh1}.

The function $g$ has the Fourier expansion
\[
g(r,t) = \int_0^\infty \xi^{-1} \mathcal{F}_{\tilde{H}}
\tilde{\psi}(\xi,t) \phi_\xi(r) d \xi = \int_0^\infty \xi^{-1} e^{it
  \xi^2} \mathcal{F}_{\tilde{H}} \psi_0(\xi) \phi_\xi(r) d \xi
\]
where $\psi_0=\psi(0)$. We denote by $g_k, \tpsi_k$ the dyadic pieces
of $g$ respectively $\tpsi$ in the $H$, respectively $\tilde{H}$
calculus. We have the following

\begin{l1}
Let $q$ be as in \eqref{repphi}. Then we have 
\begin{equation} \label{gg1}
\|g_k\|_{\dot H^1} \lesssim \| \tpsi_k(t)\|_{L^2} +  \left | \int   \xi^{-1} q(\xi) \mathcal{F}_{\tilde{H}}
\psi_k(\xi,t)  d\xi\right|
\end{equation}
\end{l1}
\begin{proof}
 
For $g_k$ we have the straightforward $L^2$ relations
\[
 \| g_k(t)\|_{L^2} \lesssim 2^{-k} \|\tilde{\psi}_k(t)\|_{L^2}, \qquad 
\| (\partial_r + \frac{h_3}{r}) g_k(t) \|_{L^2} = \|\tilde{\psi}_k(t)\|_{L^2}
\]
Combining them we obtain
\[
 \| \chi_{r \gtrsim  2^{-k}} g_k(t)\|_{\dot H^1}  \lesssim \|\tilde{\psi}_k(t)\|_{L^2}
\]

It remains to estimate the part of $g_k$ in the region $\{r \lesssim
2^{-k}\}$. We consider the more difficult case $k < 0$. A similar but
simpler argument applies in the case $k \geq 0$. In the above region
we use \eqref{repphi} to obtain:
\[
g_k(r,t) = \int_0^\infty \xi^{-1} q(\xi) \mathcal{F}_{\tilde{H}}
\tilde{\psi}_k(\xi,t)  (h_1 + O( \xi^2
r \log r)) d \xi
\]
 respectively
\[
 \partial_r g_k(r) = \int_0^\infty \xi^{-1} q(\xi) \mathcal{F}_{\tilde{H}} \tilde{\psi}_k(\xi,t) 
 (h'_1 + O( \xi^2  \log r))
 d \xi
\]
The $O$ term in both formulas  admits the same bound as before.
The contribution of the principal part corresponds to the second 
term on the right-hand side of \eqref{gg1}. 
\end{proof}
Applying the above lemma  gives 
\[
\begin{split}
 \| g(t)\|_{\dot H^1} & \lesssim \sum_{k} \| \psi_k(0) \|_{L^2} + \sum_{k} 
 \left| \int   \xi^{-1} q(\xi) e^{it \xi^2} \mathcal{F}_{\tilde{H}} \psi_k(0,\xi)
 d\xi\right|
\end{split}
\]
The first term is easily estimated by $\gamma \e$ using \eqref{psi0est}.
For the integrals we use stationary phase together
with \eqref{qest} and \eqref{psi0est} to obtain 
\[
 \left| \int   \xi^{-1} e^{it \xi^2} \mathcal{F}_{\tilde{H}} \psi_k(0,\xi)
\frac{1}{\xi^{\frac12}\log \xi} d\xi\right | \lesssim_N 
\g \frac{|\ln \epsilon|}{|k|^2} \langle 2^k \e^{-1} \rangle^{-N} 
(1+ 2^{2k} |t|)^{-N} 
\]
Hence for large $t$ we obtain the bound
\[
\| g(t) \|_{\dot{H}^1} \lesssim \g \left(\e + \frac{\ln \epsilon}{|\ln t|}\right)
\]
which concludes the proof of \eqref{winh1}.

\end{document}